\documentclass[article, 11pt,reqno]{amsart}
\usepackage[english]{babel}
\usepackage[utf8]{inputenc}
\usepackage{amsmath, amssymb, amsthm}
\usepackage[table]{xcolor}
\usepackage{thmtools}
\usepackage{thm-restate}
\usepackage{graphics}
\usepackage[]{todonotes} 
\usepackage{physics}
\usepackage[margin=1in]{geometry}
\usepackage[hidelinks]{hyperref}
\usepackage{cleveref}
\usepackage[para]{footmisc}
\usepackage{caption}
\usepackage{subcaption}
\usepackage{floatrow}
\usepackage{xpatch}

\usepackage{graphicx}
\usepackage{mathrsfs}
\usepackage{verbatim}

\usepackage{tikz-cd}

\DeclareMathOperator*{\good}{good}
\DeclareMathOperator*{\is}{is}
\DeclareMathOperator*{\ext}{ext}
\DeclareMathOperator*{\int1}{int}
\DeclareMathOperator*{\Past}{Past}
\newcommand{\UEG}{\text{UEG}}
\DeclareMathOperator*{\rad}{diam}
\DeclareMathOperator*{\UEGop}{UEG}
\DeclareMathOperator{\perc}{perc}
\DeclareMathOperator{\arcsinh}{arcsinh}
\DeclareMathOperator{\arctanh}{arctanh}
\DeclareMathOperator{\TV}{TV}
\DeclareMathOperator{\FKG}{FKG}
\DeclareMathOperator{\clust}{clust}

\newcommand{\R}{\mathbb{R}}

\newcommand{\N}{\mathbb{N}}
\newcommand{\E}{\mathbb{E}}

\newcommand{\Z}{\mathbb{Z}}
\newcommand{\cc}{\leftrightarrow}
\newcommand{\Prb}{\mathbb{P}}
\newcommand{\Prbcur}{\mathbf{P}}

\renewcommand{\P}{\mathbb P}
\newcommand{\Even}{\Omega_{\emptyset}}
\newcommand{\Wired}{\Omega_{\emptyset}^1}

\usepackage{tikz}
\usepackage{pgfplots}
\usepackage{tikz-3dplot} 
\usetikzlibrary{positioning}
\usetikzlibrary{arrows}

\newcommand{\id}{1\! \!1}

\newcommand{\nn}{\mathbf{n}}

\newtheorem{theorem}{Theorem}[section]
\newtheorem{definition}[theorem]{Definition}
\newtheorem{lemma}[theorem]{Lemma}

\newtheorem{remark}[theorem]{Remark}
\newtheorem{corollary}[theorem]{Corollary}

\newtheorem{conjecture}[theorem]{Conjecture}

\newtheorem{question}[theorem]{Question}
\newtheorem{proposition}[theorem]{Proposition}

\makeatletter
\def\subsection{\@startsection{subsection}{2}%
  \z@{.5\linespacing\@plus.7\linespacing}{.3\linespacing}%
  {\normalfont\bfseries}}
\makeatother

\newlength\tindent
\title{The Uniform Even Subgraph and Its Connection to Phase Transitions of  Graphical Representations of the Ising Model}

\date{\today}

\pgfplotsset{compat=1.18}

\begin{document}

\author{Ulrik Thinggaard Hansen}
\address{Ulrik Thinggaard Hansen \\ Department of Mathematics
Université de Fribourg, Chem. de Musée 9, 1700 Fribourg, Switzerland }
\email{ulrik.hansen@unifr.ch}

\author{Boris Kj\ae r}
\address{Boris Kj\ae r \\ QMATH, Department of Mathematical Sciences, University of Copenhagen, Universitetsparken 5, 2100 Copenhagen, Denmark}
\email{bbk@math.ku.dk}

\author{Frederik Ravn Klausen}
\address{Frederik Ravn Klausen \\ QMATH, Department of Mathematical Sciences, University of Copenhagen, Universitetsparken 5, 2100 Copenhagen, Denmark}
\email{klausen@math.ku.dk}


\maketitle 
\begin{abstract}
The uniform even subgraph is intimately related to the Ising model, the random-cluster model, the random current model, and the loop $\mathrm{O}$(1) model. In this paper,  we first prove that the uniform even subgraph of $\Z^d$ percolates for $d \geq 2$ using its characterisation as the Haar measure on the group of even graphs. We then tighten the result by showing that the loop $\mathrm{O}$(1) model on $\Z^d$ percolates for $d \geq 2$ for edge-weights $x$ lying in some interval $(1-\varepsilon,1]$. Finally, our main theorem is that the loop $\mathrm{O}$(1) model and random current models corresponding to a supercritical Ising model are always at least critical, in the sense that their two-point correlation functions decay at most polynomially and the expected cluster sizes are infinite.
\end{abstract}

\section{Introduction}

The Ising model has been extensively studied for the past 100 years. A central tool in the recent study of the model has been its graphical representations, the random-cluster model \cite{fortuin1972random}, the random current model \cite{aizenman1982geometric, griffiths1970concavity}  and the high-temperature expansion \cite{van1941lange}, the latter also known as the loop $\mathrm{O}$(1) model. 

Even more recently, it has turned out that the uniform even subgraph of a graph is intimately related to the graphical representations of the Ising model. The uniform even subgraph of a finite graph $G$ is the uniform measure on the set of (spanning) subgraphs of $G$ with even degree at every vertex.  It holds that the loop $\mathrm{O}$(1) model can be sampled both as a uniform even subgraph of the random-cluster model \cite{GJ09} and as a uniform even subgraph of the (traced) double random current model \cite{klausen2021monotonicity}. 

These couplings may serve as motivation for the study of the uniform even subgraph and we will see that this perspective does, in fact, give new information about the percolative properties of the graphical representations of the Ising model. 
Meanwhile, the uniform even subgraph appears naturally on its own.
Indeed, we may consider the group of subgraphs of a given graph with symmetric difference of sets of edges as the group operation. The uniform even subgraph is then nothing but the Haar measure on the subgroup of even graphs.

In between the loop $\mathrm{O}(1)$ model and the random-cluster model, we find the (traced single) random current model, a graphical representation which has been central in the latest developments of the Ising model both in the planar case \cite{klausen2022mass}, in higher dimensions \cite{aizenman2021marginal,duminil2020exponential} and in even higher generality \cite{aizenman1987phase,duminil2016new,raoufi2020translation}. As with the random-cluster model, the random current model has recently become an object of inherent interest \cite{duminil2019double,duminil2021conformal2,duminil2021conformal}.

In \cite[Question 1]{DC16}, it was asked whether the single random current has a phase transition at the same point as the random-cluster model on $\Z^d$. 
From one point of view, since it is known that the random-cluster model can be obtained from the single random current by adding edges independently at random, a positive answer would imply that the added edges do not shift the phase transition. 
This, in turn, would run counter to common heuristics in statistical mechanics. 

 On the other hand, we will exploit a combinatorial fact about even subgraphs of the torus to prove that the correlations of the loop $\mathrm{O}(1)$ model, and hence the single random current, decay at most polynomially fast below the critical temperature of the random-cluster model. In models with positive association, this often signifies criticality. In particular, if the single random current were known to have a sharp phase transition, then our results would imply that the phase transition of the random current would coincide with that of the random-cluster model (and therefore, also the Ising model and the double random current). As such, this may be taken as evidence towards a positive answer to \cite[Question 1]{DC16}.

\subsection{Overview of the results and sketch of proofs}
In this paper, we let $\ell^{\xi}_{x,G},\Prbcur_{\beta,G}$ and $\phi^{\xi}_{p,G}$ denote the loop $\mathrm{O}$(1), the random current, and the random-cluster model, respectively. These are all models of random graphs, formal definitions of which will be provided in Section \ref{sec:preliminaries}.  For the time being, all we need to know about them is the following:

For a finite graph $G=(V,E),$ the uniform even subgraph of $G$ (henceforth $\UEG_G$) is the uniform probability measure on the set of spanning even subgraphs of $G$. If $\omega\sim \phi^{\xi}_{p,G}$ and $\eta$ is a uniform even subgraph of $\omega,$ then there exists $x=x(p)$ such that $\eta\sim \ell^{\xi}_{x,G}$. It is a classical result that the random-cluster model on subgraphs of $\mathbb{Z}^d$, for $d \geq 2$ undergoes a phase transition at a parameter $p_c\in (0,1),$ separating a regime of all clusters being small from one where there exists an infinite cluster \cite{peierls1936ising}. In the latter case, we say that the model \textit{percolates}. We denote $x_c=x(p_c)$. 

Our first contribution is to give an abstract characterisation of the uniform even subgraph of infinite graphs, an application of which is Theorem \ref{d3percolates} below. To fix notation, for a vertex $v$ of a random graph, we write $\mathcal{C}_v$ for the connected component of $v$, write $v\cc w$ for the event $\mathcal{C}_v=\mathcal{C}_w$ and write $v\cc \infty$ for the event that $|\mathcal{C}_v|=\infty$. 
Furthermore, it turns out that there is a natural infinite volume measure for the uniform even subgraph on $\mathbb{Z}^d,$ which we will denote by $\UEG_{\mathbb{Z}^d}$ (in the usual parlance, it is the wired infinite volume measure). See Section 3 for further explication.
 \begin{restatable}{theorem}{uegdpercolates}\label{d3percolates}
    For $d\geq 2$,
    $$
    \operatorname{UEG}_{\mathbb{Z}^d}[0\cc \infty]>0.
    $$
\end{restatable}
Previously, this was only known for $d=2$ \cite[Theorem 1.3]{GMM18}. Proving percolation of the uniform even subgraph was suggested to the authors as a toy problem, as a solution might shed light on the question about the phase transitions of random currents\footnote{UTH and FRK are grateful to Franco Severo and Aran Raoufi for posing the toy problem to them.}.
We solve the toy problem by showing a criterion for marginals of the uniform even subgraph to be Bernoulli distributed at parameter $\frac{1}{2}$ (see  \Cref{lemma:restriction-separation}).
The proof technique extends to showing that the phase transition of the loop $\mathrm{O}$(1) model is non-trivial for $\Z^d$, $d \geq 3$. 
In the spirit of the Grimmet-Janson coupling \cite[Theorem 3.5]{GJ09}, we define $\ell_{\mathbb{Z}^d}$ as the (wired) $\UEG$ of the associated random-cluster model (see Table \ref{table:parametrization_overview} for the coupling of the parameters).
\begin{restatable}{theorem}{nontrivphasetransition} \label{thm:non_tri_phase_transition}
Suppose that $d \geq 2$. 
Consider the loop $\mathrm{O}$(1) model $\ell_{x, \Z^d}$ with parameter $x \in \lbrack 0,1 \rbrack$. Then there exists an $x_0 < 1$ such that 
$$
\ell_{x,\Z^d}\lbrack 0 \cc \infty \rbrack > 0
$$
for all $x \in (x_0,1)$. 
\end{restatable}

Towards proving our main theorem, the main technical contribution is the insight that the loop $\mathrm{O}$(1) model is insensitive to boundary conditions. In the following, $\Lambda_n:=[-n,n]^d\cap \mathbb{Z}^d$ denotes the box of size $n$  around $0$ and $\mathbb{T}_n^d=\Lambda_n/(2n\mathbb{Z}^d)$ the associated torus.
Furthermore, we say a supergraph $G$ extends a graph $G_0$ with vertex boundary $\partial_vG_0$ if $G=G_0\cup G_1$ and the vertices of $G_1$ intersect $G_0$ only on its boundary.

\begin{restatable}{theorem}{mixing} \label{thm:Loop O(1) mixing}
    For $d\geq 2$ and $x\in(x_c,1]$, there exists $c>0$ such that for any $n\in \mathbb N$ and any event $A$ which only depends on edges in $\Lambda_n$, and $G$ is any supergraph extending $\Lambda_{4n}$,
    then
    \begin{equation*}
        |\ell^{\xi}_{x,G}[A]-\ell_{x,\mathbb{Z}^d}[A]| \le \exp(-cn),
    \end{equation*}
    for any boundary condition $\xi$. In particular, for $x\in (x_c,1]$  and any sequence $\xi_k$ of boundary conditions, $\lim_{k\to\infty}\ell^{\xi_k}_{x,\Lambda_{k}}=\ell_{x,\mathbb{Z}^d}$ in the sense of weak convergence of probability measures.
\end{restatable}

Due to the different behaviour of boundary conditions in the loop $\mathrm{O}(1)$ model (cf. \Cref{rem:boundary_conditions}), mixing does not follow immediately as would be expected from a model satisfying the Domain Markov Property \Cref{DMP}.
However, to complete the picture, mixing is proven in \Cref{theorem:the_real_mixing}. 

In the second step towards the main theorem, we exploit the topology of the torus to produce an essential lower bound.

\begin{restatable}{theorem}{lowerbound} \label{thm:Torus wrap}
Suppose $d\geq 2$ and let $x\in(x_c,1]$. Then, there exists $c>0$ such that
$
\ell_{x,\mathbb{T}_n^d}[0\cc \partial \Lambda_n]\geq \frac{c}{n}
$
for all $n$.
\end{restatable}

Combining  \Cref{thm:Loop O(1) mixing} and \Cref{thm:Torus wrap} proves our main theorem:
\begin{theorem}  \label{main theorem}
    Let $d \geq 2$ and $x\in (x_c,1].$ Then, there exists a $C>0$ such that for every $k$ and every $N\geq 4k$ and any boundary condition $\xi$, $\ell^{\xi}_{x,\Lambda_N}[0\cc \partial \Lambda_k]\geq \frac{C}{k}$. 
    It follows that, $\ell_{x,\mathbb{Z}^d}[\abs{\mathcal{C}_0}] = \infty$. 
\end{theorem}
From the couplings that we state in \Cref{thm:couplings}, the same result follows for the (sourceless, traced) single random current $\Prbcur_\beta$ for any $ \beta > \beta_c$, where $\beta_c$ is the Ising 
critical inverse temperature. 
\begin{corollary}
     Let $d \geq 2$ and $\beta>\beta_c.$ Then, there exists a $C>0$ such that $\Prbcur_{\beta,\Lambda_N}[0\cc \partial \Lambda_k]\geq \frac{C}{k},$ for every $k$ and every $N\geq 4k$.
    It follows that $\Prbcur_{\beta, \Z^d} [\abs{\mathcal{C}_0}] = \infty.$
\end{corollary}
Together with previous results (sharpness, couplings, $d=2$), the main theorem establishes the almost complete phase diagram for the loop $\mathrm{O}$(1) and (single) random current model on $\Z^d$ (see Figure \ref{fig:Phase diagram}).  We also argue the case of the hexagonal lattice $\mathbb{H}$ in Section \ref{sec:hexagonal}.

The overall strategy of our approach is simple. We consider the random-cluster model $\phi_p$ for $p >p_c$ on the $d$-dimensional torus of size $n$.  When $p$ is supercritical, one may prove the existence of a simple path wrapping around the torus once - henceforth called a wrap-around.
Whenever a wrap-around $\gamma$ exists in a random-cluster configuration $\omega$, if $\eta$ is a uniform even subgraph of $\omega,$ then $\gamma \triangle \eta$ also has the law of the uniform even subgraph of $\omega$. However, using the topology of the torus, we can prove that the number of wrap-arounds of the torus modulo 2 of the two configurations $ \eta, \gamma \triangle \eta$ are different. Therefore,  with probability at least $\frac{1}{2}$, there is at least one wrap-around in the uniform even subgraph  (and hence, in the loop $\mathrm{O}$(1) model and therefore also the random current model). Since $\eta$ has the distribution of the loop $\mathrm{O}$(1) model, this lets us provide lower bounds for connection probabilities for the loop $\mathrm{O}$(1) model on the torus.

The main technical part of the paper then lies in proving \Cref{thm:Loop O(1) mixing}, which allows us to transfer the connections of the torus into the space $\Z^d$. We, in turn, prove this theorem by a) proving that the uniform even subgraph is generically not very sensitive to boundary conditions and b) recalling some classical literature on the connectivity of the supercritical random-cluster model.

 This trick of exhibiting large clusters through the topology of the torus was previously employed to prove a polynomial lower bound for the escape probability of the Lorenz mirror model in \cite{kozma2013lower}. One may view the XOR-trick employed in \cite{crawford2020macroscopic} as another instance. However, the trick comes in many disguises. For example, it rears its head as the four-fold degeneracy in the ground state of the toric code \cite{kitaev2003fault}, which is an important candidate for the implementation of quantum error correction. 

 As a final aside, one may wonder just how sensitive the loop $\mathrm{O}(1)$ model is to the topology in which it is embedded. 
 The above trick makes explicit use of the topology of the torus, and as such, its application in proving \Cref{thm:Torus wrap} prevents us from generalising directly to a proof on $\Z^d$.
 On the other hand, \Cref{thm:Loop O(1) mixing} seems to tell us that the topology largely does not matter. In \Cref{sec:zipper}, we investigate this topological sensitivity and prove that removing certain edges in the hexagonal lattice can change the loop $\mathrm{O}$(1) model from not percolating to percolating.

\begin{figure}
\begin{center} 
\begin{tikzpicture}
\draw[dashed] (0,0)--(5,0);
\draw[line width=0.5mm] (5,0)--(10,0);
\node (X) at (-1,0) {$\Z^2$}; 
\draw node[fill,circle, scale=0.6]  at (10,0) {};
\draw node[fill, circle, scale=0.6]  at (5,0) {};

\node (X) at (-1, -2) {$\Z^d, d \geq 3$}; 
\draw[dashed] (0,-2)--(4,-2);
\draw[thin] (4,-2)--(8,-2);
\draw[line width=0.5mm] (8,-2)--(10,-2);
\draw node[fill,circle, scale=0.6]  at (10,-2) {};
\draw node[fill,circle, scale=0.6]  at (4,-2) {};
\draw node[fill,circle, scale=0.6]  at (8,-2) {};

\node (X) at (2,0.5) {exp decay}; 
\node (X) at (8,0.5) {percolation}; 
\node (X) at (2,-1.5) {exp decay}; 
\node (X) at (6,-1.5) {$\E[\abs{C(0)}] = \infty$ }; 
\node (X) at (9,-1.5) {percolation}; 
\node (X) at (5.2,-0.5) {$x_c^{\perc}$}; 
\node (X) at (8.2,-2.5) {$x_c^{\perc}$}; 
\node (X) at (4.2,-2.5) {$x_c^{\clust}$}; 

\node (X) at (0,0.5) {$x=0$}; 
\node (X) at (10,0.5) {$x=1$};

\node (X) at (-1, -4) {$\mathbb{H}$}; 
\draw[dashed] (0,-4)--(7,-4);
\draw[thin] (7,-4)--(10,-4);
\draw node[fill,circle, scale=0.6]  at (7,-4) {};

\node (X) at (2,-3.5) {exp decay}; 
\node (X) at (8.5,-3.5) {$\E[\abs{C(0)}] = \infty$}; 
\node (X) at (7.2,-4.5) {$x_c^{\clust}$}; 
\end{tikzpicture}
\end{center} 
\caption{The phase diagram of the loop $\mathrm{O}$(1) model and the single random current on $\Z^2$ and $\Z^d, d \geq 3$ as well as the loop $\mathrm{O}$(1) model on the hexagonal lattice $\mathbb{H}$,  all three in the ferromagnetic regime $x \in \lbrack 0,1 \rbrack$. What remains to be proven in the phase diagram is whether $x_c^{\perc} = x_c^{\clust}$ for $\mathbb Z^d, d \geq 3$ (as conjectured in Conjecture  \ref{conj:overall}). Note also that the case $x=1$ corresponds to the uniform even subgraph $\UEG$, which percolates for $\Z^d, d \geq 2$, but not for $\mathbb{H}$.  \label{fig:Phase diagram}}

\end{figure}

\tableofcontents

\section{Preliminaries: Graphical representations of the Ising model and their connections}
\label{sec:preliminaries}
In this section, we introduce the classical ferromagnetic Ising model $\mathbf{I}$ and its graphical representations: the random-cluster model $\phi$, the (sourceless, traced) random current model
$\Prbcur$ and the loop $\mathrm{O}$(1) model $\ell$. 

Furthermore, we will be concerned with the uniform even subgraph $\UEGop$, Bernoulli percolation $\Prb$ and the double random current $\Prbcur^{\otimes 2},$ the latter of which is obtained as the union of two independent copies of the single random current $\Prbcur$. 

The models are defined, first on finite graphs, then suitably extended to models on infinite graphs. For any graph $G=(V(G),E(G))=(V,E),$ we denote the space of \emph{percolation configurations}
$\Omega(G) = \{0,1\}^E,$ which is identified with $\mathcal{P}(E)$ under the map 
$\omega\mapsto E_\omega:=\omega^{-1}(\{1\})$. 
For $\omega\in\Omega(G)$, the associated spanning subgraph is $(V,E_\omega)$.
A measure on $\Omega=\Omega(G)$ will be called a percolation measure.
The Ising model $\mathbf{I}$ is a measure on $\{-1,1\}^V$, also called a spin model, while its graphical representations $\phi,\Prbcur,\ell$, as well as  $\Prb$, are percolation measures.
The models are related in a myriad ways, the most prominent of which are:\footnote{Details including parametrisations are given in \Cref{thm:couplings} and \Cref{table:parametrization_overview}.}
\begin{itemize}
    \item  The random current model $\Prbcur $ and the random-cluster model $\phi$ are obtained from the loop $\mathrm{O}$(1) model $\ell$ by adding additional edges independently at random (that is, according to $\Prb$). 
    \item The loop $\mathrm{O}$(1) model $\ell$ is obtained as a uniform even subgraph of either the random-cluster model $\phi$ or the double random current model $\Prbcur^{\otimes 2}.$ 
    \item For any $v,w \in V$,
     $$ \Prbcur^{\otimes 2}\lbrack v \cc w \rbrack = \langle \sigma_v \sigma_w \rangle^2 = \phi \lbrack v\cc w \rbrack^2.$$  
     These equalities mean that $\Prbcur^{\otimes 2}, \mathbf{I}$ and $\phi$ have a common phase transition at the critical inverse temperature $\beta_c$ of the Ising model.
     \item In two dimensions, the loop $\mathrm{O}$(1) model $\ell$ is the law of the interfaces of the Ising model.
\end{itemize}
    
We summarise the couplings in Figure \ref{fig:couplings} (partially borrowed from \cite{klausen2021monotonicity}).

\subsubsection{The Ising model} 
The much celebrated Ising model, introduced by Lenz \cite{lenz1920beitrag}, is a paradigmatic example of a model which undergoes a phase transition.  
The Ising model on a finite graph $G=(V,E)$ is a probability measure on the configuration space $\{-1,+1\}^V$. The energy of a configuration $\sigma\in\{-1,+1\}^V$ is
$$
H(\sigma) = - \sum_{(v,w)\in E} \sigma_v \sigma_w. 
$$
At parameter $\beta \in \lbrack 0, \infty \rbrack$, called the inverse temperature, the probability of a configuration is 
$$
\mathbf{I}_{\beta,G}[\sigma] = \frac{ \exp( - \beta H(\sigma))}{Z_{\beta,G}},
$$
where $Z_{\beta,G}=\sum_{\sigma} \exp(-\beta H(\sigma))$ is a normalisation constant called the partition function. This extends to the case $\beta=\infty$ by weak continuity. 

In keeping with the literature, we write $\langle \sigma_v \sigma_w \rangle_{\beta,G}$ for the correlation function $\mathbf{I}_{\beta,G}[\sigma_v\sigma_w],$ i.e. the expectation of the random variable $\sigma_v \sigma_w$ under the measure $\mathbf{I}_{\beta,G}$. 
\subsection{The graphical representations}
We start off by fixing some terminology. 
For a percolation configuration $\omega\in \Omega$, we say that $e$ is \emph{open} in $\omega$ if $\omega(e) = 1$ and $e$ is \emph{closed} if $\omega(e) = 0$. 
There is a canonical partial order on $\Omega$ given by pointwise comparison, that is, $\omega \preceq \omega'$ if $\omega(e) \leq \omega'(e)$ for all $e \in E$. 
We say that an event $\mathcal{A} \subset \Omega$ is \emph{increasing} if for all pairs $\omega, \omega' \in \Omega,$ it holds that if $\omega \in  \mathcal{A}$ and $\omega\preceq \omega',$ then $\omega' \in \mathcal{A}$. 
The notion of increasing events enables us to define a  partial order on percolation measures on $\Omega$ by declaring $\nu_1\preceq \nu_2$, if $\nu_1(\mathcal{A}) \leq \nu_2(\mathcal{A})$ for all increasing events $\mathcal{A}$. In this case, we say that $\nu_2$ \emph{stochastically dominates} $\nu_1.$ 

Further, for two percolation measures $\nu_1$ and $\nu_2$, we let $\nu_1 \cup \nu_2$ denote the measure sampled as the union of two independent copies of $\nu_1$ and $\nu_2$. 
That is, if $(\omega_1,\omega_2) \sim \nu_1\otimes \nu_2,$ then $\nu_1 \cup \nu_2$ is the law of  $\omega_1 \cup \omega_2$ which is defined by $E_{\omega_1 \cup\omega_2}=E_{\omega_1}\cup E_{\omega_2}$. 
Note that $\nu_1 \cup \nu_2\succeq\nu_1,\nu_2$.

\subsubsection{Boundary conditions}\label{rem:boundary_conditions}
If $G=(V,E)$ is a subgraph of the (typically infinite) graph $\mathbb{G}= (\mathbb{V,E})$, the edge boundary is $\partial_eG=\{(v,w)\in \mathbb{E}|\; v \in V, w\notin V\}$. Similarly, the vertex boundary is $\partial_vG=\{v\in V|\; \exists w:\; (v,w)\in \partial_eG\}$.

For a graph with boundary $\partial_vG\subset V$, a (topological) boundary condition is a partition of $\partial_vG$.
In physics parlance, the boundary vertices belonging to the same class in the partition are wired together.

There is a  partial order on boundary conditions given by fineness. If $\xi$ is a finer partition than $\xi'$, we write $\xi \preceq \xi'$.
There is a maximal and minimal boundary condition with respect to this order, namely the \emph{free} boundary condition where all vertices of  $\partial_v G$ belong to distinct classes (which is minimal), and the \emph{wired} boundary condition consisting of only one class (which is maximal). 
We denote these by $0$ and $1$ so that $0\preceq\xi\preceq 1$ for any boundary condition $\xi$. 
Given a boundary condition $\xi$, we may define the quotient multigraph $G/\sim_\xi$ by identifying vertices according to the partition. So for $\omega\in \Omega(G)$, we obtain $\omega^\xi\in\Omega(G/\sim_\xi)$. Define $\kappa^\xi(\omega)$ to be the number of connected components of $\omega^\xi$.
The connected components are often referred to as \emph{clusters}.

Typically, a percolation model is defined by a family of measures, say $\nu_G$, on the configurations on finite graphs $G$. 
Given a topological boundary condition $\xi$, define the measure $\nu^\xi_G$ via the natural identification between edges of $G$ and edges of $G/\sim_{\xi}$:
\begin{equation*}
    \nu^\xi_G[\omega]\propto \nu_{G/\sim_\xi}[\omega^\xi].
\end{equation*}
We refer to $\nu^{\xi}_G$ as the model with \emph{topological with boundary conditions} $\xi$. 
The Domain Markov Property (\Cref{DMP}) for the random-cluster model states that if $G\subset G'$ and $\psi$ denotes the restriction to $G\setminus G'$, then
\begin{equation*}
    \phi_{G'}[\omega | \psi] = \phi_{G}^{\xi_{\psi}},
\end{equation*}
where $\xi_{\psi}$ is defined such that $v$ and $w$ belong to the same class if and only if $v\overset{\psi}{\cc} w$. However, it is not generally true that these topological boundary conditions correspond to the usual DLR-type boundary conditions for models satisfying a Markov property. For instance, this is not true for the loop O($1$) model.

We show in \Cref{prop:marginal-is-boundary-condition} that boundary conditions arise in a rather different way for the uniform even subgraph (denoted $\UEG$, see \Cref{sec:loop O(1) model}). 
There is a natural marginal of $\UEG_{G'}$ on the configurations of $G$ satisfying
\begin{equation*}
    \UEG_{G'}|_{G} = \UEG_{G}^{G'}.
\end{equation*}
Here $G'$ is by abuse of notation the equivalence relation given by connection\footnote{
When $G$ is infinite we require a notion of $\infty$-connectivity, see \Cref{sec:marginals-of-UEG}.
} in $G'\setminus G$. Hence, marginals of the UEG behave analogously to conditional measures of the random-cluster model.

\subsubsection{The random-cluster model}
The most well-studied graphical representation of the Ising model is the random-cluster model. For $p\in \lbrack 0,1 \rbrack$, finite graph $G$, and boundary condition $\xi$, it is the probability measure given by
\begin{align*}
\phi^{\xi}_{p,G} \lbrack \omega \rbrack \propto 2^{\kappa(\omega^{\xi})} \left( \frac{p}{1-p} \right)^{o(\omega)}, 
\end{align*}
where $o(\omega)=|E_{\omega}|$ denotes the number of open edges in $\omega$. 
Whenever the boundary condition is omitted, we assume the free boundary condition, i.e. $\phi_{p,G} = \phi^{0}_{p,G}$.

The random-cluster model is related to the Ising model through the Edwards-Sokal coupling \cite{edwards1988generalization}: 
Suppose $\omega\sim \phi_{p,G}$ and $\sigma$ is obtained from $\omega$ by independently assigning $+$ and $-$ spins to the clusters of $\omega$. Then, $\sigma\sim \mathbf{I}_{\beta,G}$. Conversely, $\omega$ can be sampled from $\sigma$ by taking each edge $e = (v,w)$ such that $\sigma_v = \sigma_w$ and opening it with probability $p$. 
From the Edwards-Sokal coupling, it follows that  
\begin{align}\label{eq:random_cluster_representation}
\langle \sigma_v \sigma_w \rangle_{\beta, G} = \phi_{p,G} \lbrack v \cc w \rbrack
\end{align}
for $p = 1- e^{-2 \beta}$ (cf. \cite[Corollary 1.4]{DC17}). In particular, \eqref{eq:random_cluster_representation} implies that there is long-range order for the Ising model if and only if there are large clusters in the random-cluster model. Thus, the physical properties of the Ising model may be studied via the graphical properties of the random-cluster model. 
Towards this end, the random-cluster model has several monotonicity properties which make it amenable to analysis. 
While this paper is not focused on proving something new about the random-cluster model, these properties will nonetheless play a crucial role as we explore the other graphical representations through the couplings. Here, we repeat a tailored version of the more general statement in \cite{hansen2022strict}. For these and more results about the random-cluster model, we refer the reader to \cite{DC17} and \cite{Gri06}.

\begin{proposition} \label{thm:comparison theorem}
For the random-cluster model on a finite subgraph $G$ of a graph $\mathbb{G}$ (finite or infinite), the following relations hold:
\begin{enumerate}
\item[$i)$] The measure $\phi^{\xi}_{p,G}$ is monotone in $\xi$ in the sense that if  $\xi \preceq \xi'$, then
$$
\phi_{p, G}^{\xi}  \preceq   \phi_{p,G}^{\xi'}
$$
for any parameter $p\in \lbrack 0,1 \rbrack$.
\item[$ii)$] The measure $\phi^{\xi}_{p,G}$  is increasing in $p,$ i.e.\ if $p\leq p',$ then 
$$
\phi^{\xi}_{p,G}\preceq \phi^{\xi}_{p', G}
$$
for any boundary condition $\xi$.
\item[$iii)$] The random-cluster model is comparable to Bernoulli percolation in the following sense: 
$$
\mathbb{P}_{\tilde{p},G}\preceq \phi^{\xi}_{p,G} \preceq \mathbb{P}_{p,G}
$$
for any boundary condition $\xi$, where $\tilde{p}=\frac{p}{2-p}$. 
\end{enumerate}
\end{proposition}

\begin{proposition}[FKG inequality] \label{FKG}
Let $G$ be a finite graph, and  $A$ and $B$  increasing events, then
$$
\phi^{\xi}_{p,G}[A\cap B]\geq \phi^{\xi}_{p,G}[A]\phi^{\xi}_{p,G}[B]
$$
for any boundary condition $\xi$ and $p\in [0,1].$
\end{proposition}

\begin{proposition}[Domain Markov Property] \label{DMP}
If $G_1=(V_1,E_1)\subset G_2=(V_2,E_2)$ are two finite subgraphs of an infinite graph $\mathbb{G}$, we write $\omega_1:=\omega|_{E_1}$ and $\omega_2:=\omega \vert_{E_2\setminus E_1}$.
 Then, for any boundary condition $\xi$ on $G_2$ and any event $A$ depending on edges in $G_1$, it holds that
$$
\phi_{p,G_2}^{\xi}[\omega_1\in A|\;\omega_2]=\phi_{p,G_1}^{\xi_{\omega_2}}[A],
$$
where $v,w\in V_1$ belong to the same element of $\xi_{\omega_2}$, if and only if they are connected (or possibly equal) in $(V_2, E_{\omega_2})/\sim_\xi$.
\end{proposition}

\subsubsection{Construction of infinite volume measures} 

The monotonicity in boundary conditions combined with the Domain Markov Property allows us to define infinite volume measures in the following way:
For a sequence of finite graphs $G_n\Uparrow\mathbb{G} = (\mathbb{V}, \mathbb{E})$ for some infinite graph $ \mathbb{G},$ it holds that the marginals of $\phi_{G_n,p}^1$ on a fixed finite subset $\Lambda$ are monotonically decreasing. In other words, if $\mathcal{A}$ is an increasing event that depends only on edges in $\Lambda,$ then 
$
\{ \phi_{G_n,p}^1( \mathcal{A}) \}_{n\in \N}
$
is monotonically decreasing. 
Similarly, the sequence $
\{ \phi_{G_n,p}^0( \mathcal{A}) \}_{n\in \N}
$
is monotonically increasing. Since they are bounded, they have limits $\phi_{p,\mathbb{G}}^1( \mathcal{A})$ and $\phi_{p,\mathbb{G}}^0( \mathcal{A})$ respectively. 
One can check that the limit does not depend on the choice of sequence $G_n$. Since the set of increasing events is intersection-stable and generates the product $\sigma$-algebra of $\{0,1\}^{\mathbb{E}}$, the measures of such sets determine any possible limit measures $\phi^1_{p,\mathbb{G}}$ and $\phi^0_{p,\mathbb{G}}$ uniquely. Since $\{0,1\}^E$ is compact, so is the space of percolation measures by the Banach-Aloglu theorem, so the sequences $(\phi^1_{p,G_n})_{n\in \mathbb{N}}$  and $(\phi^0_{p,G_n})_{n\in \mathbb{N}}$ do have at least one accumulation point. By the uniqueness outlined above, both sequences are weakly convergent towards (possibly equal) limits $\phi^1_{p,\mathbb{G}}\succeq \phi^0_{p,\mathbb{G}}$.

On the other hand, any measure $\nu$ on  $\{0,1\}^{\mathbb{E}}$ which almost surely has the Domain Markov Property (cf. \Cref{DMP}) could be called an infinite volume random-cluster measure. However, such a measure would necessarily satisfy $\phi^0\preceq \nu \preceq \phi^1$ and by \cite[Corollary 3]{raoufi2020translation}, $\phi^0_{p,\mathbb{Z}^d}=\phi^1_{p,\mathbb{Z}^d}$ for all $p.$ Accordingly, there is a unique infinite volume measure, and we shall drop the boundary conditions from our notation in the infinite volume case and merely write $\phi_{p,\mathbb{Z}^d}$. 
A similar construction defines the infinite volume Ising measures $\mathbf{I}_{\beta,\mathbb{Z}^d}^{0},\mathbf{I}_{\beta,\mathbb{Z}^d}^{+}, \mathbf{I}_{\beta,\mathbb{Z}^d}^{-}$.
With the infinite volume measures at hand, the following is obtained from \eqref{eq:random_cluster_representation} in a more or less straightforward manner.
\begin{proposition}
Let $ \mathbb{G} = \mathbb{Z}^d$ and let $\beta \geq 0$ be given. Then there is long-range order in the Ising model if and only if there is an infinite cluster in the random-cluster model percolates almost surely.
That is,
\begin{align*}
\inf_{v,w\in \mathbb{Z}^d}\langle \sigma_v \sigma_w \rangle_{\beta,\mathbb{Z}^d} >0,
\end{align*}
if and only if
\begin{align} \label{eq:phi_Ising} 
\phi_{\beta,\mathbb{Z}^d}\lbrack 0 \cc \infty \rbrack > 0. 
\end{align}
\end{proposition}
On the hypercubic lattice $\Z^d$, for $d \geq 2$, there exists a unique sharp phase transition \cite{aizenman1987phase,peierls1936ising}. This means that there exists a $\beta_c$ such that $\langle \sigma_v \sigma_w \rangle $ decays exponentially in the distance between $v$ and $w$ for all $\beta < \beta_c$ and that there is long-range order for all $\beta > \beta_c$.

\subsubsection{Even subgraphs}
Both the random current model and loop $\mathrm{O}$(1) model are defined in terms of even subgraphs of a graph. 
A graph is said to be \emph{even} if every vertex degree is even.
We  denote the set of all percolation configurations corresponding to even graphs by $\Even$.

A related notion often used in the context of multigraphs in the random current literature is that of sources. For any (multi-)graph $H = (V_H, E_H) $ we say that the \emph{sources} of $H$, denoted $\partial H,$ is the set of vertices of odd degree. An even graph is then a graph such that $\partial H = \emptyset$. 
Just like we identified spanning subgraphs of a given graph $G = (V,E)$ with the space of percolation configurations $\Omega,$ so we identify a configuration $\nn \in \N_0^{E}$ with a multigraph, and such a multigraph is called a \emph{current}. If $\partial \nn = \emptyset,$ we say that the current is \emph{sourceless}.

\subsubsection{The random current model}
We now briefly introduce the random current model. For a more complete exposition, see \cite{DC16} or \cite{DC17}.
For a finite graph $G = (V,E)$, define the weight 
$$
w_\beta(\nn) = \prod_{e\in E} \frac{ \beta^{\nn_e}}{\nn_e!}.
$$
Then the random current with source set $A$, denoted $\Prbcur^{A}_{\beta,G}$,  is the probability measure on $\mathbb{N}_0^E$, given by 
\begin{align*}
\Prbcur^{A}_{\beta,G} \lbrack \nn  \rbrack \propto  w_\beta(\nn) \id_{\{\partial \nn = A\}}.
\end{align*}
Since we can view any deterministic boundary condition as a free boundary condition on an appropriate graph, we will mostly work with free random currents $\Prbcur^{A}_{\beta,G} = \Prbcur^{A,0}_{\beta,G}.$
The following relation  (cf. \cite[(4.5)]{DC17}) provides the first relation of the random current model to the Ising model
\begin{align*}
    \langle \sigma_v \sigma_w \rangle_{\beta,G}
    = \frac{ \sum_{\nn \mid \partial \nn = \{v,w\}} w_\beta(\nn) }{ \sum_{\nn \mid \partial \nn = \emptyset} w_\beta(\nn)}. 
\end{align*}

From the multigraph valued random current, one derives a percolation measure called the traced random current. For any current $\nn \in \N_0^E,$ the corresponding traced current $\hat \nn$ is defined by $\hat \nn(e) = \id_{\nn(e) > 0}$ for each $e \in E$. Of prime importance to the trace operation is the fact that it does not change connectivities, i.e. $v\cc w$ in $\nn$ if and only if $v \cc w$ in $\hat\nn$. In the following, we will only discuss the traced random current without sources, which we will denote $\Prbcur_{\beta,G}$.

Finally, a key player in the random-current literature is the \emph{(sourceless) double random current}, defined by $$
\Prbcur^{\otimes 2}_{\beta,G}:= \Prbcur_{\beta,G} \cup \Prbcur_{\beta,G}.
$$
 As a celebrated consequence of the switching lemma (see \cite[Lemma 4.3]{DC17}), one can prove that
\begin{align} \label{eq:dc_Ising} 
    \langle \sigma_v \sigma_w \rangle_{\beta,G}^2
    =\Prbcur^{\otimes 2}_{\beta,G} \lbrack v \cc w \rbrack . 
\end{align}
Thus, the double random current, just like the random-cluster model, has the onset of large clusters at the critical point of the Ising model. 
One of the main motivations for this paper is to investigate for which infinite graphs $\mathbb{G}$ this is also true of the single random current $\Prbcur_{\beta, \mathbb{G}}$. A positive answer would be implied by a positive answer for the loop $\mathrm{O}(1)$ model, which we shall now introduce.

\subsubsection{The loop \texorpdfstring{$\mathrm{O}$}{O}(1) model}\label{sec:loop O(1) model}
For a finite graph $G = (V,E),$ we now define the loop $\mathrm{O}$(1) model $\ell_{x,G}^{\xi}$ for parameter $x \in \lbrack 0, 1 \rbrack$ as a measure on $\Omega(G)$.
For $x \in \lbrack 0, 1 \rbrack,$ the loop $\mathrm{O}$(1) model $\ell_{x,G}^{\xi}[\eta]$ is defined by
$$
\ell_{x,G}^{\xi} \lbrack \eta \rbrack \propto x^{o(\eta)} \id_{\partial \eta^{\xi} = \emptyset}.
$$

In particular, for $x=1$, $\ell_{1,G}^0$ is the uniform measure on $\Omega_\emptyset(G)$, which we denote by $\UEG_G$. The expert reader may note that $\ell_{x,G}^0$ is the high temperature expansion of the Ising model for $x = \tanh(\beta)$. 

We may choose to view the sourceless random current $\Prbcur_{\beta,G}$ and the loop $\mathrm{O}$(1) model $\ell_{x,G}$ as respectively independent Poisson and a Bernoulli variable on each edge conditioned on the sum over the valences on edges adjacent to a given vertex being even. This point of view was used in e.g. \cite{tassion_notes}. 

For planar graphs, the loop $\mathrm{O}$(1) model is the law of the interfaces of a corresponding Ising model on the faces of the graph. This is discussed further in Section \ref{sec:planar}. 

\subsubsection{Bernoulli percolation} 
Although it is not a graphical representation of the Ising model, we also consider Bernoulli percolation with parameter $p \in \lbrack 0, 1 \rbrack,$ which we will denote by $\Prb_p$. This is the percolation measure where every edge $e \in E$ is open with probability $p$ independently. The model was introduced in \cite{broadbent1957percolation} and has since been subject to intense study \cite{duminil2018sixty}. Here it will mainly play an auxiliary role.

\subsubsection{Coupling the graphical representations}
To state the couplings, let us define the uniform even subgraph of a probability measure.

For a probability measure  $\mu$ on $\Omega$, consider the measure obtained by first sampling $\mu$ and then sampling the $\UEG$ of the first sample. 
That is, for every even graph $\eta \in \Even,$ consider $\UEG_\omega[\eta]$ as a function of $\omega\in \Omega$, so
\begin{align*}
\mu \lbrack  \UEG _{\omega}[\eta] \rbrack = \sum_{\omega \in \Omega}  \UEG _{\omega} \lbrack \eta \rbrack \mu \lbrack \omega \rbrack
= \sum_{\omega \in \Omega} 
    \frac{ \id_{ \eta \preceq \omega}}{ \abs{\Even \lbrack \omega \rbrack}} \mu \lbrack \omega \rbrack.   
\end{align*}

Now, we state the couplings from \cite[Exercise 36]{DC17}. The original references are  \cite[Theorem 3.5]{GJ09},  \cite[Theorem 4.1]{klausen2021monotonicity}, \cite[Theorem 3.1]{Lis} and \cite{LW16}. 
\begin{theorem} \label{thm:couplings} 
For any finite graph $G= (V,E),$ the graphical representations of the Ising model are related as follows:
\begin{itemize}
\item $\ell_{x,G}^0 \cup \mathbb{P}_{ 1 - \cosh(\beta)^{-1}, G} = \Prbcur_{\beta,G} $ 
\item $\ell_{x,G}^0 \cup \mathbb{P}_{\tanh(\beta),G} = \phi_{p,G}^0$
\item  $\Prbcur^{\otimes 2}_{\beta,G} \lbrack \UEGop_\omega \lbrack \cdot \rbrack \rbrack \overset{}{=} \ell_{x,G}^0\lbrack \cdot \rbrack \overset{}{=} \phi_{p,G}^0 \lbrack \UEGop_\omega \lbrack \cdot \rbrack \rbrack $
\end{itemize} 
Here, the parametrisations are related by $x=\tanh(\beta)$ and $p=1-e^{-2\beta}$, cf. \Cref{table:parametrization_overview}. 
\end{theorem}
\begin{remark}
We note that even though the couplings here are only formulated with free boundary conditions, the measures with boundary conditions correspond to the free measure on the quotient of the original graph. As such, the couplings in the theorem also hold for e.g. wired boundary conditions. 
\end{remark}

\begin{proof}
For a proof of the first two points, see \cite[Theorem A.1]{klausen2021monotonicity}. 
The first equality in the last point is \cite[Theorem 4.1]{klausen2021monotonicity} and the other equality is the statement of \cite[Theorem 3.5]{GJ09}. 
We give a proof in our notation. For any even graph $\eta,$ it holds that
\begin{align*}
\phi_{p,G}^0 \lbrack \UEG_\omega \lbrack \eta \rbrack \rbrack  \propto  \sum_{\omega \in \Omega}  \frac{2^{\kappa(\omega)}\id_{\eta \preceq \omega}}{\abs{\Omega_{\emptyset}(\omega)}}\hspace{-0.1 cm}\left( \frac{p}{1-p} \right)^{\hspace{-0.1 cm}o(\omega)}  \hspace{-0.2cm}
 \propto  \sum_{\omega \in \Omega}  \left( \frac{p}{2(1-p)} \right)^{ \hspace{-0.1 cm}o(\omega)} \hspace{-0.4cm} \id_{\eta \preceq \omega} \propto \Prb_{x,G}[\eta \text{ open}]= x^{o(\eta)} \propto \ell_{x,G}^0 \lbrack \eta \rbrack, 
\end{align*}
where $\Prb_{x,G}$ denotes the expectation under Bernoulli percolation with parameter $x = \tanh(\beta) = \frac{p}{2-p}$ and we used the fact that
$\abs{\Omega_{\emptyset}(\omega)} = 2^{\kappa(\omega) + \abs{o(\omega)}-\abs{V}}$.

\end{proof}

The relation $\abs{\Omega_{\emptyset}(\omega)}= 2^{\kappa(\omega) + \abs{o(\omega)}-\abs{V}}$ is classical and one of its proofs goes as follows: First, note that the relation is true if $\omega$ is a forest\footnote{That is, a graph without cycles.}. Now, if $e$ is part of a loop $l$, then $\eta \mapsto \eta \triangle l$ is a bijection between the set of even subgraphs with $e$ open and the set of even subgraphs with $e$ closed. Since the latter can be identified with $\Even((V,E\setminus\{e\})),$ it follows that adding an extra edge $e$ to a connected graph doubles the number of even subgraphs. 
This bijection will be important in many arguments throughout the paper.

One might also take the coupling in \Cref{thm:couplings} as the definition of the loop $\mathrm{O}$(1) model with boundary condition, so
\begin{align}\label{eq:loop_with_bc}
\ell_{x,G}^{\xi} \lbrack \cdot \rbrack = \phi_{p,G}^\xi \lbrack {\UEGop}_\omega \lbrack \cdot \rbrack \rbrack.
\end{align}
Notice that this definition is consistent with our previous definition of $\ell_{x,G}^{\xi}$. 

We can also take this approach to define the loop $\mathrm{O}(1)$ model in infinite volume, as was noted in  \cite[Remark 3.16]{angel2021uniform}. We define the uniform even subgraph on an infinite graph $\mathbb G$ as the Haar measure on the group of even graphs\footnote{Even in finite volume, this perspective is interesting. It is not generally true that uniform measures are nice, easy objects to study. For instance, the uniform measure on self-avoiding walks of fixed length is infamously difficult to study. Groups, however, naturally come equipped with a rich family of symmetries, which, as we shall see, aids a lot in analysis.}, which allows us to define the loop $\mathrm{O}(1)$ model on an infinite graph as
\begin{align} \label{eq:def_infinite_loop_O}
\ell_{x,\mathbb G} \lbrack \cdot \rbrack = \phi_{p,\mathbb G} \lbrack {\UEGop}_\omega \lbrack \cdot \rbrack \rbrack. 
\end{align}
The equivalence of this definition to those given in \cite{angel2021uniform, GJ09} is discussed in \Cref{section:limit-constructions}. 
There, we will also see that $\UEG_{\mathbb G}$ is unique as a Gibbs measure when $\mathbb G$ is one-ended.
This is almost surely the case for the infinite cluster of $\omega\sim \phi_{p,\mathbb{Z}^d}$ when $d\ge 2$ by the Burton-Keane theorem \cite[Theorem 2]{BurtonKeane} and deletion tolerance.

The literature \cite{aizenman2015random, hutchcroft2020continuity} has yet another construction of the infinite volume loop $\mathrm{O}$(1) model which uses a relation of the loop $\mathrm{O}$(1) model to the gradient Ising model. Fortunately, the uniqueness of the infinite volume measure proven in \Cref{thm:Loop O(1) mixing} implies that the constructions agree. We remark that \eqref{eq:def_infinite_loop_O} also gives an independent construction of the infinite volume random current measure as 
$$ \Prbcur_{\beta, \mathbb{Z}^d } = \ell_{x, \mathbb{Z}^d} \cup \mathbb{P}_{ 1 - \cosh(\beta)^{-1},  \mathbb{Z}^d.}$$

\begin{figure}
{
\begin{tikzcd}
\phi_p \arrow[rd, dashed] & \Prbcur_\beta \arrow[r] \arrow[l] & \Prbcur^{\otimes2}_\beta \arrow[ld, dashed] \\
                          & \ell_x \arrow[u]        &                                  
\end{tikzcd}
}
{\caption{Overview of the couplings between the graphical representations of the Ising model. 
Dashed arrows point towards the distribution obtained from taking a uniform even subgraph. 
Full arrows indicate a union with another percolation measure. 
} 
\label{fig:couplings}}
\end{figure}

\subsection{The phase transitions of the graphical representations} 
Having introduced the various graphical representations, we now turn our attention to their phase transitions and how they are connected. Suppose that $\mathbb{G}$ is an infinite graph embedded in $\mathbb{R}^d$ on which $\mathbb{Z}^d$ acts by translation by an $\mathbb{R}$-linearly independent family of vectors $(v_j)_{1\leq j\leq d}$. We will call such a graph $d$-periodic (or bi-periodic for $d=2$). We will mainly be concerned with $\mathbb{G} = \Z^d$ for $d \geq 2$ or  $\mathbb{G} = \mathbb{H}$ where $ \mathbb{H}$ is the hexagonal lattice in two dimensions.

\begin{table}
\begin{center}
\begin{tabular}{|c|c|c|c|} 
  \hline
&  $\beta$ & $p$  & $x$  \\ 
    \hline
    $\beta$  &\cellcolor{black!10}    &  $\frac{1}{2} \log(1-p)$ &  $\arctanh(x) $ \\ 
  \hline
$p$ & $1- e^{-2 \beta}$ & \cellcolor{black!10}    &  $\frac{2x}{x+1}$ \\ 
  \hline
$x$  & $\tanh(\beta)$  &  $\frac{p}{2-p}$& \cellcolor{black!10}    \\ 
  \hline
\end{tabular}
\end{center}
\caption{\label{table:parametrization_overview}
The parameters $\beta,p,x$ in standard parametrisations of $\mathbf{I}_\beta,\Prbcur_\beta,\phi_p,\ell_x$ are always assumed to relate according to the table. This over-determination is convenient for the individual parametrisations and follows the literature standard. \label{table:parametrizations}}
\end{table}

We consider one of the parametrised families of translation invariant infinite volume measures $\nu \in \{\ell^0_{x,\mathbb{G}}, \Prbcur_{\beta,\mathbb{G}}, \Prbcur_{\beta,\mathbb{G}}^{\otimes 2}, \phi^0_{p,\mathbb{G}}\}$.
Further, whenever we parametrise a measure by a parameter that is not its natural parametrisation, e.g. using $\ell_\beta$ instead of $\ell_x$, we implicitly use the bijections between the parametrisations that are summarised in Table \ref{table:parametrizations}.
For each of the models, there is a \emph{percolative phase transition} defined by
$$
\beta_c^{\perc}(\nu) = \inf \{\beta \geq 0 \mid  \nu_{\beta}\lbrack 0 \cc \infty \rbrack>0 \},  
$$
with the convention $\inf \emptyset=\infty.$ Recall that if $ \nu_{\beta}\lbrack 0 \cc \infty \rbrack>0,$ we say that the model \emph{percolates} at $\beta$.

One may also consider related but a priori different phase transitions corresponding to expected cluster sizes and regimes of exponential decay. These are defined by
$$
\beta_c^{\clust}(\nu) = \inf \{\beta \geq 0 \mid  \nu_{\beta} \lbrack \abs{\mathcal C_0} \rbrack = \infty \}, 
$$
as well as
$$
\beta_c^{\exp}(\nu) = \sup \{\beta \geq 0 \mid  \exists c,C >0,  \forall v \in \mathbb{G}:  \nu_{\beta} \lbrack 0 \cc v \rbrack \leq C e^{-c \abs{v}} \}. 
$$
For a translation invariant measure $\nu$, it is straightforward to see that
\begin{align*}
\beta_c^{\exp}(\nu) \leq \beta_c^{\clust}(\nu) \leq \beta_c^{\perc}(\nu).
\end{align*}
Furthermore, one says that a phase transition is \emph{sharp} if  $\beta_c^{\exp}(\nu) =  \beta_c^{\perc}(\nu)$. 
In general, the phase transition of the random-cluster model on a lattice is sharp (see \cite{DCsharpness}), but examples exist of lattice models in which the phase transition is not sharp, with the famous example of the intermediate phase described in the Berezinski-Kosterlitz-Thouless phase transition predicted in  \cite{berezinskii1971destruction, kosterlitz1973ordering} and rigorously proven in 
\cite{frohlich1981kosterlitz}.

The relations \eqref{eq:phi_Ising}, \eqref{eq:dc_Ising} together with sharpness of the Ising correlation function, proven in \cite{aizenman1987phase}, imply that 
\begin{align} \label{eq:phase_transtion_relations}
\beta_c^{\perc}(\phi) = \beta_c^{\clust}(\phi) = \beta_c^{\exp}(\phi)
= \beta_c
= \beta_c^{\perc}(\Prbcur^{\otimes 2}) = \beta_c^{\clust}(\Prbcur^{\otimes 2}) = \beta_c^{\exp}(\Prbcur^{\otimes 2}), 
\end{align}
where $\beta_c$ is the critical inverse temperature of the Ising model.

The main concern of this paper is the phase transition of the two remaining graphical representations, $\ell_x$ and $\Prbcur_\beta$. The couplings from \Cref{thm:couplings} immediately imply stochastic domination
\begin{align} \label{eq:elementary_stochastic_domination}
\ell_\beta \preceq \Prbcur_\beta \preceq \phi_\beta. 
\end{align}
Thus, for $\# \in \{\perc,\clust, \exp \}, $ it holds that 
$$
\beta_c^{\#}(\ell) \geq \beta_c^{\#}( \Prbcur) \geq  \beta_c^{\#}(  \phi) = \beta_c. 
$$
Duminil-Copin posed the question, 
 whether the percolative phase transition for the random current is the same as for the Ising model.
\begin{question}[\cite{DC16} Question 1]
    For infinite $\mathbb{G}$, does it hold that 
$
\beta_c^{\perc}( \Prbcur) =  \beta_c
$?
\end{question}

As one of the findings of this paper, we note that percolative loop $\mathrm{O}$(1) phase transition depends a lot on the lattice. A more nuanced picture appears from additionally considering the phase transitions $\beta_c^{\clust}$ and $\beta_c^{\exp}$. 
Our main result (Theorem \ref{main theorem}) implies that for $\mathbb{G} = \Z^d$ and $d \geq 2$ or  $\mathbb{G} = \mathbb{H}$,   
\begin{align*}
\beta_c^{\clust}(\ell) = \beta_c. 
\end{align*} 
In particular, 
$$
\beta_c^{\clust}(\ell)=  \beta_c^{\exp}(\ell)
= \beta_c^{\exp}(\Prbcur) = \beta_c^{\clust}(\Prbcur) = \beta_c. 
$$

\begin{remark}[On generality]
    The techniques we employ mostly rely on the fact that $\mathbb{Z}^d$ descends to a graph on the torus, which allows us to utilise the non-trivial topology of the torus to control the existence of large clusters. As such, the result should carry over to locally finite $d$-periodic graphs embedded in $\R^d$.

    The reason for not reflecting this in our statement is a) making it more accessible and b) the fact that a lot of the literature on the Ising model that we use is stated specifically for $\mathbb{Z}^d$ - for instance, \cite{Bod05,duminil2020exponential,Pis96}. We expect no particular difficulty in extending those results to lattices with suitable symmetries but prefer not to get bogged down on the possible generality of the results on which we rely.

    One may note that none of the above papers are necessary for treating the planar case where planar duality (introduced in Section \ref{sec:planar}) and sharpness of the phase transition does the necessary work for us. In particular, our results apply to the hexagonal lattice, which we shall discuss further in Section \ref{sec:hexagonal}.

    The assumption of having some structure is necessary for the result, as there are examples of graphs where the phase transitions for the loop $\mathrm{O}$(1) model $\ell$ and the random current model $ \Prbcur$ are non-unique \cite{monster_paper}. 

\end{remark}

\subsection{All phase transitions coincide for $\Z^2$}
On the square lattice, the overall picture is well-understood. For the Ising model, the existence of a phase transition was proven by Peierls \cite{peierls1936ising} and the exact value of $\beta_c = \frac{1}{2} \log\left(1+ \sqrt{2}\right) $ was proven by Onsager \cite{onsager1944crystal}. 
The uniqueness of the phase transition at the point $\beta_c$ for the loop $\mathrm{O}$(1) model was noted in \cite{GMM18} and the corresponding result for the single random current follows directly by coupling as per \Cref{thm:couplings}. 
For completeness, we present the proof here. It builds on a result of Higuchi for which we include an independent proof inspired by \cite[Proposition 2]{camia2020exponential}.

\begin{figure}
    \centering
    \includegraphics[scale=0.75]{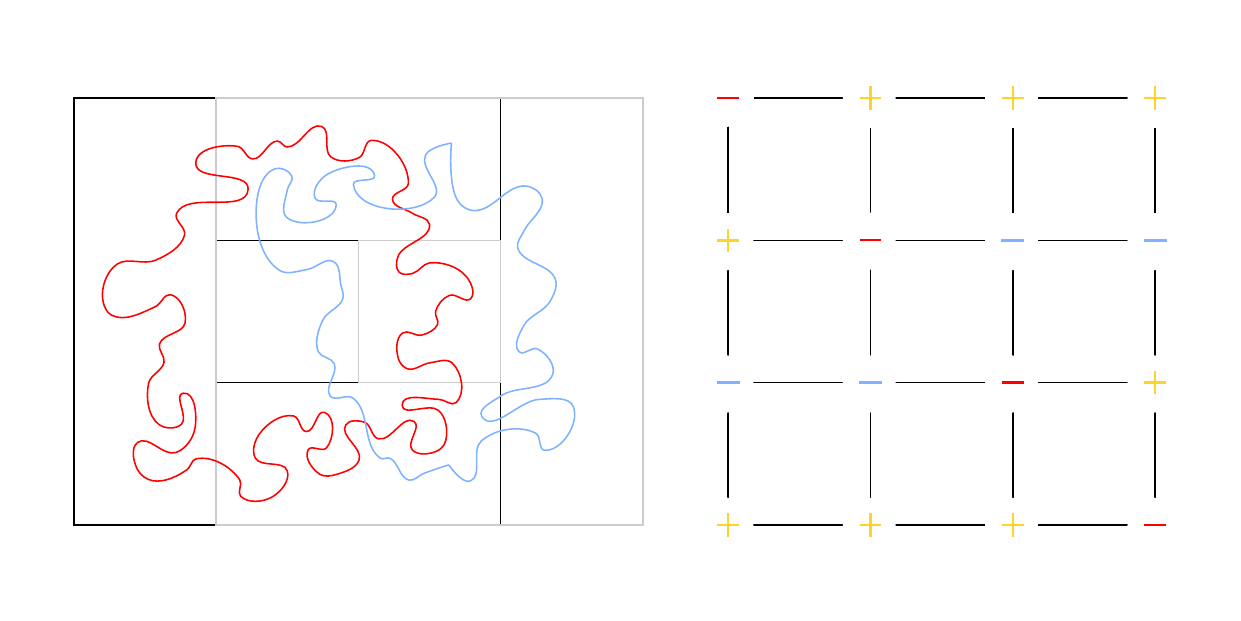}
    \caption{On the left: Two adjacent annuli (one in black, one in gray), both of which are good, as witnessed by the bold red and pale blue $*$-circuits of $-$'s. \\
    On the right: A spin configuration near an intersection between two circuits. It is not true that the red and blue $*$-circuits must share a vertex. However, by the Jordan Curve Theorem, their piecewise linear extensions must intersect at a point, which possibly lies inside a face of $\mathbb{Z}^2$. But this implies that the two circuits belong to a common $*$-cluster.}
    \label{fig:Star_circuits}
\end{figure}

First, recall the notion of $*$-connectivity for planar graphs. Two vertices $v,w$ of $\mathbb{Z}^2$ are said to be $*$-adjacent if they are adjacent to the same face. 
Compare the ordinary notion where two vertices are adjacent if they are adjacent to the same edge. 
Accordingly, we have a notion of $*$-paths (see Figure \ref{fig:Star_circuits} for examples). A set of vertices $W$ is said to be $*$-connected if any pair of points $v,w\in W$ can be joined by a $*$-path contained in $W$. For an Ising configuration $\sigma\in \{-1,+1\}^{V(\mathbb{Z}^2)},$ a $*$-connected $+$-cluster is a maximal $*$-connected set of vertices that are coloured $+$.
Ordinary clusters are similarly defined with respect to ordinary connectivity. Since ordinary connectivity implies $*$-connectivity, a $*$-cluster decomposes into a disjoint union of ordinary clusters. 

Let $\mathcal{C}_\infty^{+,*}$ be the event that there is an infinite  $*$-connected $+$-cluster and let $\mathcal{C}_\infty^{-,*}$ be the event that there is an infinite  $*$-connected $-$-cluster.

Finally, we remind the reader that a boundary condition for the Ising model on a finite subgraph $G$ of $\mathbb{Z}^2$ is given by fixing a configuration on the vertices adjacent to $G$. Thus, we get an Ising model on $G$ conditioned on these extra boundary spins. Just as in the case of boundary conditions for percolation configurations, these boundary conditions are partially ordered (see e.g. \cite[Section 3.6.2]{friedli2017statistical}). 

\begin{proposition}\label{prop:higuchi}
    Let $\beta < \beta_c$ then 
    $
     \mathbf{I}_{\beta, \Z^2} [\mathcal{C}_\infty^{-,*}]= \mathbf{I}_{\beta,\Z^2} [\mathcal{C}_\infty^{+,*}] = 1.
$
\end{proposition}

\begin{proof} 
We let $A(n) = \Lambda_{2n} \backslash \Lambda_n$ be an annulus. 
If there is a $*$-circuit in $A(n)$ of $-$ spins encircling the inner boundary, then we say that the annulus $A(n)$ is \emph{good}. Note that $A(n)$ is not good if and only if there is a (usual) path of $+$-spins connecting the inner boundary to the outer boundary.

Let us first prove that 
\begin{align} \label{eq:old_claim}
\mathbf{I}^+_{\beta, A(n)} \left[ A(n) \text{} \is \text{} \good \text{}\right]  \to 1. 
\end{align}
Let $v \in \partial \Lambda_n$ and let $\mathcal{C}^+_v$ be the $+$-cluster of $v$. 
Then, since $\beta < \beta_c,$ by sharpness \cite[Theorem 3]{higuchi1993coexistence}, there exists a $c> 0$ such that  $\mathbf{I}_{\beta, A(n)}^{+} [ \abs{\mathcal{C}^+_v} \geq n ] \leq e^{-cn} $. Thus, by a union bound, we obtain \eqref{eq:old_claim} as follows:
\begin{align*}
\mathbf{I}^+_{\beta, A(n)} \left[A(n)   \text{ is not good}\right] 
& \leq \sum_{v\in \partial \Lambda_n} \mathbf{I}^+_{\beta, A(n)} \lbrack \abs{\mathcal{C}^+_v} \geq n \rbrack \\
& \leq 8n e^{-cn}\to 0.
\end{align*} 
By monotonicity in boundary conditions, the same convergence holds for arbitrary boundary conditions.

Now, we look at all $(n,2n)$-annuli with centers in $n \mathbb{Z}^2$, i.e. $\{ A(n)+nk \}_{k \in \Z^2} $. For every $k \in  \Z^2$ define $X(k) = \id \lbrack A(n)+nk\text{ is good}\rbrack$. Then, we consider $\{X(k)\}_{k \in \Z^2}$ as a spin model and show that it percolates for sufficiently large $n$. 
It follows from \eqref{eq:old_claim} and DMP for the Ising model that
\begin{align*}
\mathbf{I}_{\beta, \Z^2} \left[ X(k) = 1 \mid  \{ X(j) \}_{j \in \Z^2: \abs{j-k}\geq 4} \right] \geq \mathbf{I}^+_{\beta, A(n)} \left[A(n) \text{ is good}\right]   \to 1, 
\end{align*}
so we can use \cite[Theorem 0.0]{liggett1997domination} to dominate the process $\{X(k)\}_{k \in \Z^2}$ from below by independent Bernoulli random variables with some parameter $p_n$ where  $p_n \to 1$ as $n\to \infty$. Hence, we can choose $n$ such that $p_n$ is above the threshold for site percolation (which is strictly smaller than 1 \cite{peierls1936ising}). 
By planarity, for every edge $(v,w)$ of $\mathbb{Z}^2$, if $A(n)+nv$ and $A(n)+nw$ are both good, any choice of corresponding circuits of $-$'s must belong to the same $*$-cluster of $-$ spins (a piecewise linear interpolation of the circuits can be seen to intersect), see Figure \ref{fig:Star_circuits}. 
In particular, percolation of the good annuli implies an infinite $*$-cluster of $-$ spins. 
Thus, $\mathbf{I}_{\beta, \Z^2} [\mathcal{C}_\infty^{-,*}] = 1$ and by spin-flip symmetry of the free Ising measure, it follows that 
 $\mathbf{I}_{\beta, \Z^2}[\mathcal{C}_\infty^{+,*}] = 1$.
\end{proof}

\begin{theorem}[{\cite[Theorem 1.3]{GMM18}}] \label{thm:two_dimensions} 
It holds that
$\beta_c^{\perc}(\ell_{\mathbb{Z}^2}) = \beta_c$. It follows that 
$\beta_c^{\#}(\nu) = \beta_c$ for all  $\nu \in \{\ell_{\mathbb{Z}^2}, \Prbcur_{\mathbb{Z}^2}, \Prbcur_{\mathbb{Z}^2}^{\otimes 2}, \phi_{\mathbb{Z}^2} \} $ and  
$\# \in \{\perc,\clust, \exp \} $. 
\end{theorem} 
\begin{proof}
The argument relies on planar duality, discussed in \Cref{sec:planar}. Namely, for a dual Ising model on $(\Z^2)^*$ with $\beta^*<\beta_c$ (corresponding to primal parameter $\beta>\beta_c$), the interfaces between spins of different signs form a percolation configuration with law of the loop $\mathrm{O}$(1) model $\ell_{\beta,\Z^2}$ by \Cref{prop:ES}. Here, parameters are related according to  \eqref{eq:dual_parameter}.

The interface of a finite spin cluster is a disconnected union of simple cycles with a uniquely determined outer interface which is the unique component which encircles the spin cluster.
Observe that $*$-connections imply connectivity of outer interfaces. More precisely, consider a $*$-cluster which does not contain an infinite ordinary cluster. The outer interfaces of all ordinary clusters belonging to the $*$-cluster are connected.
This implies that the interfaces of an infinite $*$-cluster which does not contain an ordinary infinite cluster contains an infinite connected set of edges. 
Therefore, by \Cref{prop:higuchi}, $\ell_{\beta,\Z^2}$ percolates. Combining this with \eqref{eq:phase_transtion_relations} and  \eqref{eq:elementary_stochastic_domination}, we now have
$$
\beta_c=\beta^{\perc}_c(\phi_{\mathbb{Z}^2})\leq \beta^{\perc}_c(\ell_{\mathbb{Z}^2})\leq \beta_c,
$$
which proves the first item.

The second statement also follows from \eqref{eq:phase_transtion_relations} and \eqref{eq:elementary_stochastic_domination}. 
\end{proof}
Since the uniform even subgraph is the loop $\mathrm{O}$(1) model $\ell_x$ for $x=1,$ the results settle percolation of  ${\UEGop}_{\Z^2}$.
\begin{corollary} \label{Z2_percolates}
 ${\UEGop}_{\Z^2}$ percolates. 
\end{corollary}

\section{The uniform even subgraph}
In the analysis of the uniform even subgraph we take a group-theoretic approach. 
This approach has merit in that it produces rather concise proofs and
brings to the fore the structure that makes the proofs work.
The relation to Gibbs measures of the UEG in the DLR formalism common in statistical mechanics is explained in Section \ref{sec:Gibbs of UEG}.
To follow the main story of this paper, \Cref{prop:separated-graphs} is the key takeaway from this section. It
roughly states that the uniform even subgraph gets decoupled by the existence of a separating surface. The proposition, which concludes \Cref{sec:marginals-of-UEG}, will be vital in the sequel to compare the loop O($1$) model on the torus to the loop O($1$) model on $\mathbb{Z}^d$. As such, \Cref{sec:Main arguments}  can be read independently of \Cref{sec:even_perco_Z_3}-\ref{section:UEG measures}.

We begin this section by collecting general properties of the uniform even subgraph. We then apply the theory to prove \Cref{d3percolates}, that the uniform even subgraph of $\Z^d$  percolates for $d \geq 3$. Further remarks on even percolation follow in \Cref{sec:further remarks}, before we extend the results to percolation of the loop $\mathrm{O}(1)$ model in \Cref{thm:non_tri_phase_transition}. 
Finally, we give a detailed comparison between different constructions of uniform even subgraphs on infinite graphs in  \Cref{section:UEG measures}, and characterise infinite volume Gibbs measures.

In the treatment of the uniform even subgraph, we take a generalist point of view and regard "uniform" as synonymous with "invariant with respect to a group action".
For a given graph $G =(V,E),$ the space $\Omega(G)$ is a group under point-wise addition modulo 2 - indeed, a $\Z_2$-vector space.
This addition corresponds to symmetric difference on sets of edges, denoted $\triangle$, and we will use these notions interchangeably. 
The symmetric difference of two even graphs is again even, and so is the empty subgraph, so the set of even subgraphs $\Even(G)$ is a $\Z_2$-linear subspace of $\Omega(G).$
For a connected finite graph, a basis for this subspace (sometimes denoted the fundamental cycle basis) is indexed over the edges in the complement of a spanning tree.
Furthermore, it follows by \eqref{eq:free def} below that $\Even(G)$ is closed. 
We define the uniform even subgraph to be the Haar probability measure on this group of even subgraphs.

The measure thus constructed is known in the literature as the free uniform even subgraph when the graph is finite, and as the wired uniform even subgraph when the graph is infinite.
While previously studied constructions of limit measures coincide with the Haar measure (see section \Cref{section:limit-constructions}), 
this property has, to the best of our knowledge, not been emphasised before, e.g. in \cite{angel2021uniform, GJ09}.
We demonstrate its merits in \Cref{d3percolates}, where we prove that the uniform even subgraph of $\Z^d$ percolates for $d \geq 3$.

\subsection{Marginals of the \texorpdfstring{$\UEG$}{UEG}}\label{sec:marginals-of-UEG}

We will repeatedly let $G = (\mathbb V, E)$ be a spanning subgraph of a locally finite, connected graph $\mathbb{G} = (\mathbb{V,E})$. $\mathbb{G}$ may be finite but at most countably infinite. By slight abuse of terminology, we say that $G$ is finite when $E$ is finite.

A standard approach to infinite volume measures in statistical mechanics is to consider conditional distributions (e.g. the conditional distribution of a configuration $\omega$ in $G$ given a configuration $\omega_0$ in $\mathbb{G}\setminus G=(\mathbb{V,E}\setminus{E})$)  rather than marginals \cite[p.270]{friedli2017statistical}. 
This approach has merit in the treatment of the random-cluster model with \Cref{DMP} as a cornerstone of the theory.
However, for the uniform even subgraph, and in turn the loop $\mathrm{O}(1)$ model, conditional probabilities are less useful. Specifically, $\ell_{x,G}[\cdot|\;\omega_1^c]$
is not, in general, a loop $\mathrm{O}$(1) model with  (topological) boundary conditions as defined in \eqref{eq:loop_with_bc}, since the sources on the boundary may force edges within $G$. 
Instead, it is more tractable  to consider the marginal of the measure $\UEG_\mathbb G$ on $G$ -  that is, $\UEG_\mathbb{G}\vert_{G}$. As an indication that this is natural, the uniform even subgraph of $G$ with wired boundary conditions can (under suitable conditions) be realised as $\UEG_\mathbb{\mathbb G}\vert_{G}$ whereas wired boundary conditions are extremal for the random-cluster model (see \Cref{rem:boundary_conditions}).

The $\Z_2$-vector space $\Omega$ is compact in the product topology\footnote{equivalently, the topology of pointwise convergence or the topology generated by cylinder events.} 
 so it admits a unique Haar measure normalised to probability. 
Notice that this Haar probability measure is $\P_{1/2}$ (the relation to construction via Kolmogorov's  extension theorem is discussed in \Cref{remark:Haar-Kolmogorov}).

Every closed subgroup of $\Omega$ is also Abelian and compact, so it has a corresponding Haar probability measure. 
Since $G$ is locally finite, we can define the source map 
$\partial : \{0,1\}^E \to\{0,1\}^V\cong\mathcal{P}(V)$  by mapping edges to incident vertices and extending linearly and continuously. 
The source map $\partial$ is a continuous homomorphism so the group of even subgraphs, is identified as
\begin{equation}\label{eq:free def}
    \Even(G) = \{\omega\in \Omega(G) |\; \partial \omega = \emptyset\} = \ker \partial,
\end{equation}
and under this identification it follows that it is a closed subgroup.
The uniform measure on this group, i.e. the Haar probability measure, is denoted $\UEG_G$. 

For a fixed subgraph $G=(\mathbb V,E)\subset\mathbb{G=(V,E)}$ the inclusion $\iota_{E}:{E}\to \mathbb E$ induces the projection $\pi_{G}:\{0,1\}^\mathbb E\to \{0,1\}^{E}$ also known as the restriction to $E$. 
By construction of the product topology, the projection is continuous
and the marginal on $G$ of any Borel measure $\mu$ on $\{0,1\}^E$ is defined to be the pushforward along $\pi_{G}$,
denoted $\mu|_G$ \footnote{This is not be confused with the plain restriction of a measure to a subspace, viz. $ \Omega(G)\subset \Omega(\mathbb G)$. That is, $\mu(A)\ne \mu|_G(A)=\mu(\pi_G^{-1}(A))$ for general measurable $A\subset \Omega(G)$.}.

Crucially, $\pi_G$ is a homomorphism, which, along with the following general fact\footnote{Commutativity in this context is only used to ensure uniqueness of invariant measures. A more general statement is true for invariant measures on non-Abelian groups but this is not relevant here.}, 
explains the good behaviour of marginals of Haar measures.

\begin{lemma}\label{pushhom}
    Let $\Gamma,\Gamma'$ be Abelian 
    compact topological groups and $f:\Gamma\to \Gamma'$ a continuous homomorphism.  
    Let $\mu$ be a Haar probability measure on $\Gamma$. If $f$ is surjective, then $f_*\mu$ is a Haar probability measure on $\Gamma'$. 
    For a general continuous homomorphism, we obtain the Haar measure on the image $f(\Gamma)$.
\end{lemma}
\begin{proof}
    The homomorphism property of $f$ implies translation invariance of $f_*\mu$ on its support. The rest of the claim follows from the uniqueness of the Haar measure. 
\end{proof}

As a consequence of \Cref{pushhom}, in order to characterise marginals of the $\UEG$ we need only determine the range of the restriction map. 
To this end, we need the following variant of Kőnig's lemma. 
Define $(\mathbb{G}\setminus G)^\circ$ to be the subgraph of $\mathbb G\setminus G$ with isolated vertices removed.

\begin{lemma}\label{lem:konig surjective}
    For any locally finite, connected, countably infinite graph, the boundary map $\partial$ is surjective. 
\end{lemma}
\begin{proof}
 Let $\mathbb{G}=\mathbb{(V,E)}$ be as hypothesised. For any $v\in \mathbb{V}$, by Kőnig's lemma, there is a spanning subgraph $\gamma_v\subset \mathbb{E}$ such that $\partial\gamma_v=\{v\}$.  
    Choosing such $\gamma_v$ for each $v$ in a finite subset $W\subset \mathbb{V}$ and producing the symmetric difference yields $\partial\gamma_W=W$.
    Now, assume $U\subset \mathbb{V}$ is infinite. In order to construct $\gamma$ with boundary $U$, proceed as follows. 
    Choose an exhausting sequence of finite subgraphs $G_n \Uparrow \mathbb{G}$ in such a way that $(\mathbb{G}\setminus G_n)^\circ$ has no finite components. This is possible since, if $G_n$ is given, replace it 
  by the union of $G_n$ with the finitely many finite components of $(\mathbb{G}\setminus G_n)^\circ$.
  
    Then, recursively define $\gamma_{n}$ with 
    $$\partial\gamma_{n}=W_n=W\cap V(G_n) \quad{\text{and}}\quad \gamma_n|_{G_{n-1}}=\gamma_{n-1}.$$
    Say $\gamma_{-1}=\emptyset$. Given $n$ and $\partial\gamma_{n-1}=W_{n-1}$, we observe that $G\setminus G_{n-1}$ has only infinite components, so there is $\tilde\gamma_n\subset E(G\setminus G_{n-1})$ with $\partial\tilde\gamma_{n}=W_n\setminus W_{n-1}$. Thus $\gamma_n=\tilde\gamma_n\triangle\gamma_{n-1}$ has $\partial\gamma_n=W_n$ and satisfies $\gamma_n|_{G_{n-1}}=\gamma_{n-1}$. 
    Then there is a limit, $\gamma_n\to\gamma$, and if $v\in V(G_n)$ then $\partial\gamma(v)=\partial\gamma_n(v)=\id_U(v)$. 
\end{proof}

As we are about to see, the marginal of the UEG can be understood in terms of topological boundary conditions. However, one has to  reinterpret some definitions when dealing with infinite graphs. 
Recall from \Cref{rem:boundary_conditions} that a topological boundary condition $\xi$ is a partition of the set $\partial_vG$ for $G=(V,E)$. 


In the random-cluster model in finite volume, an exploration $\mathbb{E}'\subset \mathbb{E}$ induces a boundary condition by the equivalence relation that $x$ is connected to $y$ by a finite path in $\mathbb{E}'\setminus E$. Here instead, we use the equivalence relation  that $x,y$ are \emph{$\infty$-connected} if there is a subgraph $\gamma\subset \mathbb{E}'\setminus E$ with $\partial \gamma= \{x,y\}$ (this is what is called the wired augmentation in \cite{Unique_Forest}). Note that $\gamma$ might consist of two disjoint rays.
For our purposes, we will assume $\mathbb{E}'=\mathbb{E}$, so we do not need the distinction. We  refer to the induced equivalence relation on $\partial_vG$ as $\sim_\xi$ with equivalence classes $\xi=\xi_\mathbb{E}$.
Recall also, that when considering  $\omega\subset E$ as a subgraph of $G/\sim_{\xi},$ it is denoted by $\omega^\xi$.

There is precisely one infinite $\infty$-connected component in an infinite graph. Denote this by $\delta_\infty$, and consider it equally as a vertex in $G/\sim_{\xi}$. The boundary map $\partial$ on $G/\sim_{\xi}$ is modified so that $\partial(\omega^\xi)(\delta_\infty)=0$ for all $\omega$. 
In this definition, $\partial$ is still a homomorphism.

We can now extend \cite[Theorem 2.6]{GJ09} to also allow infinite subgraphs and by adding a description of the support of the marginal:
\begin{proposition}\label{prop:marginal-is-boundary-condition}    
    Let $G=(\mathbb V, E)\subset \mathbb{G}$. 
    The group of even subgraphs of $\mathbb{G}$ restricted to $G$ is
    \begin{equation}\label{eq:evens-with-boundary-condition}
        \pi_{G}(\Even(\mathbb G))= \Even^\xi(G)
        :=\{\omega\in\Omega(G)|\,\partial\omega^\xi=0\}.  
    \end{equation}
    Moreover, the marginal,  ${\UEGop}_\mathbb G|_G$, is given by the uniform (Haar) measure on $\Even^\xi(G)$.
\end{proposition}

\begin{proof}
    Consider $\eta\in \Even(\mathbb{G})$ and let $\omega=\eta|_G$. 
    To see that $\partial\omega^\xi=0$, we only have to check vertices $\delta$ in $G/\sim_\mathbb E$ such that $\delta\ne\delta_\infty$.
    If an $\infty$-connected component of $\eta|_{G^c}$ is finite, then it has an even number of sources in $\partial_vG$, i.e. vertices with $\partial(\eta|_{G^c})(v)=1$. Indeed, the sources may be connected pairwise in $\eta|_{G^c}$. The corresponding vertices of $\partial\omega$ are identified in $G/\sim_\mathbb E$.

    For the opposite inclusion, assume $\eta\subset E$ and $\partial\eta^\xi=0$. That is, for any equivalence class $\delta\subset\partial_vG$ corresponding to an $\infty$-connected component of $\mathbb{G}\setminus G,$ we have that $\partial\eta^\xi(\delta)=0$ with $\delta$ considered as a vertex of $G/\sim_{\xi}$. 
    We will construct $\tilde{\eta}\in \Even(\mathbb{G})$ such that $\tilde{\eta}\cap E=\eta.$ 
    
    For $\delta \in \partial_v G/\sim_{\xi},$ define $\eta_{\delta}\subset \mathbb{E}\setminus E$ as follows:  If $\delta\ne\delta_\infty$, there is an even number of vertices in $\partial_vG$ in the class $\delta$ with odd degree in $\eta$. We may pair those odd vertices such that they are connected in $\mathbb{G}\setminus G$ due to the assumption that the pair belongs to a finite connected component therein. 
    The symmetric difference of these subgraphs we take as $\eta_\delta$ and note $\partial(\eta_{\delta})=(\partial\eta)|_{\delta}$. 
    For $\delta_\infty$, we apply \Cref{lem:konig surjective} separately to the (ordinary) infinite components of $\mathbb{G}\setminus G$ corresponding to $\delta_\infty$ to find $\eta_{\delta}$  with $\partial(\eta_{\delta_\infty})=(\partial\eta)|_{\delta_\infty}$.

    By construction $\eta_\delta\cap\eta_{\delta'}=\emptyset$ if $\delta\ne\delta'$, so $\tilde{\eta}=\eta \triangle \left(\triangle_{\delta\in \partial_v G/\sim_{\mathbb{E}}} \eta_{\delta}\right)$ is well-defined. As desired, $\tilde{\eta}\cap E=\eta$ and $\partial \tilde{\eta}=0$.
\end{proof}

Recall, from \Cref{rem:boundary_conditions}, the wired boundary condition which identifies all vertices on the boundary. The set of wired even subgraphs is  
\begin{equation}\label{eq:def wired evens}
    \Wired(G)=\{\omega\in\Omega(G)|\,\partial\omega^1=0\}.
\end{equation}
For finite graphs $G$, then $\Wired(G)$ is the set of subgraphs of $G$ with an even number of sources on the boundary and where vertices in the interior are sourceless.

We write $\UEG_{{G}}^1$ for the uniform probability measure on $\Wired({G})$. 
Similarly, we may consider the free boundary condition, but only in the case of finite graphs do we refer to 
$\Even^0(G)=\Even(G)$ as the set of free even subgraphs. 
Free even subgraphs of infinite graphs are defined as limits of finite free even subgraphs, as elaborated upon in  \Cref{section:limit-constructions}.

\begin{lemma}\label{lemma:restriction-separation}
    Let $G=(\mathbb V,E)\subset \mathbb{G}$ and assume $(\mathbb{G}\setminus G)^\circ$ is $\infty$-connected. Then, 
    \begin{equation}\label{eq:marginal-wired}
        \pi_{G}(\Even(\mathbb G)) = \Wired(G), 
    \end{equation}
    and hence, 
    ${\UEGop}_\mathbb{G}|_{G}={\UEGop}_{G}^1$. 
    
    If moreover $(\mathbb{G}\setminus G)^\circ=\mathbb{G}\setminus G$, then 
    $
        {\UEGop}_\mathbb{G}|_G = \P_{\frac{1}{2},G}.
    $
\end{lemma}
\begin{proof}
The first statement is a special case of \Cref{prop:marginal-is-boundary-condition}. The  second follows from the observation that $(\mathbb{G}\setminus G)^\circ=\mathbb{G}\setminus G$ implies that every vertex of $G$ is considered boundary, hence $\Wired(G)=\Omega(G)$.
\end{proof}

\begin{figure} 
\tdplotsetmaincoords{70}{120} 
\tdplotsetrotatedcoords{0}{0}{0} 
\begin{tikzpicture}[scale=2,tdplot_rotated_coords,
                    rotated axis/.style={->,purple,ultra thick},
                    blackBall/.style={ball color = black!80},
                    borderBall/.style={ball color = white,opacity=.25}, 
                    very thick]

\foreach \x in {0,1,2}
   \foreach \y in {0,1,2}
      \foreach \z in {0,1,2}{
           \ifthenelse{  \lengthtest{\x pt < 2pt}  }{
             \draw (\x,\y,\z) -- (\x+1,\y,\z);
             \shade[rotated axis,blackBall] (\x,\y,\z) circle (0.025cm); 
           }{}
           \ifthenelse{  \lengthtest{\y pt < 2pt}  }{
               \draw (\x,\y,\z) -- (\x,\y+1,\z);
               \shade[rotated axis,blackBall] (\x,\y,\z) circle (0.025cm);
           }{}
           \ifthenelse{  \lengthtest{\z pt < 2pt}  }{
               \draw (\x,\y,\z) -- (\x,\y,\z+1);
               \shade[rotated axis,blackBall] (\x,\y,\z) circle (0.025cm);
           }{}

}

\colorlet{lightblue}{blue!50}

\path[fill=orange,fill opacity=0.5] (1,1,1) -- (0,1,1) -- (0,2,1) -- (1,2,1) -- (1,1,1); 

\path[fill=orange,fill opacity=0.5] (1,1,1) -- (2,1,1) -- (2,0,1) -- (1,0,1) -- (1,1,1); 

\path[fill=orange,fill opacity=0.5] (1,1,1) -- (0,1,1) -- (0,0,1) -- (1,0,1) -- (1,1,1); 

\path[fill=orange,fill opacity=0.5] (1,1,1) -- (2,1,1) -- (2,2,1) -- (1,2,1) -- (1,1,1); 

\draw[lightblue] (1,1,1) -- (1,1,2)--(1,2,2)--(1,2,1)--(1,1,1); 
\draw (1,1.8,1) node[scale=1.2, above] {\textbf{e}};   
\shade[rotated axis,borderBall] (1,1,1) circle (0.05cm);
\shade[rotated axis,borderBall] (1,2,1) circle (0.05cm);
\end{tikzpicture}
\caption{A subset of the lattice $\Z^3$, with a $\Z^2$ sheet highlighted with orange. For an edge $e$ in the sheet, we show the even graph (loop in blue) that we use to prove that $\Even(\Z^3)$ separates edges of $\Z^2$. \label{fig:plaquette_sticking_out} }

\end{figure}

The condition that $(\mathbb{G}\setminus G)^\circ=\mathbb{G}\setminus G$ and $(\mathbb{G}\setminus G)^\circ$ is $\infty$-connected is equivalent to saying that every edge in $G$ can be completed in $\mathbb{G}\setminus G$ to a possibly infinite loop in $\mathbb{G}$. In that case, we say that $\Even(\mathbb{G})$ \emph{separates edges} of $G$.  
Our first result follows as an application of \Cref{lemma:restriction-separation}.

\uegdpercolates*

\begin{proof} 

    The case of $d=2$ is \Cref{Z2_percolates}. 
   For $d \geq 3$ we show that $\Z^d$ contains a subgraph $H$ such that $\Even(\Z^d)$ separates edges of $H$ and such that $\Prb_{\frac{1}{2},H}$ percolates. Then, by \Cref{lemma:restriction-separation}, the marginal $\UEG_{\Z^d}|_H$ percolates and hence, so does $\UEG_{\Z^d}$. 

    In the case $d\ge 4$, $\Prb_{\frac{1}{2}, \Z^{d-1}}$ percolates (this, for example, follows from \Cref{prop:three_complications}, but was originally proven in \cite{campanino1985upper}) and therefore, we can choose $H$ to be a hyperplane of co-dimension 1.
    Every edge in the hyperplane $H$ is part of a plaquette\footnote{Recall that a plaquette of $\mathbb{Z}^d$ is a cycle of length 4.} which is not contained in the hyperplane, so $\Even(\Z^d)$ separates edges of $H$, see Figure \ref{fig:plaquette_sticking_out}.
    
    When $d=3,$ a hyperplane will not suffice, since $\Prb_{\frac{1}{2},\Z^2}$ does not percolate \cite{kesten1980critical}. 
However, if we let $H$ be the induced subgraph of $\Z^2\times \{0,1\}\subset \Z^3$ (that is, the graph with vertex set $\Z^2\times \{0,1\}$ and an edge $(v,w)$, whenever $(v,w)$ is an edge of $\Z^3$, see Figure \ref{fig:3d}),
    then  $\Prb_{\frac{1}{2},H}$ percolates. We defer the proof of this fact to \Cref{prop:three_complications}.

    We can separate edges on each sheet of $\Z^2\times \{0,1\}$ as before (with the loop around the plaquette sticking out upwards on the upper plane and with the loop around the plaquette sticking out downwards for each edge on the lower plane). Meanwhile, the connecting edges between the two sheets can be separated by doubly infinite cycles orthogonal to the slabs as sketched in Figure \ref{fig:3d}. 
\end{proof}


\begin{figure}

\tdplotsetmaincoords{70}{120} 
\tdplotsetrotatedcoords{0}{0}{0} 
\begin{tikzpicture}[scale=2,tdplot_rotated_coords,
                    rotated axis/.style={->,purple,ultra thick},
                    blackBall/.style={ball color = black!80},
                    borderBall/.style={ball color = white,opacity=.25}, 
                    very thick]

\foreach \x in {0,1,2}
   \foreach \y in {0,1,2,3}
      \foreach \z in {0,1,2,3}{
           \ifthenelse{  \lengthtest{\x pt < 2pt}  }{
             \draw (\x,\y,\z) -- (\x+1,\y,\z);
             \shade[rotated axis,blackBall] (\x,\y,\z) circle (0.025cm); 
           }{}
           \ifthenelse{  \lengthtest{\y pt < 2pt}  }{
               \draw (\x,\y,\z) -- (\x,\y+1,\z);
               \shade[rotated axis,blackBall] (\x,\y,\z) circle (0.025cm);
           }{}
           \ifthenelse{  \lengthtest{\z pt < 2pt}  }{
               \draw (\x,\y,\z) -- (\x,\y,\z+1);
               \shade[rotated axis,blackBall] (\x,\y,\z) circle (0.025cm);
           }{}

}

\foreach \x in {0,1,2}
   \foreach \y in {0,1,2,3}
      \foreach \z in {1,2}{
           \ifthenelse{  \lengthtest{\x pt < 2pt}  }{
             \draw[orange] (\x,\y,\z) -- (\x+1,\y,\z);
             \shade[rotated axis,blackBall] (\x,\y,\z) circle (0.025cm); 
           }{}
           \ifthenelse{  \lengthtest{\y pt < 2pt}  }{
               \draw[orange]  (\x,\y,\z) -- (\x,\y+1,\z);
               \shade[rotated axis,blackBall] (\x,\y,\z) circle (0.025cm);
           }{}
           \ifthenelse{  \lengthtest{\z pt < 2pt}  }{
               \draw[orange]  (\x,\y,\z) -- (\x,\y,\z+1);
               \shade[rotated axis,blackBall] (\x,\y,\z) circle (0.025cm);
           }{}

}

\foreach \x in {0,1,2}
   \foreach \y in {0,1,2,3}{
  \draw (\x,\y,2) -- (\x,\y,3);}

\foreach \x in {0,1,2}
   \foreach \z in {0,1,2,3}{
  \draw (\x,2,\z) -- (\x,3,\z);}

\foreach \x in {0,1,2}
   \foreach \z in {1,2}{
  \draw[orange] (\x,2,\z) -- (\x,3,\z);}

\path[fill=orange,fill opacity=0.5] (1,1,2) -- (0,1,2) -- (0,2,2) -- (1,2,2) -- (1,1,2); 

\path[fill=orange,fill opacity=0.5] (1,1,2) -- (2,1,2) -- (2,0,2) -- (1,0,2) -- (1,1,2); 

\path[fill=orange,fill opacity=0.5] (1,1,2) -- (0,1,2) -- (0,0,2) -- (1,0,2) -- (1,1,2); 

\path[fill=orange,fill opacity=0.5] (1,1,2) -- (2,1,2) -- (2,2,2) -- (1,2,2) -- (1,1,2);

\path[fill=orange,fill opacity=0.5] (1,2,2) -- (0,2,2) -- (0,3,2) -- (1,3,2) -- (1,2,2); 
\path[fill=orange,fill opacity=0.5] (1,2,2) -- (2,2,2) -- (2,3,2) -- (1,3,2) -- (1,2,2);

\colorlet{lightblue}{blue!50}

\path[fill=orange,fill opacity=0.5] (1,1,1) -- (0,1,1) -- (0,2,1) -- (1,2,1) -- (1,1,1); 

\path[fill=orange,fill opacity=0.5] (1,1,1) -- (2,1,1) -- (2,0,1) -- (1,0,1) -- (1,1,1); 

\path[fill=orange,fill opacity=0.5] (1,1,1) -- (0,1,1) -- (0,0,1) -- (1,0,1) -- (1,1,1); 

\path[fill=orange,fill opacity=0.5] (1,1,1) -- (2,1,1) -- (2,2,1) -- (1,2,1) -- (1,1,1); 

\path[fill=orange,fill opacity=0.5] (1,2,1) -- (0,2,1) -- (0,3,1) -- (1,3,1) -- (1,2,1); 
\path[fill=orange,fill opacity=0.5] (1,2,1) -- (2,2,1) -- (2,3,1) -- (1,3,1) -- (1,2,1);

\foreach \x in {0,1,2}
   \foreach \y in {0,1,2,3}{
  \draw[orange] (\x,\y,1) -- (\x,\y,2);}


\draw[lightblue] (1,2,2) -- (1,2,3)--(1,3,3)--(1,3,2)--(1,2,2); 
\draw[thick,lightblue] (1,1,0) -- (1,1,1)--(1,1,2)--(1,1,3);

\draw (1,1,1.2) node[scale=1.2, above] {$e_v$};   
\draw (1,2.8,1.9) node[scale=1.2, above] {$e_h$};   

\end{tikzpicture}

\caption{The situation in the proof in the case $d=3$ in Theorem \ref{d3percolates} with the edges in the subgraph $\Z^2 \times \{0,1\}$ coloured orange (except for the two edges $e_h$ and $e_v$). The horizontal edge $e_h$ is separated by the light blue cycle, whereas the edge $e_v$ is separated by the bi-infinite path, the first part of which is also coloured light blue. \label{fig:3d}}
\end{figure}

Another application of \Cref{lemma:restriction-separation} is the following. Note that $\infty$-connectivity reduces to ordinary connectivity when $\mathbb{G}$ is finite. 
\begin{corollary} \label{Surprising fact}
For all $d\geq 2$ and any $k,$
$\UEGop_{\Lambda_{k+1}}|_{\Lambda_k}=\UEGop^{1}_{\Lambda_{k+1}}|_{\Lambda_k}=\UEGop^1_{\Lambda_k}$.
\end{corollary}
\begin{proof}
It is a direct consequence of \Cref{lemma:restriction-separation} that $\UEGop_{\Lambda_{k+1}}|_{\Lambda_k}=\UEGop^1_{\Lambda_k}$, since $(\Lambda_{k+1}\setminus\Lambda_k)^\circ$ is connected. 
We have $\pi_{\Lambda_{k}}(\Wired(\Lambda_{k+1})) = \pi_{\Lambda_{k}}(\Even(\Lambda_{k+1}))$ 
since $\Even(\Lambda_{k+1})\subset \Wired(\Lambda_{k+1})$ and $\pi_{\Lambda_k}(\Wired(\Lambda_{k+1}))\subset \Wired(\Lambda_{k})=\pi_{\Lambda_k}(\Even(\Lambda_{k+1}))$.
So the last equality holds by appeal to \Cref{pushhom}.
\end{proof}
This result, while simple in itself, runs counter to common intuition in statistical mechanics because it shows that the uniform even subgraph is extremely insensitive to boundary conditions despite being strongly constrained.

We continue by showing a general condition under which this insensitivity to boundary conditions occurs.
On a graph $\mathbb{G}=(\mathbb{V,E}),$ we say that a path $\eta\subset E$ connects $G_1=(\mathbb V,E_1)$ and $G_2=(\mathbb V,E_2)$ if its end-points are non-isolated vertices in $G_1,G_2$ respectively.
Recall that $W\subset\mathbb V$ is called a \emph{separating surface} between disjoint edge sets $E_1, E_2$ if any path connecting $(\mathbb V,E_1)$ and $ (\mathbb V,E_2)$ visits $W$.

\begin{figure}
    \centering
\includegraphics[scale=0.65]{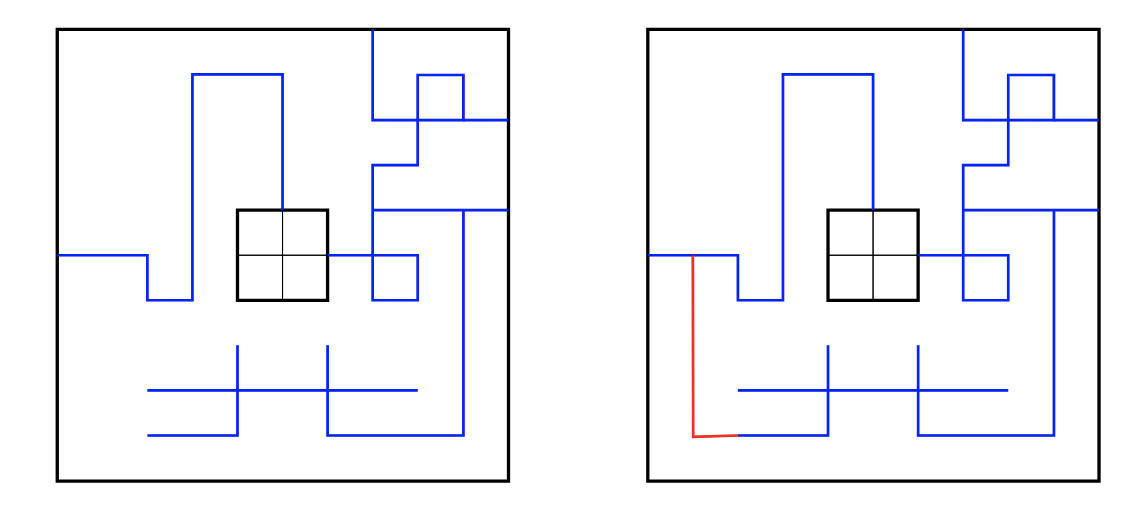}
    \caption{An example of two subgraphs of $\Lambda_{5}$ that only differ from each other outside of $\Lambda_1$. On the left, the marginal of $\UEG^1$ on $\Lambda_1$ is different from the marginal of $\UEG$ due to the lack of a \emph{connected} separating surface between the inner and outer boundary. On the right, the two marginals are the same due to the presence of the red path. One may note that the existence of a connected separating surface in a percolation configuration is not an increasing event.}
    \label{fig:Boundaries_matter}
\end{figure}

\begin{proposition}\label{prop:separated-graphs}
    Assume $E_1, E_2, E_3\subset E$ are pairwise disjoint and that $E_1,E_2$ are finite.
    Let $G_j=(V,\bigcup_{i\le j} E_i)$.
    If the induced graph $W$ of $E_2$ is connected 
    and a separating surface between $E_1$ and $E_3,$
    then 
    $
        \pi_{G_1}(\Even(G_2))=\pi_{G_1}(\Even(G_3)).
    $

    Hence, ${\UEGop}_{G_2}|_{G_1}={\UEGop}_{G_3}|_{G_1}.$
\end{proposition}
\begin{proof}
    It suffices to check $\pi_{G_1}(\Even(G_3))\subset \pi_{G_1}(\Even(G_2))$. The only non-trivial case to check is $\gamma\in\Even(G_3)$ intersecting both $E_1$ and $E_3$.
    Without loss of generality, assume $\gamma\in \Even(G_3)$ is connected. 
   
    Since $\gamma$ is even and connected, it may be parametrised as an edge self-avoiding closed path such that the first edge lies in $E_1$. This gives rise to a partition, $S$, of $\gamma_{|E_2\cup E_3}$ consisting of maximal subpaths contained in $E_2\cup E_3$.
    Now, we show that for each $s\in S,$ there is $h_s\in\Even((V,E_2\cup E_3))$ such that $s\triangle h_s\subset E_2$. 
    The situation is illustrated in Figure \ref{fig:Maximal_subpaths}.

    \begin{figure}
        \centering
        \includegraphics[scale =0.6]{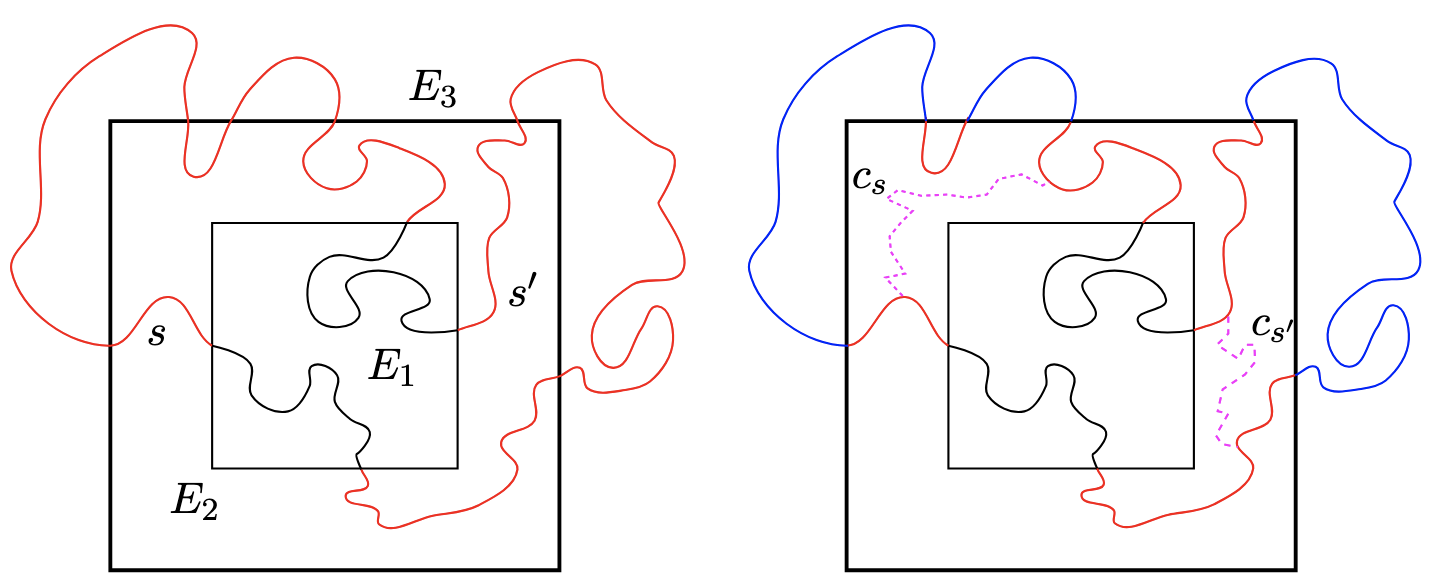}
        \caption{An illustration of the division of a loop $\gamma$ into maximal subpaths. On the left, we partition $\gamma$ into two paths $s$ and $s'$ connecting the boundary of $G_1$, coloured red. To the right, these may then be further subdivided into paths completely contained in either $E_2$ or $E_3$, coloured red and blue respectively. By the assumption that $E_2$ contains a connected separating surface, any two red paths crossing from $E_3$ to $E_1$ must be connected to each other inside $E_2$, illustrated by the paths $c_s$ and $c_{s'}$ in dashed purple. This, in turn, allows us to construct $h_s$ by following $c_s$ and then following the part of $s$ going through $E_2$ and $E_3$. If we take the symmetric difference of $\gamma$ with $h_s$ and $h_{s'}$, we construct a loop in $E_2\cup E_1$, which intersects $E_1$ in the same edges as $\gamma$.}
        \label{fig:Maximal_subpaths}
    \end{figure}
  
    If $s\subset E_2$, $h_s=\emptyset$. Otherwise, we partition $s_{|E_2}$ into maximal subpaths of $\gamma$ with respect to the parametrisation. The first and last segments both connect between $E_1$ and $E_3$, so they both visit  $W$ as it is a separating surface. 
    Since $W$ is connected, these two segments are connected by a 
    (possibly trivial)   path $c\subset E_2$ which 
    may be chosen so that 
        it only intersects the given two segments of $s$ in its end-points.
    By assumption, the endpoints of $c$ lie on $\gamma$ so they are parametrised as $\gamma(t_0),\gamma(t_1)$, for some $t_0< t_1$. Since $\gamma(0)\in E_1$, $\gamma([t_0,t_1])\subset s$.
    Then, $h_s=\gamma([t_0,t_1])\triangle c$ is even and satisfies $s\triangle h_s\subset E_2$.
    Let $ h=\triangle_{s}h_s$ which is even and satisfies $h|_{E_3}=\gamma|_{E_3}$ and $h|_{E_1}=\emptyset$. So $\tilde \gamma=\gamma\triangle h\in\Even(G_2)$ and $\tilde\gamma_{|E_1}=\gamma_{|E_1}$
    which completes the proof.
\end{proof}

This tells us that, in the presence of a uniform even subgraph of $G_3$, observing the edges of $G_1$ gives no information about the even subgraph on the edges of $G_3$. This holds in the strongest possible sense:

\begin{corollary} \label{UEG independence}
Assume $E_1, E_2, E_3\subset E$ are pairwise disjoint and that $E_1,E_2$ are finite.
    Let $G_j=(V,\bigcup_{i\le j} E_i)$.
    If the induced graph $W$ of $E_2$ is connected and a separating surface $W$ between $E_1$ and $E_3$, and $\eta\sim \UEGop_{G_3},$ then $\eta|_{E_1}\perp \eta|_{E_3}$.
\end{corollary}
\begin{proof}
For $j\in \{1,3\}$ fix $\eta_j\in \pi_{(V,E_j)}(\Even(G_3))$ and define the set $\Omega_{\emptyset}^{\eta_1,\eta_3}$ as the subset of $\Omega_{\emptyset}(G_3)$ which agrees with $\eta_j$ on $E_j$. Since $E_1$ and $E_2$ are finite, one may note that 
$$
\UEG_{G_3}[\id_{\eta|_{E_1}=\eta_1}|\; \eta|_{E_3}]=\frac{|\Even^{\eta_1,\eta|_{E_3}}|}{\sum_{\eta_1'\in \pi_{(V,E_1)}(\Omega_{\emptyset}(G_3))} |\Even^{\eta_1',\eta|_{E_3}}|}.
$$
Thus, the claim follows if we can prove that $\Even^{\eta_1,\eta_3}$ is a translate of $\Even^{0,\eta_3},$ where $0$ denotes the empty graph.
By \Cref{prop:separated-graphs}, there exists $\eta \in \Even^{\eta_1,0}$. The translation map, $\eta'\mapsto \eta'\triangle \eta$ on $\Even(G_3)$, is injective and maps $\Even^{\eta_1,\eta_3}$ into $\Even^{0,\eta_3}$ and  $\Even^{0,\eta_3}$ into $\Even^{\eta_1,\eta_3}$. 
\end{proof}

\Cref{prop:separated-graphs} applies to $\Lambda_k$ and $\Z^d\setminus\Lambda_{2k}$ with the annulus in between as the separating surface, illustrated in \Cref{fig:Boundaries_matter}.
To extend this application to a supercritical random-cluster configuration, our first order of business in \Cref{sec:Main arguments} will be 
to prove that the absence of a connected separating surface  in an annulus is exponentially unlikely in the supercritical random-cluster model. 
This is instrumental to proving \Cref{thm:Loop O(1) mixing}.

\subsection{Percolation on \texorpdfstring{$\Z^3$}{Z3}}\label{sec:even_perco_Z_3}
\begin{figure}
    \centering
    \includegraphics[scale=1.5]{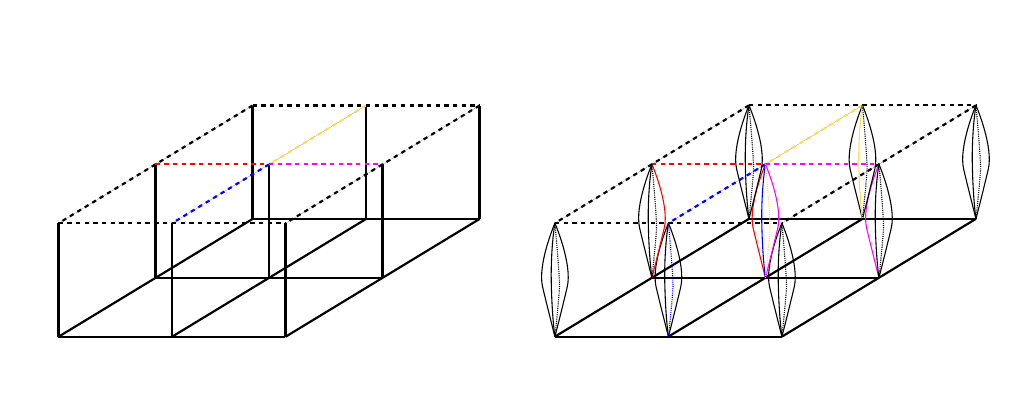}
    \caption{In order to get independence of $T\omega$, we replace every vertical edge of $\mathbb{Z}^2\times\{0,1\}$ by a parallel family of four vertical edges with lower edge weights. In order to construct $\tilde{T}\omega,$ we then colour the new vertical edges according to the horizontal ones.}
    \label{fig:Two_layered_surgery}
\end{figure}

The proof of percolation of $\UEG_{\Z^3}$ relies on the following result on Bernoulli percolation:
\begin{proposition} \label{prop:three_complications}
    The critical parameter for Bernoulli percolation on $\mathbb{Z}^2\times\{0,1\}$ is strictly smaller than $\frac{1}{2}$. In particular, 
    $\Prb_{\frac{1}{2},\Z^2\times \{0,1\}}$ percolates. 
\end{proposition}
\begin{proof}
We shall see that this is a relatively straightforward corollary of the fact that $p_c(\mathbb{P}_{\mathbb{Z}^2})=\frac{1}{2}$ \cite{kesten1980critical}. To see this, we define a map $T:\{0,1\}^{E(\mathbb{Z}^2\times \{0,1\})}\to \{0,1\}^{E(\mathbb{Z}^2)}$ as follows:
For $\omega\in \{0,1\}^{E(\mathbb{Z}^2\times \{0,1\})},$ and an edge $e=(v,w)\in E(\mathbb{Z}^2),$ we set $e$ to be open in $T\omega$ if either $e^0:=((v,0),(w,0))$ is open in $\omega$ or each of the three edges
$$
\{v^{\uparrow},e^1,w^{\uparrow}\}
$$
is open in $\omega$. Here, $e^1:=((v,1),(w,1))$ and $v^{\uparrow}=((v,0),(v,1))$.
    
    Now, for $\omega\sim \Prb_{p,\Z^2\times \{0,1\}},$ the distribution of $T\omega$ is not quite Bernoulli percolation on $\mathbb{Z}^2,$ since each vertical edge appears in multiple plaquettes.  Nonetheless, the upshot of the construction is that any two vertices $v,w$ of $\mathbb{Z}^2,$ $v$ and $w$ are connected in $T\omega$ only if $(v,0)$ and $(w,0)$ are connected in $\omega$. 

To get rid of the dependence, let $G$ be the multigraph obtained from $\mathbb{Z}^2\times \{0,1\}$ by replacing every vertical edge $v^{\uparrow}$ by four parallel edges $(v)^{\uparrow,e}$ indexed by the four non-vertical edges $e$ joining $(v,0)$ to $(v,1)$ (see Figure \ref{fig:Two_layered_surgery}). We consider the measure $\mathbb{P}^{\boldsymbol{p}}_{G}$ on $\{0,1\}^{E(G)}$ where each edge is open independently and with probability $p$ for all non-vertical edges and probability $\frac{p}{4}$ for the new vertical edges. For each $v\in \mathbb{Z}^2,$ by a union bound, the probability that some edge above $(v,0)$ is open is at most $p$. Therefore, for every pair of vertices $v,w\in V(G)=V(\mathbb{Z}^2\times \{0,1\}),$ we have
    \begin{equation} \label{Parallel_law}
    \Prb_{p,\Z^2\times \{0,1\}}[v\cc w]\geq \Prb^{\boldsymbol{p}}_{G}[v\cc w].
    \end{equation}
    Using this, we shall define a similar map $\tilde{T}:\{0,1\}^{E(G)}\to \{0,1\}^{E(\mathbb{Z}^2)}$ by declaring an edge $e=(v,w)$ open in $\tilde{T}\omega$ if either $e^0$ is open in $\omega$ or each of the three edges
    $$
    \{v^{\uparrow,e^1},e^1,w^{\uparrow,e^1}\}
    $$
    is open in $\omega$. Now, because of the splitting of the vertical edges, we get that $\tilde{T}\omega\sim \mathbb{P}_{\tilde{p},\mathbb{Z}^2}$
    for $\omega\sim \Prb^\mathbf{p}_{G}$. Inclusion-exclusion yields that
    $$\tilde p= p+p\left(\frac{p}{4}\right)^2-p^2\left(\frac{p}{4}\right)^2.
    $$
    For $p=\frac{1}{2}$ we find $\tilde p=\frac{1}{2}+\frac{1}{256}>\frac{1}{2}$, and by continuity, $\tilde p>\frac{1}{2}$ in a neighbourhood of $p=\frac{1}{2}$. 
    Accordingly, $\tilde{T}\omega$ percolates for all $p$ in this neighbourhood. As before, connections in $\Tilde{T}\omega$ imply connections in $\omega$, which, in turn, implies that $\omega$ percolates. By (\ref{Parallel_law}), so does $\Prb_{p,\Z^2\times \{0,1\}}$, which is what we wanted.
\end{proof}

\subsection{Furher remarks on UEG percolation}\label{sec:further remarks}

\subsubsection{Vertices of degree 4}
The proofs that $\UEG_{\Z^2}$ percolates (\Cref{Z2_percolates}) and that $\UEG_{\Z^d}$ percolates for $d\ge 3$ (\Cref{d3percolates}) are rather different, and neither of the proofs extends to cover the other case.
Let us take a moment to ponder the differences between the two methods.

First notice that one can sample the uniform even subgraph of $\Z^2$ by placing a fair coin on each plaquette, flipping them, and taking the symmetric difference of all the plaquettes where the coin landed heads up. As $\Z^2$ is self-dual, the distribution of the coins is site percolation on $\Z^2$.  Since site percolation on $\Z^2$ does not percolate at parameter $p = \frac{1}{2},$ the clusters of coins showing heads are all finite. Therefore, the infinite component of the UEG arises as a union of finite clusters of heads meeting at plaquettes that share only a vertex and not an edge. This is intimately related to the existence of vertices of degree $4$. Attempting to investigate the importance of vertices of degree 4, we prove that the uniform even subgraph of any bi-periodic trivalent planar graph does not percolate in \Cref{prop:trivalent_loop}. However, in \Cref{prop:jump_graph}, we construct a one-ended (non-amenable, non-planar) trivalent graph $\mathbb{J}$ such that ${\UEGop}_\mathbb{J}$ percolates.

\subsubsection{Odd percolation}\label{sec:odd percolation}
In parallel to the case of the uniform even subgraph, we may consider the uniform odd subgraph. Again, we consider spanning subgraphs but with the condition that each vertex degree is odd. The caveat is that some graphs do not admit any odd subgraphs, e.g. a cycle of odd length.
One criterion sufficient for the existence of odd subgraphs  of the graph $G$ is that it allows a dimerisation - that is, a perfect matching. Since $\Z^d$ allows a dimerisation, one way to define the uniform odd subgraph of $\Z^d$ is by taking the symmetric difference of the UEG and a fixed dimerisation. A more general characterisation of the uniform odd subgraph is that it is the unique probability measure on the co-set of odd subgraphs (supposing this is non-empty) which is invariant under the action of the even subgraphs.

Since Bernoulli-$\frac{1}{2}$ percolation is the Haar measure on the space of percolation configurations, it is invariant under taking the symmetric difference with any deterministic set. Therefore, the symmetric difference of a dimerisation with the edges in a hyperplane as in the proof of \Cref{d3percolates} is still $\Prb_{\frac{1}{2}}$ distributed. Thus, the proof of \Cref{d3percolates} generalises and the uniform odd subgraph percolates on $\Z^d$ for $d \geq 3$.

On the contrary, the proof of \Cref{Z2_percolates} does not generalise to the uniform odd subgraph of $\Z^2$ and to our knowledge, it is still open whether odd percolation has the same phase transition as even percolation in $\Z^2,$ the consequences of which are discussed in \cite[Section 6]{hansen2022strict}. Furthermore, we note in \Cref{prop:hexagonal_trivial} and \Cref{prop:hexagon_odd_no_perco} that neither the uniform even nor the uniform odd subgraph of the hexagonal lattice percolates. We summarise our knowledge of even and odd percolation in Table \ref{table:percolation_overview}. 

\subsubsection{Monotonicity in the ambient graph}
Recall from \Cref{thm:two_dimensions} that the uniform even subgraph of $\Z^2$ percolates. However, if we consider the inclusion $\Z^2 \subset \Z^3,$ then the marginal of $\UEGop_{\mathbb{Z}^3}$ on $\Z^2$ is $\Prb_{\frac{1}{2}, \Z^2}$ by Lemma \ref{lemma:restriction-separation}. 
So we find that  $\UEG_{\Z^3}\vert_{\Z^2}=\Prb_{\frac{1}{2}, \Z^2}$ does not percolate although $\UEG_{\Z^2}\vert_{\Z^2}=\UEG_{\Z^2}$ does percolate. This form of non-monotonicity of percolation with respect to inclusion is studied further in \cite[Corollary 1.2]{monster_paper} where we show that there exists a subgraph inclusion $\mathbb{G} \subset \mathbb{G}'$ such that $\UEG_{\mathbb{G}}$ percolates, but $\UEG_{\mathbb{G}'}$ does not.

\begin{table}
\begin{center}
\begin{tabular}{|c|c|c|} 
  \hline
   & even  & odd \\ 
    \hline
$\Z^2$  & \checkmark & ? \\ 
  \hline
$\Z^d, d \geq 3$ &  \checkmark  &  \checkmark  \\ 
  \hline
$\mathbb{H}$ & $\times$ & $\times$ \\ 
    \hline
\end{tabular}
\end{center}
\caption{Overview of percolation of even and odd percolation on the hypercubic and hexagonal lattices. \label{table:percolation_overview}}
\end{table}

\subsection{The loop \texorpdfstring{$\mathrm{O}$}{O}(1) model percolates for large \texorpdfstring{$x \in \lbrack 0,1 \rbrack$}{x in [0,1]}}

For monotone measures, given a single point of percolation, there is an interval of parameters for which there is percolation. 
This does not apply to the loop $\mathrm{O}$(1) model, which has negative association, but we can still bootstrap the strategy from \Cref{d3percolates} to get an interval of percolation points.

\nontrivphasetransition*

\begin{remark}
A similar statement for the single random current $\Prbcur_{\beta}$ follows from the fact that  $\Prbcur_{\beta}$ stochastically dominates $\ell_x$ cf. \Cref{thm:couplings}. 
However, this also follows without \Cref{thm:non_tri_phase_transition} since $\Prbcur$ dominates a Bernoulli percolation. 
\end{remark}

 \begin{proof}
For $d=2,$ the statement follows from the stronger result in Theorem  \ref{thm:two_dimensions}. 

The rest of the argument follows the strategy of Theorem \ref{d3percolates}. We are going to show that the marginal of $\ell_{x,\Z^d}$ on a suitably chosen set is bounded from below by a supercritical Bernoulli percolation. Again, we divide into cases according to whether $d=3$ or $d\geq 4$.

For $d \geq 4,$ consider a hyperplane $\Z^{d-1}\subset \Z^d$. Let $p\in (0,1)$ and consider the random-cluster model $\phi_{p,\Z^d}$. Then, by \Cref{thm:comparison theorem}, we have that $ \Prb_{\tilde p, \Z^d} \preceq \phi_{p,\Z^d}$ for $\tilde{p}=\frac{p}{2-p}$. Equivalently, there exists a coupling $(\omega, \tilde \omega)$ such that $\omega \sim  \phi_{p,\Z^d},$ $\tilde{\omega} \sim \Prb_{\tilde p, \Z^d}$ and $\tilde \omega \preceq  \omega$ almost surely.

Now, for any edge $e \in \Z^{d-1},$  we say that $e$ is \emph{good} if $e$ is open, the loop around the plaquette containing $e$ just above $e$ in $\Z^d$ is open, see Figure \ref{fig:plaquette_sticking_out}. 
The probability that this loop is open  in $\tilde{\omega} \sim \Prb_{\tilde p}$  is $\tilde p^4$. Define the process $Z_e \left( ( \omega, \tilde \omega) \right)  = \id_{ \{e \text{ is good}\}}(\tilde{\omega})$ for each $e \in E$. Then, $Z_e \perp Z_{e'} $ if $e$ and $e'$ do not share an endpoint. 
Therefore, by \cite[Theorem 0.0]{liggett1997domination}, the process $Z$ stochastically dominates some $W \sim \Prb_{q, \Z^d}$ where $q \to 1$ if $p \to 1$. 
Defining now $Q_e\left( ( \omega, \tilde \omega) \right) =  \id_{ \{e \text{ is good} \}}(\omega)$, then $Q_e \geq Z_e \geq W_e$  almost surely and therefore $Q \succeq W$. 

Next, we apply the relation \eqref{eq:def_infinite_loop_O} to get a coupling $(\omega,\eta)$ between $\omega\sim\phi_{p,\Z^d}$ and $\eta\sim\ell_{x,\Z^d}$  such that the conditional distribution of $\eta$ given $\omega$ is $\UEG_{\omega}$. Under this coupling, conditionally on $ \omega$, edges that are separated by cycles become independent by \Cref{lemma:restriction-separation}. Therefore, we must have that 
$\ell_{x,\Z^d}|_{\mathbb{Z}^{d-1}\times\{0\}} \succeq \mathbb{P}_{\frac{q}{2},\mathbb{Z}^{d-1}\times\{0\}}$.

Thus, for any fixed $a\in (0,\frac{1}{2})$ if we pick $p$ close enough to $1$ then $\frac{q}{2} > a$. Since $d-1\geq 3,$ we can once again use the fact that the edge percolation threshold for Bernoulli percolation on $\Z^3$ is strictly less than $\frac{1}{2}$ and so, $\ell_{x,\Z^d}|_{\mathbb{Z}^{d-1}\times\{0\}},$ and therefore also $\ell_{x,\Z^d}$, percolates.

For $d=3,$ we apply the same argument as in the proof of the $d=3$ in Theorem \ref{d3percolates}, where we replace the hyperplane by $\mathbb{Z}^2\times \{0,1\}$ and use \Cref{prop:three_complications}.
 \end{proof}

\subsection{Infinite volume constructions of uniform even graphs}\label{section:UEG measures}

Defining the UEG as a Haar measure immediately provides a construction of the uniform even subgraph of an infinite graph.
However, it is standard to regard percolation measures as limits of finitely supported measures \cite[Chap. 6]{friedli2017statistical} and this is, indeed, the approach of \cite{angel2021uniform,GJ09} for the uniform even subgraph. 
For the remainder of this section, we let $\mathbb G$ be a locally finite, countably infinite, connected graph. We show how both the limit of the free and wired measures may be realised as limits of Haar measures. The group of free even subgraphs is a limit of an increasing sequence of subgroups while the group of wired even subgraphs is determined by a decreasing sequence of quotients. The last part of this section ties our constructions to the notion of infinite volume Gibbs measures and provides a characterisation of these.

\subsubsection{The free and the wired \UEG}\label{section:limit-constructions}

To begin with, let us briefly review a standard approach.
It was observed in \cite{GJ09} that any locally finite connected graph admits a finitary generating set 
$\mathcal C\subset \Even(\mathbb{G})$\footnote{If $\mathbb{G}$ is finite and connected, and $T$ is a spanning tree in $\mathbb{G}$, a map $\mathbb{E}\setminus T\to \Even (\mathbb{G}),$ the image of which is a basis,  can be constructed as follows: The end-points of a given edge $e\in\mathbb{E}\setminus T$ are connected through a unique path in $T$ which together with $e$ defines a cycle $C_e$, an element in the basis. Since each $e\in \mathbb{E}\setminus T$ appears in exactly one such cycle, the set $\{C_e\}$, indeed, linearly independent. One then checks that any cycle has coefficients in the basis given by its intersection with $\mathbb{E}\setminus T$.}.
A set $\mathcal F\subset \Even(G)$ is finitary if any given edge is contained in finitely many $f \in \mathcal F$.
Since $\mathcal{C}$ is finitary, the sums 
\begin{equation}\label{eq:sample-cycles}
    \sum_{C\in \mathcal A}C  = \sum_{C\in \mathcal C}\id_{\mathcal A}(C)C,
\end{equation}
are pointwise convergent for any $\mathcal{A\subset C}$. 
The set is generating if $\Even(\mathbb{G})=\overline{\text{span}(\mathcal C)}$.
The set $\mathcal C$ can be chosen as a Schauder basis of $\Even$ and properties of such sets are studied in \cite{angel2021uniform}. The observation is used for sampling a uniform even graph by (uniformly) sampling the subsets $\mathcal{A\subset C}$. This is done by replacing the coefficients $\id_{\mathcal A}(C)$ in the sum \eqref{eq:sample-cycles} with independent random Bernoulli-1/2 variables $\epsilon_C$. 
Observe, that the construction of this measure amounts to defining a surjective homomorphism $\Phi^\mathcal{C}:\Z_2^\N\to\Even(\mathbb{G})$ (by enumerating $\mathcal C$) and applying \Cref{pushhom} to obtain the Haar measure on $\Even(\mathbb{G})$ as $\Phi^\mathcal C_*\mu$, where $\mu$ is just a product measure of Bernoulli-1/2 measures.
One speaks of sampling the coefficients $\epsilon_C$ sequentially, implying that one is really taking the weak limit of the measures $(\Phi^{\mathcal{C}}|_{\N_{\le N}})_*\mu$ for increasing $N$.
These measures on $\Omega(\mathbb{G})$ are supported on finitely many elements of $\Even(G)$ and approximate $\Phi^\mathcal C_*\mu$. Indeed, since $\mathcal{C}$ is finitary, the expectation of any local event (or so-called cylinder event on a finite cylinder) is eventually constant as $N\to\infty$.

With $\UEG$ defined a priori as the Haar measure (and not as a limit or as $\Phi^\mathcal C_*\mu$ for specific $\mathcal C$), we are led to consider approximations of $\UEG$ by finitely supported measures.
In addition, we may ask whether the approximations are \emph{local} in the sense that it agrees with the pushforward of a measure $\mu$ on $\Omega(G)$ along the natural inclusion $\Omega(G)\subset \Omega(\mathbb{G})$  for a finite set $E$ and $G=(\mathbb V,E)$.
If the finitary basis $\mathcal C$ consists of finite cycles, we may consider $\Phi^\mathcal C_*\mu$ to be a limit of local measures.
In \cite{angel2021uniform}, it is recognised that when $\mathcal C$ is a subset of the finite subgraphs, it can at most  generate the set of free uniform graphs $\Omega^0_\emptyset(\mathbb{G})\subset \Even(\mathbb{G})$ (see eq. \eqref{eq:def-free}), which is a proper subset of $\Even(\mathbb{G})$ if $\mathbb{G}$ has more than one end (a fact which will be discussed at the end of this section, no pun intended). 

A natural approach to locally approximating $\UEG$ is to consider $\UEG_{G}$ for all $G=(\mathbb V,E)\subset \mathbb{G}$ where $E$ is finite.
The simple observation that $\Even(G)\subset \Even(G')$ whenever $G\subset G'$ shows that these groups form an increasing net ordered by inclusion of finite subgraphs\footnote{The reader may choose to think of a sequence of finite subgraphs $G_1\subset G_2\subset\ldots$ such that $\mathbb{G}=\bigcup_{n\in\N}G_n$ rather than the net of all finite subgraphs.}. Consequently, the Haar measures $\UEG_{G}$ converge due to a general result.

\begin{theorem} \label{thm:Haar_first} 
    Let $H$ be a compact Abelian group and let $(\Gamma_\alpha)_{\alpha\in\mathcal{I}}$ be an increasing net of closed sub-groups of $H$.  The Haar measures $\mu_\alpha$, on $\Gamma_\alpha$ respectively, converge weakly to the Haar measure on $\Gamma=\overline{\cup_{\alpha\in\mathcal{I}} \Gamma_{\alpha}}$.
\end{theorem}
\begin{proof}
    The measures $\mu_\alpha$ extend to measures on $\Gamma$ by push-forward under the inclusion $\Gamma_\alpha\subset \Gamma$. 
    By compactness of the space of probability measures on $H$, it suffices to establish that all accumulation points of $(\mu_\alpha)$ agree. Thus, let $\nu$ be an accumulation point of $\mu_\alpha$.    Since $\Gamma_\alpha\leq \Gamma_{\beta}$ for all $\alpha\preceq \beta$, $\nu$ is invariant under translation by any element of $\cup_{\alpha\in\mathcal{I}}\Gamma_\alpha$.

    Now let $(g_j)_{j\in\mathcal{J}}\subset \cup_{\alpha\in\mathcal{I}}\Gamma_\alpha$ be a net converging to $g\in \Gamma$. 
    Let $f$ be a continuous function on $\Gamma$, and note that
 the uniform continuity of $f$ implies that $f(x-g_j)\to f(x-g)$ uniformly. Therefore,
    $$
    \lim_{j\in\mathcal{J}}\nu[f(x-g_j)]=\nu[f(x-g)].
    $$
    Since $\nu[f(x-g_j)]=\nu[f]$ for all $j$, we conclude that $\nu$ is a probability measure with support on $\Gamma$ which is invariant under translation by elements of $\Gamma$. Therefore, $\nu$ is the Haar probability measure on $\Gamma$.
 \end{proof}

Define the set of free even subgraphs of $\mathbb{G}$
\begin{equation}\label{eq:def-free}
    \Omega_\emptyset^0(\mathbb{G})
    =\overline{\bigcup_{G\subset \mathbb{G} \text{ finite}}\Even(G)}
    =\overline{\bigcup_{n\in \N}\Even(G_n)},
\end{equation}
where $(G_n)_{n\in\N}$ is any sequence of finite subgraphs, $G_n\Uparrow\mathbb{G}$. The Haar measure on $\Even^0(\mathbb{G})$, denoted $\UEG^0_{\mathbb{G}}$, is called the free uniform even subgraph. \Cref{thm:Haar_first} shows that $\UEG^0_{\mathbb{G}}=\lim_{}\UEG_{G_n}$ weakly.
Furthermore, the existence of a finitary basis implies that this approximation is eventually constant on local events.

On the other hand, we refer to $\Even(\mathbb{G})$ as the set of \emph{wired} even subgraphs \cite{angel2021uniform}. We shall now see how $\Even(\mathbb G)$ is approximated by local wired even subgraphs, giving some justification to this name.
Recall that  the set of wired even subgraphs $\Wired(G)$ for a subgraph  $G\subset \mathbb{G}$ is defined in terms of the boundary $\partial_vG$ with respect to $\mathbb G$. 
If $G\subset G'\subset\mathbb G$, typically $\Wired(G)\nsubseteq\Wired(G')$ in contrast to the situation for $\Even$.
Instead, the arrows are reversed. That is, $\Wired(G)$ is (or at least contains) a quotient of $\Wired(G')$.
Recall that by the isomorphism theorem, finite quotients correspond to homomorphisms with finite range such as the projections $\pi_{G}|_{\Even(\mathbb{G})}$ onto $\pi_{G}({\Even(\mathbb{G})})$ in our setting.
Indeed, $\Even(\mathbb{G})$ is profinite, which is a way of saying that it is determined (up to isomorphism) by its (category of) finite quotients.
In particular, write
\begin{equation*}
    \Even(\mathbb{G})=\varprojlim\pi_{G}(\Even(\mathbb{G}))
\end{equation*}
for the cofiltered projective limit over the net of finite subgraphs ordered by inclusion. We do not unravel the definition here but remark that it is the smallest object admitting suitable projections and refer the reader to \cite{Rao_1971} for a presentation within the category of measure spaces.
It was noted in \cite[Theorem 2.6]{GJ09} that the projections $\pi_{G}|_{\Even(\mathbb{G})}$ with corresponding marginals $\UEG_{\mathbb{G}} \vert_{G}$ for all finite $G\subset\mathbb G$ determine $\Phi^\mathcal{C}_*\mu$ by Kolmogorov's extension theorem. 
Again, this is a consequence of a general result for Haar measures.

\begin{theorem}\label{limit-profinite}
    Let $\Gamma=\varprojlim \Gamma_\alpha$ be a profinite group.
    Then, the projective limit of the normalised Haar measures on $\Gamma_\alpha$ exists and identifies with the consequently unique normalised Haar measure on $\Gamma$.
\end{theorem}
\begin{proof}
    Let $\mu_\alpha$ denote the Haar probability measure on $\Gamma_\alpha$, and $\pi_\alpha:\Gamma\to\Gamma_\alpha$ the projection for all $\alpha$. The conditions of \cite[Theorem 3.4]{Rao_1971} are trivially satisfied for finite sets $\Gamma_\alpha$ so the existence of a unique regular Borel probability measure $\mu$ on $\Gamma$ is granted, such that $\mu_\alpha=(\pi_\alpha)_*\mu$ for all $\alpha$ and satisfying inner regularity with respect to cylinder sets, which amounts to 
    \begin{equation}\label{eq:proj-convergence}
        \mu[U] = \lim_{\alpha}\mu_\alpha[\pi_\alpha(U)],
    \end{equation}
    for $U\subset \Gamma$ open.
    The left and right invariance of $\mu$ follows from that of $\mu_\alpha$ and \eqref{eq:proj-convergence} since $\pi_\alpha$ is a homomorphism for each $\alpha,$ so 
    $\mu_\alpha(\pi_\alpha(gUh))=\mu_\alpha(\pi_\alpha(g)\pi_\alpha(U)\pi_\alpha(h))=\mu_\alpha(U)$ for $g,h\in\Gamma$.
    Thus, $\mu$ is a Haar measure on $\Gamma$.
    Any Haar measure on $\Gamma$ shares the properties of $\mu$, so the uniqueness of $\mu$ according to \cite[Theorem 3.4]{Rao_1971} implies the uniqueness of the Haar measure.
\end{proof}
\begin{remark}\label{remark:Haar-Kolmogorov}
    There are several approaches to constructing Haar measures, and \Cref{limit-profinite} with the groundwork in \cite[Theorem 3.4]{Rao_1971} is among them. This result is not new but included for completeness of the presentation.
    Note also that \cite[Theorem 3.4]{Rao_1971} generalises Kolmogorov's extension theorem. In particular, Kolmogorov's extension theorem provides a construction of $\mathbb P_{\frac{1}{2}},$ which coincides with the Haar measure on $\Omega(\mathbb G)=\varprojlim \Omega(G)$ by \Cref{limit-profinite}.
\end{remark}

Returning to the setting $G\subset \mathbb G$, observe that since $\pi_{G}(\Even(\mathbb{G}))\not \subset \Even(\mathbb{G}),$ the uniform measure on 
$\pi_{G}(\Even(\mathbb{G}))$
does not push forward to a measure on $\Even(\mathbb{G}),$ but it does push forward to $\Omega(\mathbb{G})$ along the inclusion 
$\pi_{G}(\Even(\mathbb{G}))
\subset \Omega(\mathbb{G})$. 
With this in mind, the convergence in \eqref{eq:proj-convergence} can be realised as weak convergence of measures on $\Omega(\mathbb{G})$.
Therefore, we obtain a local approximation of $\UEG_{\mathbb{G}}$ by uniform measures on the groups $\pi_{G}(\Even(\mathbb{G}))$. As it turns out, on a connected graph, we may take approximating measures to be wired measures, $\UEG_{G}^1$.
Indeed, if $G_n\Uparrow \mathbb{G}$, such that each $G_n$ is finite and every component of $(\mathbb {V,E}\setminus E_n)^\circ$ is infinite, then it follows from \Cref{lemma:restriction-separation} that $\UEG_{\mathbb{G}}$ is approximated by a sequence of local wired uniform even subgraphs.
We summarise the discussion in a theorem:
\begin{theorem}
    Let $\mathbb G$ be a locally finite, infinite connected graph. The Haar measures on $\Even(\mathbb{G})$ and $\Even^0(\mathbb{G})$ as probability measures on $\Omega(\mathbb{G})$ are weak limits of local wired and free uniform even subgraphs respectively.
\end{theorem}

Finally to compare the free and the wired even subgraphs,
let $G_n$ be as before and consider the commuting diagram consisting of inclusion and restriction maps
\begin{equation*}
\begin{tikzcd}
\dots \arrow[r] & \Omega_{\emptyset}(G_n) \arrow[d] \arrow[r] & \Omega_{\emptyset}(G_{n+1}) \arrow[d] \arrow[r] & \dots \arrow[r] & \Omega_{\emptyset}^0(\mathbb{G}) \arrow[d]  \\
\dots           & \Wired(G_n) \arrow[l]          & \Wired(G_{n+1}) \arrow[l]          & \dots \arrow[l] & \Omega_{\emptyset}(\mathbb{G}) \arrow[l]
\end{tikzcd}
\end{equation*}
This diagram characterises the limiting groups. 
In particular, the comparison map between them is the inclusion giving rise to the quotient $\Even(\mathbb{G})/\Even^0(\mathbb{G})$ which characterises infinite volume Gibbs measures of the $\UEG$ as shown in \Cref{thm:Gibbs of UEG} below.

\subsubsection{Gibbs measures of $\UEG$}\label{sec:Gibbs of UEG}

A percolation measure $\mu$ on $\Omega(\mathbb{G})$ is said to be \emph{Gibbs} for the $\UEG$ if for every finite graph $G\subset \mathbb{G}$, if $\xi_G$ denotes the random variable $\eta|_{\mathbb{G}\setminus G}$, we have 
$$\mu[ \eta|_G=\eta_0 |\; \xi_G]\propto \id_{\partial \eta_0 =\partial \xi_G} \hspace{1cm} \text{         $\mu$-a.s.}
$$ 
 Note that any Gibbs measure for the $\UEG$ is supported on $\Even(\mathbb{G})$. We call a Gibbs measure extremal if it cannot be written as a non-trivial convex combination of distinct Gibbs measures.
 There is a binary operation of measures given by the pushforward of the product measure by symmetric difference. This is referred to as symmetric difference of percolation measures.
\begin{theorem}\label{thm:Gibbs of UEG} 
    For any infinite, locally finite graph $\mathbb{G}$, 
    the set of extremal Gibbs measures of the $\UEGop$ on $\mathbb{G}$ forms a group under symmetric difference isomorphic to 
    $\Even(\mathbb{G})/\Even^0(\mathbb{G})$.
\end{theorem}
\begin{proof}
Fix $\eta_H\in H\in \Even(\mathbb{G})/\Even^0(\mathbb{G})$ and note that the map $\eta\mapsto \eta \triangle \eta_H $ pushes $\UEG^0_{\mathbb{G}}$ forward to an $\Even^0(\mathbb{G})$-invariant measure on $H$. We denote this measure by $\UEG^H$ as it does not depend on the choice of $\eta_H$.
There can only be one measure on $H$ which is $\Even^0(\mathbb{G})$-invariant, since pushing forward along $\eta\mapsto\eta\triangle\eta_H$ again maps any such measure to the unique Haar measure on $\Even(\mathbb{G})$ (and this map is an involution). The measure $\UEG^H$ may be seen to be Gibbs.

Conversely, any Gibbs measure of the $\UEG$ is invariant with respect to the action of $\Even^0(\mathbb{G})$, since it is invariant under the action of the group of finite even subgraphs. We saw this in \Cref{thm:Haar_first}.

By ordinary theory (see e.g. \cite[Theorem 6.58]{friedli2017statistical} ), if $A$ is a tail-event and $\mu$ is any extremal Gibbs measure, then $\mu[A]\in \{0,1\}$. For a given co-set $H  \in \Even(\mathbb{G})/\Even^0(\mathbb{G})$, let $H^{\infty}=H\triangle \Omega^0(\mathbb{G}),$ the element-wise sum, where $\Omega^0(\mathbb{G})$ is the space of finite subgraphs of $\mathbb{G}$. This is the smallest tail-event on $\Omega(\mathbb{G})$ which contains $H$. Note that if $H\neq H'$ are two different co-sets, then $H^{\infty}\cap (H')^{\infty}=\emptyset$. As $\mu$ is invariant under the action of $\Even^0(\mathbb{G})$, it is supported on at least one $H^{\infty}$. By the above, $H$ is determined uniquely and $\mu=\UEG^H$. Accordingly, at least one $\UEG^H$ is extremal. Since push-forwards under the map $\eta\mapsto \eta\triangle \eta_H$ preserve convex combinations of measures, we get that every $\UEG^H$ is extremal.

Finally, we observe by the convolution property of Haar measure, that $\UEG^H\triangle\UEG^{H'}$ is the same as the pushforward of $\UEG^{H'}$ along $\eta\mapsto\eta\triangle\eta_H$. 
\end{proof}

Since the set of all Gibbs measures is a simplex spanned by the extremal Gibbs measures \Cref{thm:Gibbs of UEG} characterizes the Gibbs measures of UEG. 
 
\subsubsection{Evens and ends}
We end this section by giving some results on $\Even(\mathbb{G})/\Even^0(\mathbb{G})$ which, by \Cref{thm:Gibbs of UEG}, characterises the extremal Gibbs measures of $\UEG$.


The intuition, which will be made precise in this section, is as follows. Loosely speaking, for a path in $\mathbb{G}$ to be an even subgraph, it must either keep turning back on itself or run off to infinity in either direction. The ways to run off to infinity are, in some sense, the $\mathbb{G}$-ends, so $\mathbb{G}$-ends can be thought of as an infinite volume analogue of boundary vertices that may act as sources.

A concrete way of thinking about this is to consider trees.
Observe that a connected graph $T$ is a tree if and only if $\Even^0(T)$ is the trivial group. In fact, $\Even(\mathbb{G})/\Even^0(\mathbb{G})$ may be probed by finding subtrees $T\subset \mathbb{G}$ such that the inclusion map descends to a well defined map,
$$\Even(T)\to\Even(\mathbb{G})/\Even^0(\mathbb{G}),$$
and maximal trees with this property give isomorphisms. We will not pursue this point any further than to motivate thinking of $\Even(\mathbb{G})/\Even^0(\mathbb{G})$ as a measure of the coarse tree-like structure of $\mathbb{G}.$

The following discussion is somewhat analogous to Häggströms discussion of the random-cluster measure on the tree in \cite{Haggstrom1996}.  

Recall that an \emph{end} in an infinite graph $\mathbb{G}$ is an equivalence class of \emph{rays}\footnote{a \emph{ray} is a vertex self-avoiding path with a single end-point.}, where two rays are equivalent if there exists another self-avoiding path intersecting both rays infinitely often.
Equivalently, two rays are equivalent if and only if, whenever a finite subgraph is excised from $\mathbb{G}$, 
the remaining infinite connected components  of the two rays (one for each) belong to the same (infinite) connected component. 
Consequently, the set of ends may be defined as the limit of the set of infinite connected components of $\mathbb{G}\setminus G$ as $G$ becomes a larger and larger finite subgraph. Formally,
$$\mathfrak{e}(\mathbb{G})= \varprojlim\mathcal{C}^\infty(\mathbb{G}\setminus G),$$
where the projective limit\footnote{We only consider this as a limit of sets and disregard the natural topology of the limit as a totally disconnected space.} is taken over all finite $G\subset \mathbb G$, and $\mathcal{C}^\infty(\mathbb{G}\setminus G)$ denotes the set of infinite components of the graph $\mathbb{G}\setminus G$, and $\mathcal{C}_x(\mathbb G\setminus G')$ is mapped to $\mathcal{C}_x(\mathbb G\setminus G)$ for all $x\in \mathbb V$ when $G\subset G'$.

Specifically,for each $G$ there is a map $\mathfrak{e}(\mathbb{G})\to\mathcal{C}^\infty(\mathbb{G}\setminus G)$ , which sends an end $\epsilon\in \mathfrak{e}(\mathbb{G})$ to the infinite component $\epsilon_G$ of $\mathbb{G}\setminus G$ where it lies.

We may note that two ends $\epsilon,\epsilon'$ are distinct precisely if there exists a finite $G$ such that $\epsilon_G\ne\epsilon'_G$. However, an end may not always be determined in finite volume as such. 
We say an end $\epsilon$ is \emph{finitely separated}
if there is a finite $G$ such that $\epsilon'_G=\epsilon_G$  implies $\epsilon'=\epsilon$. It is instructive to consider the $d$-regular tree, which has two finitely separated ends when $d=2$ and none when $d>2$. 

To understand whether an even subgraph uses an end as a 'source at infinity,' so to speak, we consider its sources when restricted to $\epsilon_G$. These live on the boundary $\partial_vG$ which is finite.
Let $|A|_2$ denote the cardinality of a finite set mod $2$, equivalently viewed as a finite sum in $\Z_2$.
Define the \emph{end map} $\partial_\mathfrak{e} : \Even(\mathbb{G}) \to \{0,1\}^{\mathfrak{e}(\mathbb{G})}$,
\begin{equation}\label{eq:end map}
    \partial_\mathfrak{e}\eta(\epsilon) 
    = \liminf_G|\partial(\eta|_{\epsilon_G})|_2,
\end{equation}
where $\liminf$ is taken along finite subgraphs ordered by inclusion. 
There are examples of subgraphs that intersect an end infinitely often without containing the end (referring to the natural identification $\mathfrak{e}(H)\subset \mathfrak{e}(\mathbb{G})$ when $H\subset \mathbb{G}$): in the trivalent tree, there is an even subgraph containing every second edge of some fixed ray. That goes to show that $\liminf$ cannot be replaced by $\lim$ in the definition of $\partial_\mathfrak e$. This pathological behaviour also prevents $\partial_\mathfrak{e}$ from being a homomorphism in general.

\begin{lemma}
    \label{lem:end map hom on lim}
    For all $\epsilon\in\mathfrak{e}(\mathbb{G})$ and $\eta_1,\eta_2\in\Even(\mathbb{G})$ satisfying that $(|\partial(\eta_i|_{\epsilon_G})|_2)_G$ is a convergent net, $i=1,2$, then also $(|\partial(\eta_1\triangle\eta_2|_{\epsilon_G})|_2)_G$ is convergent with limit
    $$\partial_\mathfrak{e}(\eta\triangle\eta')(\epsilon)=
    \partial_\mathfrak{e}\eta(\epsilon)+\partial_\mathfrak{e}\eta'(\epsilon)\mod 2.$$
\end{lemma}
\begin{proof}
    First note that $\eta\mapsto |\partial(\eta|_{\epsilon_G})|_2$ is a homomorphism for every $G,\epsilon$. 
    Let $\epsilon\in\mathfrak{e}(\mathbb{G})$, and assume $\eta_1,\eta_2$ satisfy the hypothesis.
    Since the net defining $\partial_\mathfrak{e}\eta_i(\epsilon)$ is convergent and $\{0,1\}$-valued for $i=1,2$, there are finite $G_i$ such that $|\partial(\eta_i|_{\epsilon_G})|_2$ is constant for $G\supset G_i$. It follows from the homomorphism property that  $\partial(\eta\triangle\eta')|_{\epsilon_G}$ is constant valued $\partial_\mathfrak{e}\eta_1(\epsilon)+\partial_\mathfrak{e}\eta_2(\epsilon)\mod 2$ for $G\supset G_1\cup G_2$.
\end{proof}
We retain the homomorphism property of $\partial_\mathfrak{e}$ in certain cases. 
Define
\begin{equation*}
    \begin{gathered}
    \mathfrak{E}^c=\{f\in \{0,1\}^{\mathfrak{e(\mathbb{G})}}|\, f\text{ has finite support},\enspace 
   |f|_2=0\}, \\
\Even^c(\mathbb{G})=\partial_\mathfrak{e}^{-1}(\mathfrak{E}^c), \\
 \mathfrak{e}^s(\mathbb{G})=\{\epsilon \in \mathfrak{e}(\mathbb{G})|\, \epsilon \text{ is finitely separated}\}, \\
 \pi_s : \{0,1\}^{\mathfrak{e(\mathbb{G})}} \to \{0,1\}^{\mathfrak{e}_s(\mathbb{G})},
\end{gathered}
\end{equation*}
with $\pi_s$ being the restriction map.
\Cref{prop:evens and ends} below extends a result of \cite{angel2021uniform} which states that the free and wired uniform even subgraphs coincide if $\mathbb{G}$ is one-ended. 

\begin{proposition}\label{prop:evens and ends}
    We have 
    \begin{equation}\label{eq:kernel_characterization}
        \Even^0(\mathbb{G}) =
        \{\eta\in \Even(\mathbb{G})|\, \partial_\mathfrak{e}\eta=0\}. 
    \end{equation}
    The end map restricts to a surjective homomorphism $\partial_\mathfrak{e} : \Even^c(\mathbb{G}) \to \mathfrak{E}^c$.
    The composition $\pi_s\circ\partial_\mathfrak{e} : \Even(\mathbb{G}) \to \{0,1\}^{\mathfrak{e}_s(\mathbb{G})}$  is a continuous homomorphism. 
    
\end{proposition}
\begin{proof}
    To prove \eqref{eq:kernel_characterization}, let $G_n$ be an exhausting sequence of finite subgraphs.  Let $\eta\in \Even^0(\mathbb{G})$ and let $(\eta_n)$ be a sequence of finite even subgraphs in $\mathbb{G}$ such that $\eta|_{G_n}=\eta_n|_{G_n}$.
    Since $\eta_n$ is finite, then $\partial_\mathfrak{e}\eta_n = 0$,\footnote{Note that $\partial_{\mathfrak{e}}$ is not continuous globally, so we cannot immediately take the limit $n\to\infty$ here.} and in fact $|\partial\eta_n|_{\epsilon_G}|_2=0$ for any $\epsilon, G$. Given $\epsilon$ and finite $G$, there is suitably large $n$ such that $\partial\eta|_{\epsilon_G}=\partial\eta_n|_{\epsilon_G}$. As this vanishes, $\partial_\mathfrak{e}\eta(\epsilon)=0$.

    Conversely, assume $\eta\in\Even(\mathbb{G})$ and $\partial_\mathfrak{e}\eta=0$. We construct a sequence of finite even subgraphs $\eta_n$ such that $\eta_n|_{G_n}=\eta|_{G_n}$ as follows. 
    For given $n$, since $\partial_\mathfrak{e}\eta(\epsilon)=0,$ it is possible, for any $\epsilon\in\mathfrak{e}(\mathbb{G})$, to pick $\tilde\eta_{\epsilon,n}\subset \mathbb{E}\setminus E_n$ finite connecting the sources of $\eta|_{G_n}$ that belong to $\epsilon_{G_n}$, that is,
    \begin{equation*}
        \partial(\eta|_{G_n})|_{\epsilon_{G_n}} =\partial\tilde\eta_{\epsilon,n}.
    \end{equation*} 
    Since the set $\{\epsilon_{G_n}| \epsilon\in \mathfrak{e}(\mathbb{G})\}$ is finite, it is only necessary to pick $\tilde\eta_{\epsilon,n}$ for one $\epsilon$ representing each of these finitely many infinite components. 
    Define $\eta_n$ by the finite sum over representatives $\eta_n = \eta|_{G_n} \triangle_\epsilon \tilde\eta_{\epsilon,n}$.

    Now consider the restriction to $\Even^c$.
    By \Cref{lem:end map hom on lim}, it suffices to show that $\liminf$ can be replaced by $\lim$ in \eqref{eq:end map} for all $\epsilon\in\mathfrak{e}(\mathbb{G})$ and all $\eta\in\Even^c(\mathbb{G})$, in order to prove the homomorphism property.

    Let $\eta\in \Omega^c_\emptyset(\mathbb{G})$ and $\delta\in\mathfrak{e}(\mathbb{G})$. If $\partial_{\mathfrak{e}}\eta(\delta)=1$ there is nothing to prove since $\limsup_G |\partial (\eta|_{\epsilon_G})|_2\leq 1.$
    So assume $\partial_\mathfrak{e}\eta(\delta)=0$. 
    We show that $\limsup_G|\partial(\eta|_{\delta_G})|_2=0$. 
    
    There is finite $G'$ such that for all $G\supset G'$, $\delta_{G}\ne\epsilon_{G}$ for every $\epsilon$ in the support of $\partial_\mathfrak{e}\eta$ since there are only finitely many such $\epsilon$.

    There is also finite $G''$ such that for all $\epsilon, \epsilon'$ in the (finite) support of $\partial_\mathfrak{e}\eta$, we have $\epsilon_{G''} = \epsilon'_{G''}$ if and only if $\epsilon=\epsilon'$ (i.e. $G''$ literally separates the ends in supp($\partial_\mathfrak{e}$)).
    From here on out, i.e. $G\supset G''$, we know for any $\epsilon'$, that $|\partial\eta|_{\epsilon'_{G}}|_2=1$ if and only if $\epsilon'_{G}=\epsilon_{G}$ for some $\epsilon$ in the support. 
    
    So for all finite subgraphs $G\supset G'\cup G''$, we conclude that $|\partial\eta|_{\delta_{G}}|_2=0$, hence 
    \begin{equation}
        \limsup_G |\partial\eta|_{\delta_{G}}|_2=0.
    \end{equation}
    
    The range of $\partial_\mathfrak{e}$ when restricted to $\Even^c(\mathbb{G})$ can be understood as follows:
    For any pair of ends, there is a simple bi-infinite path $\gamma$ with $\partial_\mathfrak{e}\gamma = \{\epsilon_1,\epsilon_2\}$. As $\gamma\in \Even^c(\mathbb{G})$ we can take symmetric differences of such paths to produce $f\in \mathfrak{E}^c$.\footnote{Note that the reasoning above can be applied to see that $|\partial_\mathfrak{e}\eta|_2=0$ is automatically satisfied if $\partial_\mathfrak{e}\eta$ is finitely supported.}

    Turning to the final statement, observe that
    whenever $G'\subset G$ are finite subgraphs, then
    \begin{equation}\label{eq:sum ends}
        |\partial(\eta|_{\epsilon_{G'}})| = \sum_{\substack{c\in \mathcal{C}^\infty(\mathbb{G}\setminus G) \\ c\subset \epsilon_{G'}}}
        |\partial(\eta|_c)| \mod 2,
    \end{equation} 
    for all $\eta\in\Even(\mathbb{G})$ and $\epsilon\in\mathfrak{e}(\mathbb{G})$. This follows by the usual parity considerations over sources of finite graphs.
    When $\epsilon$ is finitely separated by finite $G'$, then there is precisely one summand in \eqref{eq:sum ends}, and hence $\partial_\mathfrak{e}\eta(\epsilon)=|\partial(\eta|_{\epsilon_{G'}})|_2$.
    This implies in particular that $\liminf$ is $\lim$ for all such $\epsilon$ and all $\eta\in\Even(\mathbb{G})$. 
    Again, \Cref{lem:end map hom on lim} ensures the homomorphism property upon restriction to finitely separated ends.
    
    If $\eta_n\to \eta$ is a pointwise converging sequence in $\Even(\mathbb{G})$, and $\epsilon$ is finitely separated by $G$, then we can find $N$ such that for all $n>N$, $\eta_n|_{G}=\eta|_G$. Then
    \begin{equation}
        \partial_\mathfrak{e}\eta_n(\epsilon) = |\partial(\eta_n|_{\epsilon_{G}})|_2=
        |\partial(\eta|_{\epsilon_G})|_2 = \partial_\mathfrak{e}\eta(\epsilon).
    \end{equation}
    So $\partial_\mathfrak{e}\eta_n\to\partial_\mathfrak{e}\eta$ pointwise, which proves continuity upon restricting to finitely separated ends.
\end{proof}

In light of \Cref{thm:Gibbs of UEG}, the following characterisation of Gibbs measures of the $\UEG$ can be obtained for graphs with finitely many ends. The general case is left to further work. 
\begin{corollary}
     Let $\mathbb G$ be a locally finite, infinite connected graph with finitely many ends. Then the group of extremal Gibbs measures of the $\UEGop$ is isomorphic to 
     $$\{f\in \{0,1\}^{\mathfrak{e(\mathbb{G})}}|\,
   |f|_2=0\}.$$ 
\end{corollary}

Finally, some remarks are due to justify the restrictions in the formulation of \Cref{prop:evens and ends}.
The end map $\partial_\mathfrak{e}$ may not be surjective on all of $\Even(\mathbb{G})$. For example, on the trivalent tree, $T_3$, any even subgraph $\eta$ is a countable union of simple infinite paths, so $\partial_\mathfrak{e}\eta$ is supported on at most countably many ends. However, the tree has uncountably many ends.

The $4$-regular tree, $T_4,$ similarly has uncountably many ends but is in itself an even graph. Thus, we find that any countable or co-countable end configuration is obtained in the range (up to parity constraint). It is not clear whether there are $\eta\in \Even(T_4)$ which have an uncountable number of open ends and an uncountable number of closed ends.
An example of this is found in the $5$-regular tree $T_5$ by considering  an isomorphic copy of $T_4\in \Even(T_5)$.

The Dirac comb is the graph $\mathbb{D}$ obtained from a copy of $\N$ attached to every vertex of $\Z$. The two rays of $\Z$ are not finitely separated but every other end, belonging to a copy of $\N$, is. 
There are exactly two $\eta,\eta' \in \Even(\mathbb{D})$ with all finitely separated ends open. 
They satisfy $\eta\triangle\eta'=\eta_\Z$, the subgraph corresponding to the embedding of $\Z$. 
This shows that the end map does not give a complete invariant of $\Even(\mathbb{D})$.
In this example, we also see that $\partial_\mathfrak{e}$ is only additive on the finitely separated ends, explaining why the kernel of $\partial_\mathfrak{e}$ does not measure this failure of injectivity.

One can construct $\eta_\mathbb{Z}$ as the limit of $\eta_n\in \Even(\mathbb{D})$, where $\eta_n|_\mathbb{Z}=[-n,n]$. 
This gives an example of a discontinuity of the end map, which is observed by evaluating at a non-finitely separated end.

\section{Proof of Theorem \ref{main theorem}}
\label{sec:Main arguments}

Before we get our hands dirty with the details of the proof of our main theorem, we begin this section by giving a brief overview of the arguments that go into the proof. We do so to emphasise the underlying topological ideas, for the benefit of the impatient reader, as well as preparing the road ahead.
Then, we review some technical aspects of the Ising literature, and set up the machinery that we shall need. Together with the previous results on the marginals of the uniform even graph, this enables the proof of \Cref{thm:Loop O(1) mixing} as well as \Cref{thm:Torus wrap}. The latter requires an adaptation of the multi-valued mapping principle, to be discussed.

\subsection{Road map and torus basics}\label{Road map}
 We saw in Corollary \ref{Surprising fact} that for $\mathbb{Z}^d$, the boundary conditions under which we take the uniform even subgraph do not matter as soon as we take a single step away from the boundary due to the richness of loops in the graph. For the supercritical random-cluster model, we can hardly expect a result quite as powerful, but it is known that we have a local uniqueness property of the model ensuring that any 
 suitably nice finite piece of the lattice is going to have a single large cluster. Thus, with the caveat that we might have to take a bit more than a single step away from the boundary, we should also expect the boundary conditions we are working with not to matter. In the rest of the section, we fix $d\geq 3$.

At the same time, even graphs have some nice interplay with the topology of the torus. 

For 
$H$  a hyperplane  in $\mathbb{T}_n^d$ orthogonal to the $e_1$ direction we divide the set of edges normal to $H$ into outgoing edges $(v,v+e_1)$  and incoming edges $(v,v-e_1)$ where $v\in H$.
When restricting an even subgraph $G=(V,E)$ to $\mathbb{T}_n^d\setminus H$, degree considerations allows us to pair the edges normal to $H$ such that each pair of edges is connected in $G\setminus H$. If an incoming edge is connected in $G\setminus H$ to an outgoing edge, this connection is a path in $G$ winding all the way around the torus and hence a cluster of size at least $n$ in $E$\footnote{The astute reader might wonder why we do not simply make this argument on the sphere (corresponding to the wired random-cluster measure). It is exactly this lower bound on the size of topologically non-trivial clusters which fails for the wired measure.}.

If the number of outgoing edges in $E$ is odd, then at least one outgoing edge must be connected to an incoming edge in $\mathbb{T}^d_n\setminus H$. See Figure \ref{fig:wrap around pic}. 
This discussion is completely analogous to the observation for the mirror model on the cylinder made in \cite{kozma2013lower} and motivates the following definition:
\begin{figure}
    \centering
    \includegraphics[scale = 0.6]{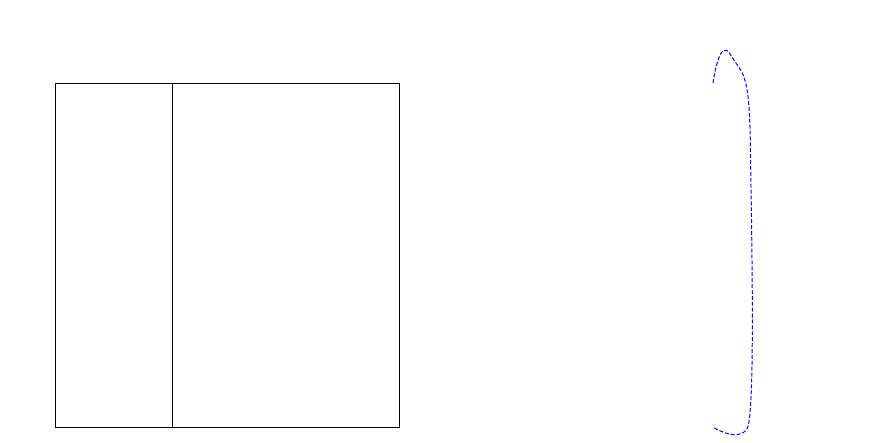}
    \caption{A simple, short, topologically trivial path pictured left and a long path wrapping all the way around the torus on the right. A path wraps all the way around the torus an odd number of times in the $e_1$ direction if and only if it contains an odd number of the outgoing edges from $H$.}
    \label{fig:wrap around pic}
\end{figure}
\begin{definition} \label{def:wrap_arounds}
We say that a loop is \textbf{simple} if it is a path from a vertex to itself such that every other vertex in the path is only visited once. 
 A simple loop $\gamma$ in $\mathbb{T}_n^d$ is a \textbf{wrap-around} if it contains an odd number of outgoing edges of a hyperplane $H$ orthogonal to the $e_1$-direction.

 We say that an even subgraph $G=(V,E)$ of $\mathbb{T}_n^d$  is \textbf{non-trivial} if $E$ contains a wrap-around. Otherwise, say that $G$ is \textbf{trivial}.
\end{definition}
\begin{remark}
Following the above discussion, if $\gamma$ contains an odd number of outgoing edges of the hyperplane $H$, then it must also do so for any hyperplane $H'$ parallel to $H$. As such, the above definition does not depend on the choice of hyperplane. 
Equivalently, $G$ non-trivial if and only if any translation of $G$ is non-trivial.

Of course, a loop might wrap around the torus in several directions, and we could define $j$-wrap-arounds for every cardinal direction in $\mathbb{Z}^d,$ but since we do not use the different choices of direction in our proofs at all, we omit them from the definition.
\end{remark}

To illustrate our use of wrap-arounds, we warm up with the following lemma:
\begin{lemma} \label{juice}
Let $G=(V,E)$ be a trivial even subgraph of $\mathbb{T}_n^d$ and let $\gamma$ be a wrap-around. Then, $G\triangle \gamma:=(V,E\triangle \gamma)$ is non-trivial.
\end{lemma}
\begin{proof}
Fix a hyperplane $H$ and note that, since the number of outgoing edges in $\gamma$ is odd, the parity of outgoing edges from $H$ is different in $G$ than it is in $G\triangle \gamma$. This immediately implies the statement.
\end{proof}
One immediate consequence hereof is the following:
\begin{corollary}
Let $G$ be a fixed, not necessarily even, subgraph of $\mathbb{T}_n^d$ which contains a wrap-around and let $\mathtt{NT}$ denote the event on $\Omega_{\emptyset}(G)$ that a percolation configuration is non-trivial. Then, 
$$
{\UEGop}_{G}[\mathtt{NT}]\geq \frac{1}{2}.
$$
\end{corollary}
\begin{proof}
Fix a wrap-around $\gamma$ in $G$ and let $\eta \sim \UEG_{G}$. By the Haar measure property, we have that $\eta\overset{d}{=}\eta\triangle \gamma$. However, by Lemma \ref{juice}, at least one of $\eta$ and $\eta\triangle\gamma$ is non-trivial. Accordingly, by a union bound,
$$
1=\UEG_G[(\eta\in \mathtt{NT})\cup (\eta\triangle \gamma\in \mathtt{NT})]\leq 2 \UEG_{G}[\mathtt{NT}].
$$
\end{proof}
Thus, if we can exhibit wrap-arounds in the random-cluster model, we get long clusters in the loop $\mathrm{O}$(1) model with positive probability. By using translation invariance of the random-cluster model on the torus, we get a lower bound on the probability of having a long cluster passing through exactly the vertex $0$. This is the main idea of our proof, although we shall be slightly more clever in our application of Lemma \ref{juice} to improve the bound we get.

\begin{definition} \label{def:non_triv}
For an even subgraph $G$ of $\mathbb{T}^d_n,$ we denote by  $\mathcal{C}_{\mathtt{NT}}$ the union of the non-trivial connected components of $G$.
\end{definition}
Note that if $G+v$ is a translate of $G$, then $\mathcal{C}_{\mathtt{NT}}(G+v)=\mathcal{C}_{\mathtt{NT}}(G)+v.$ By translation invariance of the loop O($1$) model on the torus, 
$$
\ell_{x,\mathbb{T}_n^d}[0\in \mathcal{C}_{\mathtt{NT}}]=\ell_{x,\mathbb{T}_n^d}\left[\frac{|\mathcal{C}_{\mathtt{NT}}|}{|\mathbb{T}_n^d|}\right].
$$
Hence, our goal in Section \ref{sec:torus} shall be to lower bound this quantity.

In order to exploit the uniform even subgraph to say something intelligent about the loop $\mathrm{O}$(1) model, we have to build up some technical machinery for the random-cluster model. 
The next couple of sections are dedicated to doing just that. As such, the work herein mostly consists in massaging results from the literature into a form more amenable to our needs. The arguments are slightly technical and a reader who is more eager to get to the proofs of our main theorems might choose to skip it on a first reading. Apart from the proof of \Cref{thm:Loop O(1) mixing}, the results that we shall be needing later on are Lemmata \ref{Mixing}, \ref{local uniqueness} and \ref{Many wrap-arounds}.
\subsection{Local Uniqueness}
The purpose of this section is to extract separating surfaces in the random-cluster model so that one may apply Proposition \ref{prop:separated-graphs}. As a first ingredient, one may note that the construction in \cite[Theorem 1.3]{duminil2020exponential} implies equally well the following lemma: 
\begin{lemma} \label{Mixing}
For any $p>p_c,$ there exists $c>0$ such that for any $n \in \mathbb N$ and any event $A$ depending only on edges of $\Lambda_n$,
$$
\left|\phi^{\xi}_{p,\Lambda_{2n}}[A]-\phi_{p,\mathbb{Z}^d}[A]\right|<\exp(-cn)
$$
for any boundary condition $\xi$. 
\end{lemma}

Our second input comes from Pisztora's construction of a Wulff theory for random-cluster models in dimension at least three \cite{Pis96}. Let $\mathfrak{C}_{n,L,\varepsilon,\theta} \subset \{0,1\}^{E(\Lambda_n)}$ be the event that there exists a cluster $\mathcal{C}_{\max}$ in $\Lambda_n$ such that
\begin{enumerate}
\item[$\bullet$] $\mathcal{C}_{\max}$ is the unique cluster in $\Lambda_n$ touching all faces of $\partial\Lambda_n$.
\item[$\bullet$] $|\mathcal{C}_{\max}|\geq (\theta-\varepsilon) (2n+1)^d$.
\item[$\bullet$] There are at most $\varepsilon n^{d}$ vertices in $\Lambda_n$ that do not lie on $\mathcal{C}_{\max}$ and that lie on clusters larger than $L$.
\end{enumerate}

 One should think of $\mathfrak{C}_{n,L,\varepsilon,\theta}$ as the event that the infinite cluster can be seen at scale $n$: There is one huge cluster, which is a finite-volume analogue of the infinite cluster, and the vast majority of the other clusters are of constant order. One may think of $\theta$ as the density of the infinite cluster, and $L$ as an upper bound on the expected size of a finite cluster in infinite volume (which is of order correlation length to the power of the dimension).

Combining \cite[Theorem 1.2]{Pis96} and \cite[Theorem 2.1]{Bod05} yields the following:
\begin{proposition} \label{Super glue}
For all $p>p_c,$ if $\theta:=\phi_{p,\mathbb{Z}^d}[0\cc \infty],$ for all $0<\varepsilon<\theta/2,$ there exist $L$ and $c>0$ such that
$$
\phi^{\xi}_{p,\Lambda_n}[\mathfrak{C}_{n,L,\varepsilon,\theta}]\geq 1-\exp(-cn^{d-1})
$$
for all $n\in \N$ and every boundary condition $\xi$.
\end{proposition}

A first consequence of this is a result on annular domains, which is essential for showing that the loop $\mathrm{O}(1)$ model is not very sensitive to boundary conditions.

  \begin{figure}
      \centering
      \includegraphics[scale =0.8]{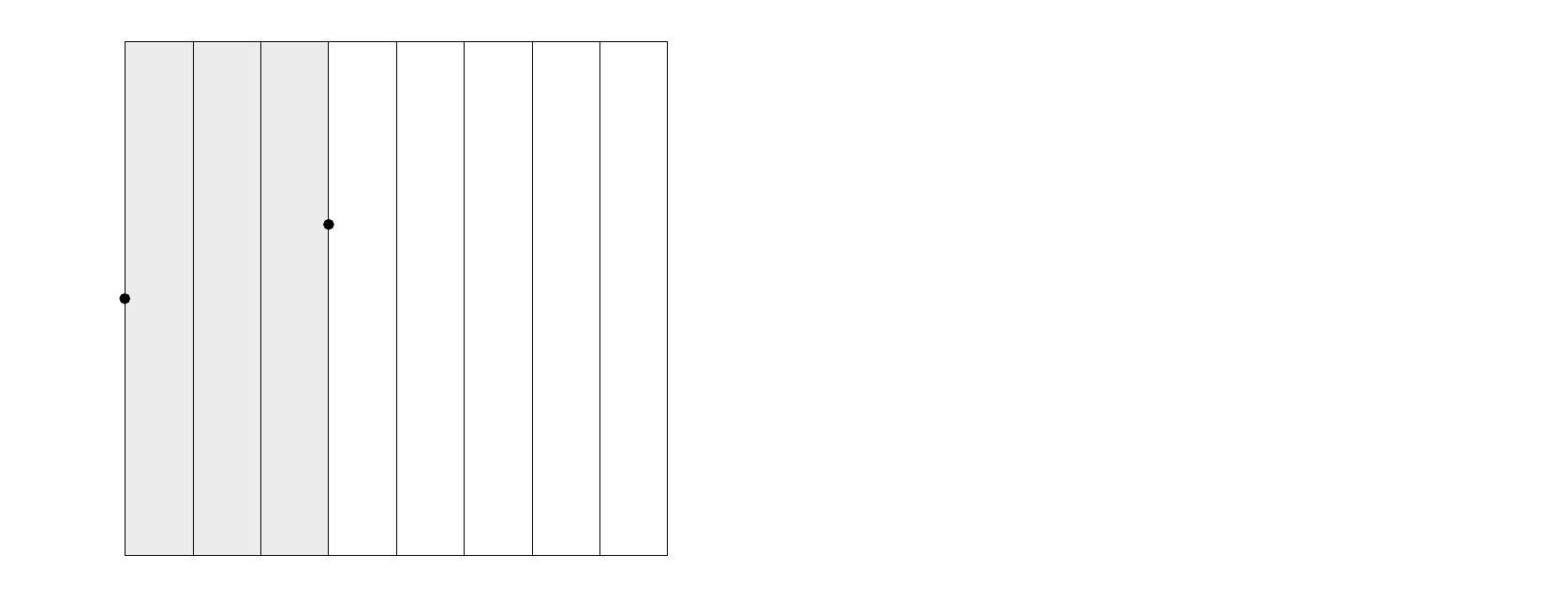}
      \caption{On the left: A sketch of the event $E_v$ from the side. When exploring a cluster traversing the annulus from the inside to a given face, we have a probability to connect to all other faces of a given box every time we enter a new strip (corresponding to a translate of $\tilde{S}_h$). This probability is uniform in the past configuration. On the right: How to apply the same argument in the last direction. Using the fact that a crossing from top to bottom of $B^{v_k}$ must intersect many transversal slabs, it is unlikely that this happens without the cluster of $v_k$ also connecting to the left and right sides of the box.}
      \label{fig:Ex_event}
  \end{figure}

\begin{lemma}[Local Uniqueness] \label{local uniqueness} For $d\geq 3$ and every $p>p_c,$ there exists $c>0$ with the following property:
If $\mathtt{UC}_n$ is the event that $\omega|_{\Lambda_{2n}\setminus \Lambda_n}$ has a unique cluster crossing from $\partial\Lambda_n$ to $\partial \Lambda_{2n}$, then
 $$
\phi^{\xi}_{p,\Lambda_{2n}}[\mathtt{UC}_n] > 1-\exp(-cn)
$$
for all $n\in \N$ and every boundary condition $\xi$.
\end{lemma}
\begin{remark}
Note that the above statement is also true for $d=2$ by sharpness and duality arguments (see \Cref{sec:planar}) and for $d=1$ since here, there is no $p>p_c$. However, the proof we give below very specifically uses the fact that $d\geq 3$.
\end{remark}
\begin{remark} In the following, we are juggling several constants. Our convention here and throughout will be to remark upon the changing of the value of a constant $c$ by denoting the new one $c'$ the first time it appears. Afterwards, to prevent notational bloat, we shall revert to simply writing $c$.
\end{remark}
\begin{proof}
The strategy for the proof comes in two steps: First, we show that with probability exponentially close to $1$, there is a unique large cluster in $\Lambda_{2n}\setminus\Lambda_n$ of large volume and then we show that, again with high probability, any crossing must be part of this one big cluster.

  For any finite set $\mathcal{B}$ of translates $B_j$ of $\Lambda_{n/2},$ observe the auxiliary graph $G_\mathcal{B}$ with vertices $j\in \{1,...,|\mathcal{B}|\}$ and an edge $(j,l)$ if $B_j\cap B_l$ contains a translate of $\Lambda_{n/4}$. For each edge $(j,l)$, let $B_{j,l}$ denote a choice of such a translate. For our purposes, $\mathcal{B}$ will be the cover of $\Lambda_{2n}\setminus \Lambda_n$ consisting of all translates of $\Lambda_{n/2}$ contained in $\Lambda_{2n}\setminus \Lambda_{n-1}.$ Since the centre of every such translate lies on $\partial \Lambda_{3n/2},$ we get that
  \begin{align} \label{eq:B_bounded} 
  |\mathcal{B}|\leq Cn^{d-1}.
  \end{align}

For each $j$, let $\mathfrak{C}^j$ denote the event that the corresponding translate of the event $\mathfrak{C}_{n/2,L,\theta/4^d, \theta}$ from Proposition \ref{Super glue} occurs in $B_j$. Similarly, for an edge $(j,l)$ of $G_\mathcal{B},$ let $\mathfrak{C}^{j,l}$ denote the event that the corresponding translate of $\mathfrak{C}_{n/4,L,\theta/4^d,\theta}$ occurs in 
$B_{j,l}$.

By a union bound, \Cref{Super glue} and \eqref{eq:B_bounded}, we see that
$$
\phi^{\xi}_{p,\Lambda_{2n}}\left[\left(\cap_{j}\mathfrak{C}^j \right)\cap \left(\cap_{(j,l)\in E(G_\mathcal{B})} \mathfrak{C}^{j,l}\right)\right]\geq 1-C n^{d-1}\exp(-c(n/2)^{d-1})-C^2n^{2(d-1)}\exp(-c (n/4)^{d-1}),
$$
which is at least $1-\exp(-c'n^{d-1})$ for an adjusted value $c'$ for $n$ sufficiently large. Possibly adjusting the value of $c'$ further, we get such a bound for all $n$.

For future reference, we shall abbreviate $\mathfrak{C}^{\max}=\left(\cap_{j}\mathfrak{C}^j \right)\cap \left(\cap_{(j,k)\in E(G_B)} \mathfrak{C}^{j,k}\right)$. In conclusion,
\begin{align}\label{eq:Cmax} 
\phi^{\xi}_{p,\Lambda_{2n}}[\mathfrak{C}^{\max}] \geq 1- \exp(-cn). 
\end{align}

Now, for $(j,k)\in E(G_\mathcal{B}),$ on the event $\mathfrak{C}^j\cap \mathfrak{C}^k,$ there is a unique large cluster $\mathcal{C}^j$ contained in $B_j$ of size at least $\theta/4^d (n/2)^d$ and likewise for $k$. However, on $\mathfrak{C}^{j,k},$ there is a cluster $\mathcal{C}^{j,k}$ in $B^{j,k}$ of size $\frac{4^d-1}{4^d} \theta\left(\frac{n}{4}\right)^d>\theta/4^d \left(\frac{n}{2}\right)^d$. Accordingly, $\mathcal{C}^{j,k}\subset \mathcal{C}^j\cap \mathcal{C}^k$. Since $G_\mathcal{B}$ is connected, we get that, on $\mathfrak{C}^{\max},$ all the $\mathcal{C}^j$ are part of one big cluster in $\omega|_{\Lambda_{2n}\setminus \Lambda_n}$.

Now, for the second part of the argument, let $v\in \partial \Lambda_n$ be given and let $\mathcal{C}_v$ denote the cluster of $v$ in $\omega|_{\Lambda_{2n}\setminus \Lambda_n}$.
Let $\mathcal{E}_v$ denote the event that $v$ is connected inside of $\Lambda_{2n}\setminus \Lambda_n$ to $\partial \Lambda_{2n}$ and for any $j,$ let $\mathfrak{A}^v_{j}$ be the event that $\mathcal{C}_v\cap B_j$ contains a cluster which touches every face of $\partial B_j$. We denote by $\mathfrak{A}^v$ the event that $\mathfrak{A}^v_{j}$ occurs for some $j$.
 We wish to show that if $\mathcal{C}_v$ crosses $\Lambda_{2n}\setminus \Lambda_n,$ then with high probability, it must touch all faces of some $B_j$. That is,
\begin{align} \label{eq:x_is_touching} 
\phi_{p,\Lambda_{2n}}^{\xi}[\cup_{v\in \partial \Lambda_n} \mathcal{E}_v\setminus \mathfrak{A}^v]
\leq e^{-cn}. 
\end{align}
Let us first see how \eqref{eq:x_is_touching} finishes the proof. 
Notice that
$$
\mathfrak{C}^{\max}\setminus \left(\cup_{v\in \partial \Lambda_n} \mathcal{E}_v\setminus \mathfrak{A}^v\right) \subset \mathtt{UC}_n, 
$$
 since on the former event, any cluster which crosses $\Lambda_{2n}\setminus \Lambda_n$ touches all faces of $B_j$ for some $j,$ any such cluster must be the same as $\mathcal{C}^j,$ and all the $\mathcal{C}^j$ are part of the same cluster. The lemma then follows by combining \eqref{eq:x_is_touching} and \eqref{eq:Cmax}.

Thus, let us establish  \eqref{eq:x_is_touching}. 
 First, by \cite[Theorem 2.1]{Bod05}, there exists a constant $h$ such that $\phi^0_{p,S_h}$ percolates\footnote{Due to the finite size of the graph in all but two directions, the boundary conditions under which we take the infinite volume limit on $S_h$ matter, unlike in the case of $\mathbb{Z}^d$.}, where $S_h$ denotes the slab $\{0,...,h\}^{d-2}\times \mathbb{Z}^2$. In particular, if $\tilde{S_h}:=\{0,...,h\}\times \mathbb{Z}^{d-1},$ we get that $\phi^0_{p,\tilde{S}_h}$ also percolates. For the rest of the proof, we shall assume, without loss of generality, that $n>h$. Combining \cite[Lemma 3.3]{Pis96} with the FKG inequality (\Cref{FKG}), there exists a $c>0,$ depending only on $p$ such that, for any hyper-rectangle $B$ in $\mathbb{Z}^d$ and any vertex $w$, if 
 $\mathfrak{A}^h_B(w)$ denotes the event that the cluster of $w$ in $\tilde{S}_h+w$ touches all faces of $B\cap (\tilde{S}_h+w)$, we have  $\phi^0_{p,\tilde{S}_h+w}[\mathfrak{A}^h_B(w)]\geq c$.

 Now, suppose $v$ is connected to $\{\langle w,e_1\rangle =2n\}$, and call this event $\mathcal{E}_v^{\uparrow}$. On $\mathcal{E}_v^{\uparrow}$, we must have that $\mathcal{C}_v$ crosses $\frac{n}{h}$ disjoint translates of $\tilde{S}_h$. To use this, we explore the cluster of $v$ from $\Lambda_k$ to $\{\langle w,e_1\rangle =2n\}$ one translate of $\tilde{S}_h$ at a time and denote by $v_k$ the first vertex of $\mathcal{C}_v$ we encounter in $v+kh e_1+ \tilde S_h$. See Figure \ref{fig:Ex_event}.

Let $B^{v_k}$ denote some $B_j$ such that $v_k,v_k+h e_1\in B^{v_k}$. Similarly, let $\Past(v_k)$ denote the state of all previously discovered edges (open or closed) and let $E(\Past(v_k))$ denote the set of discovered edges. By the Domain Markov Property (\Cref{DMP}), we have that
$$
\phi_{p,\mathbb{Z}^d}[\cdot|\;\Past(v_k)]=\phi_{p,\mathbb{Z}^d\setminus E(\Past(v_k))}^{\xi(\Past(v_k))}[\cdot],
$$
where $\xi(\Past(v_k))$ are the boundary conditions which are wired on the component of $v$ in $\Past(v_k)$ and free otherwise. Since for any probability measure $\nu$ and events $U,V$ with $\nu(U\cap V)>0$, we have $\nu[\cdot |\; U\cap V]=\nu_V[\cdot |\; U],$ where $\nu_V[\cdot]=\nu[\cdot|\; V]$, we conclude that
$$
\phi_{p,\mathbb{Z}^d}[\mathfrak{A}^h_{B^{v_k}}(v_k)|\Past(v_k),\mathcal{E}^{\uparrow}_v]=\phi_{p,\mathbb{Z}^d\setminus E(\Past(v_k))}^{\xi(\Past(v_k))}[\mathfrak{A}^h_{B^{v_k}}(v_k)|\tilde{\mathcal{E}}^{\uparrow}_v],
$$
where $\omega\in \tilde{\mathcal{E}}^{\uparrow}_v$ if $\omega\cup \Past(v_k)\in \mathcal{E}^{\uparrow}_v$ (and $\Past(v_k)$ is identified with its open edges). Since $\tilde{\mathcal{E}}^{\uparrow}_v$ is increasing, we can apply the FKG inequality to get that
$$
\phi_{p,\mathbb{Z}^d\setminus E(\Past(v_k))}^{\xi(\Past(v_k))}[\mathfrak{A}^h_{B^{v_k}}(v_k)|\; \tilde{\mathcal{E}}^{\uparrow}_v]\geq \phi_{p,\mathbb{Z}^d\setminus E(\Past(v_k))}^{\xi(\Past(v_k))}[\mathfrak{A}^h_{B^{v_k}}(v_k)].
$$
Applying the comparison between boundary conditions (cf. \Cref{thm:comparison theorem}$i$), we see that
$$
 \phi_{p,\mathbb{Z}^d\setminus E(\Past(v_k))}^{\xi(\Past(v_k))}[\mathfrak{A}^h_{B^{v_k}}(v_k)]\geq \phi^0_{p,\tilde{S}_h+v_k}[\mathfrak{A}^h_{B^{v_k}}(v_k)]\geq c
$$
Iterating on the above, we see that, conditionally on $\mathcal{E}_v^{\uparrow}$, $\mathfrak{A}^h_{B^{v_k}}(v_k)$ occurs for some $k$ with 
probability at least $1-(1-c)^{n/h} > 1- e^{-c'n}$ for some appropriate choice of $c'$. On the event $\mathfrak{A}^h_{B^{v_k}}(v_k),$ we have that the cluster of $v_k$ in $\mathcal{C}_v\cap B^{v_k}$ touches all faces of $B^{v_k}$ except possibly those orthogonal to $e_1$. However, we may apply a similar exploration argument to get that
$$
\phi_{p,\mathbb{Z}^d}[\mathfrak{A}^{v_k}|\;\mathfrak{A}^h_{B^{v_k}}(v_k)]>1-e^{-cn},
$$
see \Cref{fig:Ex_event}. Note that if $\kappa$ denotes the first $k$ such that $\mathfrak{A}^{h}_{B^{v_k}}(v_k)$ occurs, then
$$
\mathcal{E}_v^{\uparrow}\cap(\kappa<\infty)\cap (\mathfrak{A}^{v_{\kappa}})\subset \mathfrak{A}^v
$$
and therefore, on $\mathcal{E}^{\uparrow}_v\setminus \mathfrak{A}^v,$ either $\kappa=\infty$ or there is some $k$ and a $w$ on the boundary of $v+khe_1+S_h$ such that $\mathfrak{A}^h_{B^w}(w)\setminus\mathfrak{A}^w$ occurs. Therefore, a union bound shows that 
$$
\phi_{p,\mathbb{Z}^d}[\mathcal{E}_v^{\uparrow}\setminus \mathfrak{A}^v]\leq \phi_{p,\mathbb{Z}^d}[\kappa=\infty|\; \mathcal{E}_v^{\uparrow}]+\sum_{k=1}^{n/h}\sum_{w\in \partial (v+khe_1+S_h)}\phi_{p,\mathbb{Z}^d}[\mathfrak{A}^h_{B^{w}}(w)\setminus \mathfrak{A}^w]\leq (1+Cn^d)e^{-cn}\leq e^{-c'n}.
$$
The argument in the case where $v$ is connected to another face is similar. Thus, summing over the $2d$ faces of $\Lambda_{2n}$, we get that 
$$\phi_{p,\mathbb{Z}^d}[\mathcal{E}_v\setminus \mathfrak{A}^v]\leq 2de^{-c'n}\leq e^{c''n}.$$
By yet another union bound, 
$$
\phi_{p,\mathbb{Z}^d}[\cup_{v\in \partial \Lambda_n} \mathcal{E}_v\setminus \mathfrak{A}^v]\leq C n^{d-1}e^{-c''n}\leq e^{-c'''n}.
$$
  By \Cref{Mixing}, we get \eqref{eq:x_is_touching}, which is what we wanted. 
\end{proof}

\subsection{Insensitivity to boundary conditions and mixing of the loop $\mathrm{O}$(1) model}

\Cref{local uniqueness} is the random-cluster analogue of the connectivity property that we used for the uniform even graph of $\mathbb{Z}^d$ in Corollary \ref{Surprising fact}. 
Combined with Lemma \ref{Mixing}, we get that the loop $\mathrm{O}$(1) model is insensitive to boundary conditions:

\mixing*
\begin{proof} 
    Note that the result for $x=1$ follows from the stronger result in \Cref{UEG independence}. The event $\mathtt{UC}_n$ from \Cref{local uniqueness} on percolation configurations $\omega\subset \Z^d$ or $\omega\subset G$ is equal to the event that there is a cluster in $\omega|_{\Lambda_{2n}\setminus \Lambda_n}$ containing a separating surface between $\omega\vert_{\Lambda_n}$ and $\omega\vert_{\Lambda_{2n}^c}$.
    Then, \Cref{prop:separated-graphs} applies, so that whenever $\omega\in \mathtt{UC}_n$, we have $\UEG_{\omega}[A]=\UEG_{\omega\vert_{\Lambda_{2n}}}[A]$.
    It follows that $\id_{\mathtt{UC}_n}(\omega)\UEG_{\omega}[A]$ is a random variable which is measurable with respect to the state of $\omega$ on edges in $\Lambda_{2n}$. Furthermore, it is positive and bounded from above by $1$. 
    Denoting by $d_{\TV}$ the total variation distance between probability measures, we can conclude:
    \begin{align*}
        |\ell^\xi_{x,G}[A]-\ell_{x,\mathbb{Z}^d}[A]|
        & = \big|\phi_{x,G}^\xi[\UEG_\omega[A]]-\phi_{x,\mathbb{Z}^d}[\UEG_\omega[A]]\big| \\
        & \le \big|\phi_{x,G}^\xi[\id_{\mathtt{UC}_n}(\omega)\UEG_\omega[A]]-\phi_{x,\mathbb{Z}^d}[\id_{\mathtt{UC}_n}(\omega)\UEG_\omega[A]]\big| +2-\phi_{x,G}^\xi[\mathtt{UC}_n]-\phi_{x,\mathbb{Z}^d}[\mathtt{UC}_n] \\
        & \le d_{\TV}(\phi^{\xi}_{G}|_{\Lambda_{2n}} ,\phi_{\mathbb{Z}^d}|_{\Lambda_{2n}})+\exp(-cn) \\
        & \leq \exp(-c'n)+\exp(-cn)\\
        &\leq \exp(-c''n),
    \end{align*}
    where in the second inequality, we used Lemma \ref{local uniqueness} and in the third, we used Lemma \ref{Mixing}. 
\end{proof} 

By a similar argument, one may show an actual mixing result on $\mathbb{Z}^d$:
\begin{theorem} \label{theorem:the_real_mixing}
For $d\geq 2$ and $x\in (x_c,1]$, there exists $c>0$ such that for any $n \in \mathbb N$ and any events $A$ and $B$ such that $A$ depends only on the edges in some box $v_A+\Lambda_n$ and $B$ depends only on the edges in some box $v_B+\Lambda_n$, for two vertices $v_A$ and $v_B$ such that $|v_A-v_B|\geq 6n,$ then 
$$
|\ell_{x,\mathbb{Z}^d}[A\cap B]-\ell_{x,\mathbb{Z}^d}[A]\ell_{x,\mathbb{Z}^d}[B]|<\exp(-cn).
$$
\end{theorem}
\begin{proof}
Again, the result for $x=1$ follows from the stronger result in \Cref{UEG independence}. Denote by $\mathtt{UC}_n^A$ the event that $\omega|_{v_A+\Lambda_{2n}\setminus \Lambda_n}$ has a unique cluster crossing from inner to outer radius and define $\mathtt{UC}_n^B$ similarly. By \Cref{local uniqueness}, we have that
\begin{equation} \label{No 1}
\ell_{x,\mathbb{Z}^d}[A\cap B]=\phi_{x,\mathbb{Z}^d}[\UEG_{\omega}[A\cap B]]=\phi_{x,\mathbb{Z}^d}[\id_{\mathtt{UC}_n^{A}}(\omega)\id_{\mathtt{UC}_n^{B}}(\omega)\UEG_{\omega}[A\cap B]]+O(\exp(-cn)).
\end{equation}

Denoting by $\omega_A=\omega|_{v_A+\Lambda_{2n}}$ and $\omega_B=\omega|_{v_B+\Lambda_{2n}},$ we can apply \Cref{UEG independence} to get that
\begin{align}\label{No 2}
\id_{\mathtt{UC}_n^{A}}(\omega)\id_{\mathtt{UC}_n^{B}}(\omega)\UEG_{\omega}[A\cap B]&=\id_{\mathtt{UC}_n^{A}}(\omega)\id_{\mathtt{UC}_n^{B}}(\omega)\UEG_{\omega}[A]\UEG_{\omega}[B]\nonumber \\
&=\id_{\mathtt{UC}_n^{A}}(\omega_A)\id_{\mathtt{UC}_n^{B}}(\omega_B)\UEG_{\omega}[A]\UEG_{\omega}[B]. 
\end{align}

As before, $\id_{\mathtt{UC}_n^{A}}(\omega_A)\UEG_{\omega}[A]$ is a measurable function of $\omega_A$, positive, and bounded below $1$, so we shall once again attempt to bound a total variation distance. Let $(\tilde{\omega}_A,\tilde{\omega}_B)$ a coupling of two percolation configurations with the same marginals as $(\omega_A,\omega_B)$ such that $\tilde{\omega}_A$ and $\tilde{\omega}_B$ are independent. Then, by \cite[Corollary 1.4]{duminil2020exponential},
we have that there exists a coupling $P$ of the two such that $\tilde{\omega}_B=\omega_B$ almost surely and 
$$
P[\omega_A\neq \tilde{\omega}_A|\; \omega_B]<\exp(-cn).
$$
Equivalently,
$$
d_{\mathrm{TV}}((\omega_A,\omega_B),(\tilde{\omega}_A,\tilde{\omega}_B))<\exp(-cn),
$$
and accordingly,
\begin{align}\label{No 3}
&\phi_{x,\mathbb{Z}^d}[\id_{\mathtt{UC}_n^{A}}(\omega_A)\id_{\mathtt{UC}_n^{B}}(\omega_B)\UEG_{\omega}[A]\UEG_{\omega}[B]]\nonumber \\
=&\phi_{x,\mathbb{Z}^d}[\id_{\mathtt{UC}_n^{A}}(\omega_A)\UEG_{\omega}[A]]\times\phi_{x,\mathbb{Z}^d}[\id_{\mathtt{UC}_n^{B}}(\omega_B)\UEG_{\omega}[B]]+O(\exp(-cn)). 
\end{align}

Now, to finish, we note that
\begin{align}\label{No 4}
\ell_{x,\mathbb{Z}^d}[A]\ell_{x,\mathbb{Z}^d}[B]&=\phi_{x,\mathbb{Z}^d}[\UEG_{\omega}[A]]\times\phi_{x,\mathbb{Z}^d}[\UEG_{\omega}[B]]\nonumber\\
&=\phi_{x,\mathbb{Z}^d}[\id_{\mathtt{UC}^A_n}(\omega_A)\UEG_{\omega}[A]]\times \phi_{x,\mathbb{Z}^d}[\id_{\mathtt{UC}^B_n}(\omega_B)\UEG_{\omega}[B]]+O(\exp(-cn)).
\end{align}
Combining \eqref{No 1}, \eqref{No 2}, \eqref{No 3} and \eqref{No 4} yields the desired.
\end{proof}

\subsection{Existence of many wrap-arounds}

 \Cref{thm:Loop O(1) mixing} enables us to apply our observations from Section \ref{Road map} by first arguing directly on the torus and then saying that the model on $\mathbb{Z}^d$ does not look too different. 
 First, we prepare for proving the existence of sufficiently many wrap-arounds.

\begin{lemma} \label{Exclusion tolerance}
There exists a continuous function $f:(0,1)^2\to [0,1)$ with the following property:

Let $G=(V,E)$ be a finite graph and $\xi$ a boundary condition.
For $p_1<p_2,$ $f(p_1,p_2)>0$ and there exists an increasing coupling $P$ between $\omega_1\sim \phi^{\xi}_{p_1,G}$ and $\omega_2\sim \phi^{\xi}_{p_2,G}$ such that for any random finite set of edges $F\subset \omega_2$ 
measurable with respect to $\omega_2$, we have 
$$
P[\omega_1(e)=0\; \forall e\in F(\omega_2)|\; \omega_2]\geq P[f(p_1,p_2)^{|F(\omega_2)|}\; |\; \omega_2].
$$
\end{lemma}
\begin{proof}
Pick an ordering $(e_j)_{1\leq j\leq |E|}$ of the edges and let $U_j$ be an i.i.d. family of uniforms on $[0,1]$. For $i\in \{1,2\},$ define the target 
$\mathfrak{t}_{e_j,i}:\{0,1\}^{\{e_1,...,e_{j-1}\}}\to (0,1)$ as the conditional probability under $\phi_{p_i}$ that the edge $e_j$ is open given the state of the previous edges, i.e.
$$
\mathfrak{t}_{e_j,i}(\mathfrak{w})=\phi^{\xi}_{p_i,G}[\omega_{e_j}\;|\; \omega_{e_l}=\mathfrak{w}_{e_l} \; \forall l\leq j-1 ].
$$
Then, recursively setting
$$
\omega_{i}(e_j)=\id_{U_j\leq \mathfrak{t}_{e_j,i}(\omega_i|_{\{e_1,..,e_{j-1}\}})}
$$
yields an increasing coupling between the two random graphs $\omega_{i}\sim \phi_{p_i,G}^{\xi}$.

Let us first remark that if we can prove that
\begin{equation} \label{target gap}
\mathfrak{t}_{e_j,2}(\mathfrak{w})-\mathfrak{t}_{e_j,1}(\mathfrak{w}')\geq \min\left\{p_2-p_1,\frac{p_2}{2-p_2}-\frac{p_1}{2-p_1}\right\}
\end{equation}
deterministically for any $\mathfrak{w} \succeq \mathfrak{w}'$, then we are done. 

 To see this, note that, conditional on $F$, the event that $\omega_1(e)=0$ for every $e\in F$ 
 is the event that $ \mathfrak{t}_{e,1}< U_e$ conditional on $U_e\leq \mathfrak{t}_{e,2}$ for all $e\in F$.
 Since $\mathfrak{t}_{e,2}\leq p_2$ (cf.  \Cref{thm:comparison theorem}$iii$)),  we can set 
 $$
 f(p_1,p_2):=\frac{1}{p_2}\min\left\{p_2-p_1,\frac{p_2}{2-p_2}-\frac{p_1}{2-p_1}\right\},
 $$
 and the rest is merely computation.
 
 Accordingly, let us establish (\ref{target gap}). By monotonicity in boundary conditions (\Cref{thm:comparison theorem} $i)$), $\mathfrak{t}$ is increasing in $\mathfrak{w}.$ Therefore, we can assume without loss of generality that $\mathfrak{w}=\mathfrak{w}'$.
 Now, if $A_e$ denotes the event that the end-points of $e$ are connected in $\omega\setminus \{e\}$, we have for $i \in \{1,2\}$ that

 $$
 \mathfrak{t}_{e_j,i}(\mathfrak{w})=p_i \phi^{\xi}_{p_i,G}[A_{e_j}\;| \; \omega_{e_l}=\mathfrak{w}_{e_l}\;\forall l\leq j-1]+\frac{p_i}{2-p_i}\left(1-\phi^{\xi}_{p_i,G}[A_{e_j}\;| \; \omega_{e_l}=\mathfrak{w}_{e_l}\;\forall  l\leq j-1]\right),
 $$
 whence, since $p_2>p_1$ and $A_{e}$ is increasing,
 \begin{align*}
 \mathfrak{t}_{e_j,2}(\mathfrak{w})-\mathfrak{t}_{e_j,1}(\mathfrak{w}) &=\left(p_2-p_1\right)\phi^{\xi}_{p_1,G}[A_{e_j}\;| \; \omega_{e_l}=\mathfrak{w}_{e_l}\;\forall l\leq j-1]\\
 &+\left(\frac{p_2}{2-p_2}-\frac{p_1}{2-p_1}\right)\left(1-\phi^{\xi}_{p_2,G}[A_{e_j}\;| \; \omega_{e_l}=\mathfrak{w}_{e_l}\;\forall l\leq j-1]\right)\\
 &+\left(p_2-\frac{p_1}{2-p_1}\right)\left(\phi^{\xi}_{p_2,G}[A_{e_j}\;| \; \omega_{e_l}=\mathfrak{w}_{e_l}\;\forall l\leq j-1]-\phi^{\xi}_{p_1,G}[A_{e_j}\;| \; \omega_{e_l}=\mathfrak{w}_{e_l}\;\forall l\leq j-1]\right)\\
 &\geq \min\left\{p_2-p_1,\frac{p_2}{2-p_2}-\frac{p_1}{2-p_1}\right\}.
 \end{align*}
\end{proof}

This enables us to reprove the classical result  \cite[Lemma 4.2]{AllTheAuthors} in Bernoulli percolation for the random-cluster model, allowing for the control of the number of crossings in a box. 

For an event $A$ and $r>0$, let $I_r(A)$ denote the set of $\omega$ such that the Hamming distance from $\omega$ to $\Omega\setminus A$ is at least $r$, i.e. changing the state of any $r-1$ edges of $\omega$ cannot produce a configuration outside of $A$. For $A$ the event that $\Lambda_n$ contains a crossing between two opposite faces, we remark that $I_r(A)$ is exactly the event that $\omega$ contains at least $r$ edge-disjoint crossings.

\begin{lemma} \label{Robustening}
There exists a continuous function $f:(0,1)^2\to (0,1)$ such that for any increasing event $A$, finite graph $G$, boundary condition $\xi,$ $r\in \mathbb{N}$ and $p_1<p_2,$ then
$$
f(p_1,p_2)^{r}\left(1-\phi^{\xi}_{p_2,G}[I_r(A)]\right)\leq 1-\phi^{\xi}_{p_1,G}[A].
$$
\end{lemma}
\begin{proof}
Let $\omega_2\sim \phi_{p_2,G}^{\xi}$ and note that, on the event $\omega_2\not\in I_r(A),$  there exists a (possibly empty) set of edges $F$ such that $|F|\leq r,$ every edge in $F$ is open and $\omega_2\setminus F \not\in A$.
This set is not necessarily unique, but we may simply posit some rule for resolving ambiguities. Under any such choice, we see that $F$ becomes measurable with respect to $\omega_2$.

Hence, letting $P$ denote the coupling from Lemma \ref{Exclusion tolerance}, we see that  
\begin{align*}
1-\phi^{\xi}_{p_1,G}[A] &\geq P[\omega_2\not \in I_r(A), \omega_1\not \in A] \\
&\geq P[\omega_1(e)=0\;\forall e\in F|\; \omega_2\not\in I_r(A)]\left(1-\phi^{\xi}_{p_2,G}[I_r(A)]\right)\\
&\geq f(p_1,p_2)^r\left(1-\phi^{\xi}_{p_2,G}[I_r(A)]\right),
\end{align*}
which is what we wanted.
\end{proof}

\begin{lemma} \label{Many wrap-arounds} For a percolation configuration $\omega$ on $\mathbb{T}_n^d,$
let $\boldsymbol{N}$ denote the maximal number of edge disjoint wrap-arounds in $\omega$. Then, for every $p>p_c,$ there exist $\alpha,c>0$ such that it holds that
$$
\phi_{p,\mathbb{T}_n^d}[\boldsymbol{N}\leq \alpha n^{d-1}]\leq \exp(-c n^{d-1})
$$
for every $n$.
\end{lemma}
\begin{proof}
The proof essentially follows in two steps, which we outline heuristically here: Pick $\delta>0$ small enough and show that under $\phi_{p-\delta, \mathbb{T}_n^d},$ there is a crossing winding around the torus once with high probability. Then, we will use the previous lemma to argue that, under $\phi_{p, \mathbb{T}_n^d},$ there is a large number of such crossings with high probability.

Let $\delta<p-p_c$ and, for $j\in\{0,1\}$, observe the two boxes
$B_j=\left(\Lambda_{n-1}+j ne_1\right)/2n\mathbb{Z}^d$ and let $\mathfrak{C}^j$ denote the event that the relevant translate of the event $\mathfrak{C}_{n-1,L,\theta,\theta/4^d}$ from Proposition \ref{Super glue} occurs in $B_j$. Furthermore, let $\tilde{\mathfrak{C}}^j$ denote the event that the relevant translate of $\mathfrak{C}_{(n-1)/2,L,\theta,\theta/4^d}$ occurs in $\tilde{B}_j:=\Lambda_{(n-1)/2}+(-1)^j\frac{n}{2}e_1$. Note that $\tilde{B}_j\subset B_0\cap B_1$.

By a union bound and Proposition \ref{Super glue}, we have
\begin{align} \label{eq:four_big_clust}
\phi_{p-\delta,\mathbb{T}_n^d}\left[\mathfrak{C}^0\cap \mathfrak{C}^1\cap \tilde{\mathfrak{C}}^0\cap \tilde{\mathfrak{C}}^1\right]\geq 1-2\exp(-cn^{d-1})-2\exp(-c(n/2)^{d-1})\geq 1-\exp(-c'n^{d-1})
\end{align}
for a slightly smaller constant $c' > 0$. Let $\mathtt{SL}$ be the event that there is a simple loop wrapping around the torus once in the sense of \Cref{def:wrap_arounds}. Let us show that
\begin{align} \label{eq:bigclust_in_simple_loop}
\mathfrak{C}^0\cap \mathfrak{C}^1\cap \tilde{\mathfrak{C}}^0\cap \tilde{\mathfrak{C}}^1 \subset \mathtt{SL}. 
\end{align} 

On the event $\mathfrak{C}^0\cap \mathfrak{C}^1,$ for each $j \in \{0,1\},$ there is a unique large cluster $\mathcal{C}^j$ contained in each $B_j$, and it contains a path from the left side of $B_j$ to its right side. Furthermore, they are the unique clusters of size at least $\frac{\theta}{4^d} (n-1)^d$ in $B_0$ and $B_1$ respectively. Finally, on $\tilde{\mathfrak{C}}^j,$ there is a cluster $\tilde{\mathcal{C}}^j$ in $\tilde{B}_j$ of size $\frac{4^d-1}{4^d}\theta \left(\frac{n-1}{2}\right)^d>\frac{\theta}{4^d} (n-1)^d$. All in all, on $\mathfrak{C}^0\cap \mathfrak{C}^1\cap \tilde{\mathfrak{C}}^0\cap \tilde{\mathfrak{C}}^1$, then $\tilde{\mathcal{C}}^j\subset \mathcal{C}^0\cap \mathcal{C}^1$. That is, all the large clusters intersect. 

Let us now construct a simple loop of open edges wrapping around the torus once. In $\mathcal{C}^j,$ there is a path $\gamma^j$ from the left to the right face of $B^j$, which is connected by paths to both $\tilde{\mathcal{C}}^0$ and $\tilde{\mathcal{C}}^1$. By using the latter paths, we may glue $\gamma^0$ and $\gamma^1$ together, which yields a cluster containing a simple loop wrapping around the torus once and \eqref{eq:bigclust_in_simple_loop} follows. See \Cref{fig:Torus_glue}.

 \begin{figure}
     \centering
     \includegraphics[scale=0.5]{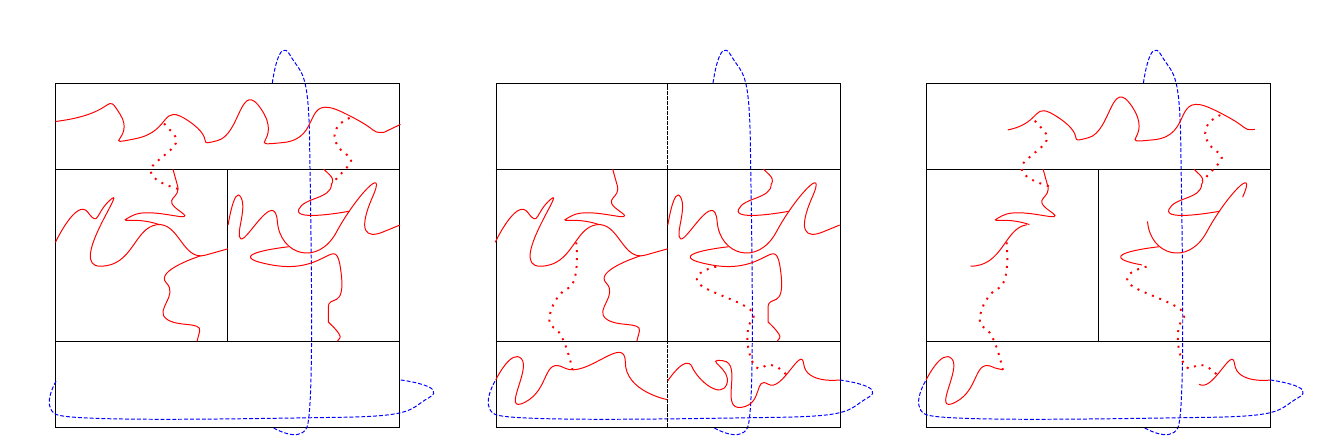}
     \caption{The construction of a simple loop wrapping around the torus once. On the left, on $\mathfrak{C}^0\cap \tilde{\mathfrak{C}}^0\cap \tilde{\mathfrak{C}}^1,$ the path $\gamma^0$ is connected to both of the clusters $\tilde{\mathcal{C}}^0$ and $\tilde{\mathcal{C}}^1$ via the dotted paths. In the middle panel, on $\mathfrak{C}^1\cap \tilde{\mathfrak{C}}^0\cap \tilde{\mathfrak{C}}^1,$ the same is true of the path $\gamma^1$. From this, a wrap-around may be extracted as displayed in the last picture.}
     \label{fig:Torus_glue}
 \end{figure}

 Since a loop cannot vanish by adding additional edges, $\mathtt{SL}$ is increasing.
Thus, by Lemma \ref{Robustening}, we have a continuous function $f: (0,1)^2 \to (0,1)$ such that, for every $\alpha>0,$
$$
\exp(-c'n^{d-1})\geq 1- \phi_{p-\delta,\mathbb{T}_n^d}[\mathtt{SL}]\geq f(p-\delta,p)^{\alpha n^{d-1}}\left(1-\phi_{p,\mathbb{T}_n^d}[I_{\alpha n^{d-1}}(\mathtt{SL})]\right),
$$
where we also used \eqref{eq:four_big_clust} and \eqref{eq:bigclust_in_simple_loop}. 
 Accordingly,
 $$
 \phi_{p,\mathbb{T}_n^d}[I_{\alpha n^{d-1}}(\mathtt{SL})]\geq 1- \exp\left(\left(\alpha \log(1/f(p-\delta,p))-c'\right) n^{d-1}\right),
 $$
 yielding that
 $$
 \phi_{p,\mathbb{T}_n^d}[I_{\alpha n^{d-1}}(\mathtt{SL})]\geq 1-\exp\left(c''n^{d-1}\right)
 $$
 for an adjusted value $c''>0$ and suitably small $\alpha$. 
 
 All that is left to do is to notice that $I_{\alpha n^{d-1}}(\mathtt{SL})$ is the event that there exist at least $\alpha n^{d-1}$ edge-disjoint simple loops wrapping around the torus exactly once, which is a sub-event of $(\boldsymbol{N}\geq \alpha n^{d-1}).$
\end{proof}

\subsection{Infinite expected cluster sizes from torus wrap-arounds} \label{sec:torus}

Finally, we are in position to prove our main result following the outline given in  \Cref{Road map}: We exploit Lemma \ref{juice} and Lemma \ref{Many wrap-arounds} to lower bound the number of large clusters in the loop $\mathrm{O}$(1) model on the torus. Then, we use Lemma \ref{Mixing} and Lemma \ref{local uniqueness} to compare the loop $\mathrm{O}$(1) model on the torus to the one on $\mathbb{Z}^d$ with arbitrary boundary conditions.

First of all, we employ a modification of the multi-valued mapping principle to make more clever use of Lemma \ref{juice}. The multi-valued mapping principle is a very general piece of combinatorial technology which allows for a rough sort of counting.
 
 In its essence, the principle generalises the idea that if there exists a $k$-to-$1$ map $f$ from a set $A$ to a set $B$,  then $|A|\leq k|B|$. A way of envisioning this is as a bipartite graph with the elements of $A$ and $B$ as vertices and an edge between $a$ and $b$ if $f(a)=b$. In general, if you have a bipartite graph with bipartition $(A,B)$ such that the degree of any vertex in $A$  is at least $n$ and the degree of any vertex in $B$ is at most $k$, then $n|A|\leq k|B|$. The following argument essentially generalises this fact to the setting where $A$ might also have internal edges. 
 Recall that $\mathcal{C}_{\mathtt{NT}}$ is the union of the clusters containing a wrap-around.

\begin{lemma} \label{MVMP}
For any $c>0$ and any $\varepsilon\in (0,\frac{c}{2d}),$ there exists $\delta>0$ such that the following holds:

For any $n$, and any subgraph $G$ of $\mathbb{T}^d_n$ with $cn^{d-1}$ edge-disjoint wrap-arounds $(\gamma_j)_{1\leq j\leq cn^{d-1}}$, we have 
$$
{\UEGop}_G[|\mathcal{C}_{\mathtt{NT}}|\geq \varepsilon n^{d-1}]\geq \delta.
$$
\end{lemma}
\begin{proof}
Let $G$ be a subgraph of $\mathbb{T}^d_n$ with $cn^{d-1}$ edge-disjoint wrap-arounds $(\gamma_j)_{1\leq j\leq cn^{d-1}}$. 
 For a subset $A\subset\{1, 2, \dots, cn^{d-1} \}$, let $\gamma_A:=\triangle_{a\in A}\gamma_a$.  As such, for $\eta_0\in\Even(\mathbb{T}^d_n)$, $\eta_0 \triangle \gamma_A$ is the symmetric difference of $\eta_0$ with $\gamma_j$ for each $j\in A$. Consider the following auxiliary graph $\mathfrak{G}_{\gamma,\eta_0}$ isomorphic to the Cayley graph with $(\gamma_j)$ as generators: 
     The vertices of $\mathfrak{G}_{\gamma,\eta_0}$ are $\eta_0\triangle \gamma_A,$ where $A$ ranges over the power-set of $\{1,..,cn^{d-1}\}$ and $\eta_1$ is adjacent to $\eta_2$ if $\eta_1\triangle \eta_2$ is equal to $\gamma_j$ for some $j$. See Figure \ref{fig:aux_graph}. 

Note that if $\eta_0\sim \UEG_G,$ $\eta_0\overset{d}{=} \eta_0\triangle \gamma_A$ for every $A$. Accordingly, if $\tilde{\eta}$ is a uniform vertex of $\mathfrak{G}_{\gamma,\eta_0}$, $\eta_0\overset{d}{=} \tilde{\eta}$. The upshot of this is that it suffices to argue that, deterministically, a high proportion of the vertices in $\mathfrak{G}_{\gamma,\eta_0}$ have the property that $|\mathcal{C}_{\mathtt{NT}}|\geq \varepsilon n^{d-1}$.

\begin{figure}
    \centering
    \includegraphics[scale=0.5]{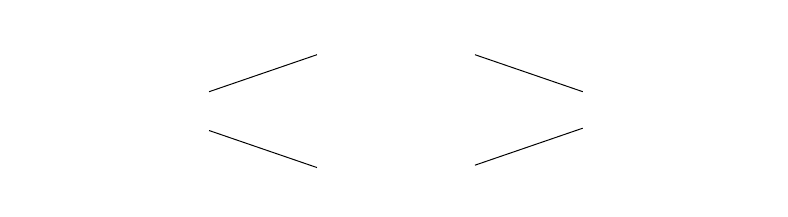}
    \caption{A small subgraph of the auxiliary graph $\mathfrak{G}_{\gamma,\eta_0}$. An initial even subgraph is exposed to the symmetric difference with all possible combinations of our initial wrap-arounds. If a vertex $\eta$ has only few wrap-arounds, then $\eta \triangle \gamma_j$ will have at least one additional wrap-around for most choices of $j$. This implies that the set of $\eta$ with few wrap-arounds must have a size which is at most comparable to the size of its complement.
    }
    \label{fig:aux_graph}
\end{figure}
Observe  that for any vertex $\eta\in \mathfrak{G}_{\gamma,\eta_0},$ and any $j$ such that $\gamma_j$ does not intersect $ E(\mathcal{C}_{\mathtt{NT}}(\eta)),$
we can apply Lemma \ref{juice} to the subgraph $(V(\mathbb{T}_n),\eta\setminus E(\mathcal{C}_{\mathtt{NT}}(\eta)))$ to get that there is at least one wrap-around in $\eta\triangle \gamma_j$ which uses no edges from $\mathcal{C}_{\mathtt{NT}}(\eta)$ while $E(\mathcal{C}_{\mathtt{NT}}(\eta))\subset \eta\triangle\gamma_j$.
So we get that
\begin{equation}\label{eq:disjoint wrap-around superedge}
    |E(\mathcal{C}_{\mathtt{NT}}(\eta\triangle \gamma_j))|\geq  |E(\mathcal{C}_{\mathtt{NT}}(\eta))|+n,
\end{equation}
for at least $cn^{d-1}-|E(\mathcal{C}_{\mathtt{NT}}(\eta))|$ values of $j$. Consequently, $|E(\mathcal{C}_{\mathtt{NT}}(\eta\triangle \gamma_j))|\leq  |E(\mathcal{C}_{\mathtt{NT}}(\eta))|$ only if $\gamma_j$ intersects $E(\mathcal{C}_{\mathtt{NT}}(\eta))$, and there are at most $|E(\mathcal{C}_{\mathtt{NT}}(\eta))|$ many such values of $j$ where this happens.

Let $\mathcal{S}_{\varepsilon}$ denote the set of vertices $\eta$ in $\mathfrak{G}_{\gamma,\eta_0}$ such that $|E(\mathcal{C}_{\mathtt{NT}}(\eta))|< d\varepsilon n^{d-1}$, 
and let $\alpha$ denote the number of edges satisfying \eqref{eq:disjoint wrap-around superedge}, that is,
\begin{equation*}
    \alpha = |\{\eta\in\mathcal{S}_{\varepsilon},\,j\le cn^{d-1}|\; |E(\mathcal{C}_{\mathtt{NT}}(\eta\triangle \gamma_j))|\geq  |E(\mathcal{C}_{\mathtt{NT}}(\eta))|+n\}|
\end{equation*}
Define $\alpha_{\ext}$ to be the number of edges with one end-point in $\mathcal{S}_{\varepsilon}$ and the other in $\mathfrak{G}_{\gamma, \eta_0}\setminus\mathcal{S}_{\varepsilon}$ and let $\alpha_{\mathrm{int}}$ denote the number of such edges with both end-points in $\mathcal{S}_{\varepsilon}$. 
The above arguments show that $\alpha\ge(c-d\varepsilon)n^{d-1}|\mathcal{S}_{\varepsilon}|$. Furthermore, for any edge $(\eta_+,\eta_-)$ contributing to $ \alpha_{\int1}$ oriented such that $|E(\mathcal{C}_{\mathtt{NT}}(\eta_+))|\geq |E(\mathcal{C}_{\mathtt{NT}}(\eta_-))|+n,$ we have that $\eta_-$ is one of the at most $d\varepsilon n^{d-1}$ neighbours of $\eta_+$ with  $|E(\mathcal{C}_{\mathtt{NT}}(\eta_-))|\le|E(\mathcal{C}_{\mathtt{NT}}(\eta_+))|$. 
Accordingly, 
$$
\alpha_{\int1}= \sum_{\eta_+\in \mathcal{S}_{\varepsilon}}\sum_{(\eta_+,\eta_-)}1\leq d\varepsilon n^{d-1}| \mathcal{S}_{\varepsilon}|.
$$
Since every vertex has $cn^{d-1}$ neighbouring edges,
$$
\alpha_{\ext}\leq cn^{d-1}|\mathfrak{G}_{\gamma, \eta_0}\setminus \mathcal{S}_{\varepsilon}|. 
$$ 

Adding this up, we get that
\begin{align*}
(c-d\varepsilon)n^{d-1}|\mathcal{S}_{\varepsilon}| &\leq \alpha= \alpha_{\ext}+\alpha_{\int1} \leq cn^{d-1}|\mathfrak{G}_{\gamma, \eta_0}\setminus \mathcal{S}_{\varepsilon}|+d\varepsilon n^{d-1}|\mathcal{S}_{\varepsilon}|,
\end{align*}
implying that
$$
\frac{|\mathfrak{G}_{\gamma, \eta_0}\setminus \mathcal{S}_{\varepsilon}|}{|\mathcal{S}_{\varepsilon}|}\geq \frac{c-2d\varepsilon}{c},
$$
or, 
equivalently, that 
$$
\frac{|\mathfrak{G}_{\gamma, \eta_0}\setminus \mathcal{S}_{\varepsilon}|}{|\mathfrak{G}_{\gamma, \eta_0}|}\geq 1-\frac{1}{2}\left(\frac{c}{c-d\varepsilon}\right)=:\delta.
$$
Now, for $\eta\in \mathfrak{G}_{\gamma, \eta_0}\setminus \mathcal{S}_{\varepsilon},$ we see that
$$
2d\varepsilon n^{d-1}\leq 2|E(\mathcal{C}_{\mathtt{NT}}(\eta))|=\sum_{v\in V(\mathcal{C}_{\mathtt{NT}}(\eta))} \deg_{\eta}(v)\leq 2d |V(\mathcal{C}_{\mathtt{NT}}(\eta))|,
$$
where $\deg_{\eta}$ denotes the degree in $\eta$. This implies that $|V(\mathcal{C}_{\mathtt{NT}}(\eta))|\geq \varepsilon n^{d-1}$.
To finish, observe that since $\tilde \eta$ was a uniform vertex of $\mathfrak{G}_{\gamma, \eta_0}$,
$$
\UEG_G[|\mathcal{C}_{\mathtt{NT}}|\geq \varepsilon n^{d-1}]\geq P[\tilde{\eta}\in \mathfrak{G}_{\gamma, \eta_0}\setminus \mathcal{S}_{\varepsilon}]\geq \delta.
$$
\end{proof}

We can now conclude the polynomial lower bound for the loop $\mathrm{O}$(1) model on the torus.
\lowerbound*
\begin{proof}
By translation invariance,
\begin{align*}
\ell_{x,\mathbb{T}_n^d}[0\cc \partial\Lambda_n]&\geq \ell_{x,\mathbb{T}_n^d}[0\in \mathcal{C}_{\mathtt{NT}}]=\ell_{x,\mathbb{T}_n^d}\left[\frac{|\mathcal{C}_{\mathtt{NT}}|}{|\mathbb{T}_n^d|}\right].
\end{align*}

Let $\boldsymbol{N}$ denote the maximal number of disjoint wrap-arounds on the torus. For any given  $c>0$, choose $\varepsilon\in (0,c/2d)$ and let $\delta$ be chosen as in as in Lemma \ref{MVMP}. Then, by Theorem \ref{thm:couplings},
\begin{align*}
\ell_{x,\mathbb{T}_n^d}\left[\frac{|\mathcal{C}_{\mathtt{NT}}|}{|\mathbb{T}_n^d|}\right]
  &\geq\phi_{x,\mathbb{T}_n^d}\left[ \id_{\boldsymbol{N}\geq cn^{d-1}}(\omega)\UEG_{\omega}\left[\frac{|\mathcal{C}_{\mathtt{NT}}|}{|\mathbb{T}^d_n|}\right]\right]\\
  &\geq \phi_{x,\mathbb{T}_n^d}\left[\id_{\boldsymbol{N}\geq cn^{d-1}}(\omega)\frac{\varepsilon n^{d-1}}{2^dn^d} \delta\right]\\
  &=\phi_{x,\mathbb{T}_n^d}[\boldsymbol{N}\geq cn^{d-1}]\cdot\frac{\varepsilon\delta}{2^d} \cdot \frac{1}{n},
\end{align*}
where, in the second inequality, we used Lemma \ref{MVMP}.
By Lemma \ref{Many wrap-arounds}, we have that $\phi_{x,\mathbb{T}_n^d}[\boldsymbol{N}\geq cn^{d-1}]\geq 1-\exp(-c'n^{d-1})$ for some value of $c'>0$, which finishes the proof.
\end{proof}
Combining the insensitivity to external topology we proved for the loop $\mathrm{O}$(1) model in \Cref{thm:Loop O(1) mixing}, this lower bound readily transfers from the torus to Euclidean space and we obtain \Cref{main theorem}.

The diameter of a cluster  $\mathcal{C}$ in $\mathbb{Z}^d$ is given by 
$
\rad( \mathcal{C}) = \sup_{x,y \in \mathcal{C}} \abs{x-y}. 
$
Recall that $\mathcal{C}_0$ is the cluster of $0$.

\begin{corollary}
For $d\geq 2$ and $x\in (x_c,1)$ then 
$$
\ell_{x,\mathbb{Z}^d}[\rad(\mathcal{C}_0)]=\infty.
$$
\end{corollary}
\begin{proof}
We have that
$$
\ell_{x,\mathbb{Z}^d}[\rad(\mathcal{C}_0)]\geq\sum_{k=1}^{\infty} \ell_{x,\mathbb{Z}^d}[\rad(\mathcal{C}_0)\geq k]= \sum_{k=1}^{\infty}\ell_{x,\mathbb{Z}^d}[0\cc \partial \Lambda_k],
$$
and the right-hand side diverges by Theorem \ref{main theorem}.
\end{proof}

\section{The loop \texorpdfstring{$\mathrm{O}$(1)}{O(1)} model on the hexagonal lattice \texorpdfstring{$\mathbb{H}$}{H} and other bi-periodic planar graphs}
\label{sec:big planar}

In this section, we first prove a no-go theorem for percolation of the loop $\mathrm{O}$(1) model on a class of planar graphs. Then, we focus on the hexagonal lattice and discuss how to adapt the arguments of \Cref{sec:Main arguments} apply to it, despite the fact (as we shall see) that the model does not percolate. This confirms that the arguments of \Cref{sec:Main arguments} alone are not strong enough to prove percolation, e.g. on $\mathbb{Z}^d$.

\subsection{Characterisation of percolation in the planar case.} \label{sec:planar} 
Using the more complete theory of planar percolation models, we can answer the question of whether $\beta_c^{\clust}(\ell)=\beta_c^{\perc}(\ell)$ for planar graphs with suitable symmetries with the aid of the arguments of \cite{GMM18} -  up to generalising some well-established results in the literature. Unfortunately, the answer varies with the graph.

First off, recall the notion of planar duality. 
To an embedded planar graph $\mathbb{G}=(\mathbb{V},\mathbb{E})$, we associate its dual graph $\mathbb{G}^*=(\mathbb{V}^*,\mathbb{E}^*),$ where $\mathbb{V}^*$ is the set of faces of $\mathbb{G}$ and every edge $e\in \mathbb{E}$ is associated to a dual edge $e^*\in \mathbb{E}^*$ between the two faces adjacent to $e$. One can check that $\mathbb{G}^*$ can be embedded into the plane by identifying a given face with a prescribed point in its interior.

 For spin models on $\{-1,+1\}^{\mathbb{V}}$ we produce the dual model on $\{-1,+1\}^{\mathbb{V}}$ by flipping spins and considering percolation in terms of *-connections instead of ordinary connections. For percolation models on $\{0,1\}^{\mathbb{E}}$ the dual model on $\{0,1\}^{\mathbb{E}^*}$ is obtained by setting $\omega^*(e^*)=1-\omega(e)$.

 In this picture, the dual model of $\phi_{p,\mathbb{G}}$ is $\phi_{p^*,\mathbb{G}^*},$ where $p$ and $p^*$ satisfy the duality relation (cf. \cite[Proposition 2.17]{DC17})
\begin{align}\label{eq:dual_parameter}
 \frac{pp^*}{(1-p)(1-p^*)}=2.
\end{align}

 The spirit of the relationship between the loop $\mathrm{O}(1)$ model and the planar Ising model goes back to Kramers and Wannier  \cite{kramers1941statistics}. Below, $\beta^*$ is the $\beta$-parameter obtained from \Cref{table:parametrizations} by plugging in $p^*$:

 \begin{proposition} \label{prop:ES}
Suppose that $\mathbb{G}$ is a planar graph. Then, for every $\beta,$ $\ell_{x,\mathbb{G}}$ is the law of the interfaces of an Ising model with $+$ boundary conditions on $\mathbb{G}^*$ at inverse temperature $\beta^*$. 
\end{proposition}

More precisely, in this picture, the loop $\mathrm{O}$(1) model on planar $\mathbb{G}$ can be coupled with the Ising model on $\mathbb{G}^*$ as the pair $(\eta,\sigma),$ where, for a given dual edge $e^*=(f,g),$ we set $\eta(e)=1$ if and only if $\sigma_f\sigma_g=-1$. This extends to $\beta^*=0,$ where $\eta\sim \UEGop_{\mathbb{G}}$ and $\sigma$ assigns $+$ and $-$ spins independently with probability $\frac{1}{2}$.

Now, we turn our attention to the class of bi-periodic embedded planar graphs, for which we have the following non-co-existence result:

\begin{theorem}[Non-co-existence] \label{thm:non_co_existence}
 For any bi-periodic planar graph $\mathbb{G} = (\mathbb{V},\mathbb{E}),$ there exists no translation-invariant measure $\mu$ on either $\{-1,+1\}^{\mathbb{V}}$ or $\{0,1\}^{\mathbb{E}}$ such that
  \begin{enumerate}
      \item[a)] $\mu$ has $\FKG$.
      \item[b)] $\mu$ has a unique infinite component and a unique infinite dual component almost surely.
  \end{enumerate}
\end{theorem}

This theorem was originally proven in \cite[Theorem 9.3.1]{She05}, and we refer the reader to \cite[Theorem 1.5]{DCsharpness} for a short, independent proof. Since any such $\mathbb{G}$ must be amenable, the Burton-Keane argument \cite{BurtonKeane} applies to all models below to show that if there is an infinite (primal or dual) component, it must be unique.
As a consequence, there is a robust subclass of bi-periodic planar graphs for which percolation of the loop $\mathrm{O}$(1) model is impossible.

\begin{proposition} \label{prop:trivalent_loop}
The loop $\mathrm{O}$(1) model on any trivalent bi-periodic planar graph $\mathbb{G}$ does not percolate, i.e. for every vertex $v \in \mathbb{V} $ and all $x \in \lbrack 0,1 \rbrack,$ then 
$$
\ell_{x, \mathbb{G}}[v \cc \infty]  = 0. 
$$
\end{proposition}
\begin{proof}
Percolation of the loop $\mathrm{O}$(1) model implies the existence of an infinite even cluster of open edges which, on a trivalent graph $\mathbb{G}$, is an infinite simple path. By \Cref{prop:ES}, this would imply the co-existence of infinite components of $+$ and $-$ for the Ising model on $\mathbb{G}^*$. Now, $\mathbb{G}^*$ is a bi-periodic planar graph, so by \Cref{thm:non_co_existence} there can be no co-existence of such infinite components in the Ising model (at finite or infinite temperature).
\end{proof}

On a bi-periodic planar graph $\mathbb{G}$ with vertices of higher degree, we would expect exponential decay of the size of $+$-clusters on $\mathbb{G}^*$, similarly to \cite{higuchi1993coexistence}. Heuristically, the subcritical Ising model on $\mathbb{G}^*$ should behave roughly like Bernoulli site percolation at parameter $\frac{1}{2}$ due to the exponential decay of correlations. In general, $*$-connections are easier than ordinary connections, so by essential enhancement (see \cite{AGenhancement,balister2014essential}), the critical parameter for site percolation, henceforth $p_c^{\mathrm{site}},$ on $\mathbb{G}$ should be strictly greater than $\frac{1}{2}$.\footnote{This is not always true, e.g. if $\mathbb{G}^*$ is obtained by periodically attaching finite graphs to the vertices of a triangulation. But $p_c^{\mathrm{site}}(\mathbb{G}^*)>\frac{1}{2}$ should be the generic case.}  Accordingly, in this range of parameters, the loop $\mathrm{O}$(1) model on $\mathbb{G}$ should percolate by the same arguments as those given in \cite{GMM18} (see \Cref{thm:two_dimensions}). On the other hand, for $x<x_c$, we  have that $\eta$ is dominated by a subcritical random-cluster model on $\mathbb{G}$ and therefore, does not percolate.

We believe that following the program outlined above, one would obtain a proof of the following conjecture:
\begin{conjecture}
Suppose that  $\mathbb{G}$ is a bi-periodic planar graph such that $p_c^{\mathrm{site}}(\mathbb{G}^*)>\frac{1}{2}$. Then,
$$
\beta_c^{\exp}(\ell_{\mathbb{G}})=\beta_c^{\perc}(\ell_{\mathbb{G}})=\beta_c(\phi_{\mathbb{G}}).
$$
\end{conjecture}

\subsection{Phase transition on the hexagonal lattice \texorpdfstring{$\mathbb{H}$}{H}} \label{sec:hexagonal} 
In order to illustrate the robustness of our arguments on the torus, we are going to apply them to the hexagonal lattice, where we know percolation of the loop $\mathrm{O}$(1) model is impossible due to \Cref{prop:trivalent_loop}.

For normalisation purposes, we embed the triangular lattice as the lattice generated by the edges $1,e^{i\frac{\pi}{3}}$, and $e^{i\frac{2\pi}{3}}$ and consider the hexagonal lattice as its dual.
Since the triangular lattice is invariant under translation by $1$ and $e^{i\frac{\pi}{3}},$ then so is the hexagonal lattice. Hence, consider the linear map $T$ which fixes $1$ and maps $i$ to $e^{i\frac{\pi}{3}}$ and the tilted box $\Lambda^{\mathbb{H}}_k=T\Lambda_k$. On this box, we can observe the quotient where $v\sim w$ if and only if $v-w\in 2k\left(\mathbb{Z} + e^{i\frac{\pi}{3}}\mathbb{Z}\right)$, which corresponds to a graph embedded on the torus on which the random-cluster model is automorphism-invariant. Since the hexagonal lattice also has a reflection symmetry in the line $\{\Im z=0\},$ we see that the toric graph thus defined is vertex-transitive and thus, $\ell_{x,\Lambda^{\mathbb{H}}_k/\sim}[0 \in \mathcal{C}_{\mathtt{NT}}] = \ell_{x,\Lambda^{\mathbb{H}}_k/\sim}[v \in \mathcal{C}_{\mathtt{NT}}]$ for any $v \in \Lambda^{\mathbb{H}}_k$.

As such, all of our arguments from above carry through so long as we can justify a version of Lemma \ref{local uniqueness} and Lemma \ref{Many wrap-arounds} for the hexagonal lattice. This, however, is downstream from sharpness of the random-cluster phase transition (see \cite{DCsharpness}) on the triangular lattice, as we shall now sketch. All the arguments are rather standard and slightly orthogonal to the themes of this paper. As such, we shall assume rough familiarity with them and only provide cursory details. We direct the interested reader to \cite{DC17} for a presentation of the necessary duality arguments.

Indeed, to get a version of Lemma \ref{local uniqueness}, note that, for $p>p_c$, the probability of having a dual crossing from $\Lambda^{\mathbb{H}}_n$ to $\partial \Lambda^{\mathbb{H}}_{2n}$ on $\mathbb{H}$ is exponentially unlikely. However, the existence of two disjoint clusters crossing the tilted annulus would imply the existence of such a dual crossing.

Similarly, the probability of having a dual top-to-bottom crossing in $\Lambda^{\mathbb{H}}_k$ is exponentially unlikely, implying that the probability of having a left-to-right crossing in $\Lambda^{\mathbb{H}}_k$ is exponentially close to $1$, which allows us to apply Lemma \ref{Exclusion tolerance} and \ref{Robustening} to get a version of Lemma \ref{Many wrap-arounds}.

Once these are off the ground, the rest of the arguments of our paper carry through without issue, yielding $\beta_c^{\clust}(\ell_{\mathbb{H}}) =\beta_c^{\exp}(\ell_{\mathbb{H}}) = \beta_c(\phi_{\mathbb{H}})$ for the loop $\mathrm{O}$(1) model on $\mathbb{H}$. 
Combined with Proposition \ref{prop:trivalent_loop}, we obtain the full phase diagram that we summarise below. We have also plotted the situation in Figure \ref{fig:Phase diagram}.
\begin{proposition} \label{prop:hexagonal_trivial}
    For the loop $\mathrm{O}$(1) model on $\mathbb{H},$
    $$\beta_c^{\clust}(\ell_{\mathbb{H}}) =\beta_c^{\exp}(\ell_{\mathbb{H}}) = \beta_c(\phi_{\mathbb{H}}),$$
    while
     $$\beta_c^{\perc}(\ell_{\mathbb{H}}) = \infty.$$
\end{proposition}
\begin{remark}  This result is largely to be found in the literature. The phase transition for the Ising and therefore also the random-cluster model on the hexagonal lattice is known to be  $\beta_c(\phi_{\mathbb{H}}) = \frac{1}{2} \arcsinh( \sqrt{3}) = \arctanh(3^{-1/2})$, see \cite[(6.5.7)]{baxter2016exactly}. In \cite{duminil2021macroscopic} a dichotomy between macroscopic loops and exponential decay is proven for the loop $\mathrm{O}(n)$ model and it is stated that the critical value for $n=1$ is known as $x = \frac{1}{\sqrt{3}}$. The critical value for general $n$ was conjectured by Nienhuis \cite{nienhuis1982exact}. 
\end{remark}
Thus, for the loop $\mathrm{O}$(1) model, we know exactly what happens on the hexagonal lattice $\mathbb{H}$. However, for the random current measure, the story is quite different. Due to stochastic domination from below by Bernoulli percolation (see \Cref{thm:couplings}), there is a percolative phase transition at some $\beta_c^{\perc}(\Prbcur_{\mathbb{H}}) < \infty$. 

Indeed, all connected even subgraphs of a trivalent graph are simple cycles, which constrains $\ell_{x,\mathbb H}^0$ rather heavily, but does not affect $\Prbcur_{x,\mathbb{H}}$ a priori. As such, we suspect that while the behaviour of the loop $\mathrm{O}(1)$ model is sensitive to the graph, the behaviour of the random current measure is generic.
\begin{conjecture}
All phase transitions of the random current and random-cluster model on the hexagonal lattice coincide: 
$$
    \beta_c^{\perc}(\Prbcur_{\mathbb{H}}) = \beta_c( \phi_{\mathbb{H}}). 
$$
\end{conjecture}

As a final aside, we add a note on odd percolation on the hexagonal lattice. Recall from \Cref{sec:odd percolation} that the uniform odd subgraph can be obtained by fixing a deterministic dimerisation and then taking the symmetric difference with a uniform even subgraph. On an odd graph, however, the uniform odd graph also arises as the complement of a uniform even graph.

\begin{proposition}\label{prop:hexagon_odd_no_perco}
For any trivalent, bi-periodic, planar graph $\mathbb{G},$ the uniform odd subgraph of $\mathbb{G}$ does not percolate. In particular, the uniform odd subgraph of $\mathbb{H}$ does not percolate.
\end{proposition}
\begin{proof}
Since the uniform odd graph of $\mathbb{G}$ is the complement of a uniform even one, we have a coupling $(\eta,\sigma)$ of a uniform odd subgraph $\eta$ of $\mathbb{G}$ and $\mathbf{I}_{0,\mathbb{G}^*}$ by setting $\eta(e)=1$ for $e^*=(f,g)$ if and only if $\sigma_f\sigma_g=1$. 
Since the spin clusters are surrounded by simple loops in the complement of $\eta$ and since  $\mathbb{G}$ is trivalent, $\eta$ does not cross these simple loops. Thus we see that the spin clusters of $\sigma$ percolates if $\eta$ percolates. 
However, by spin-symmetry, this would imply the co-existence of an infinite cluster of $+$'s and an infinite cluster of $-$'s. Since this does not happen, we conclude that $\eta$ does not percolate.
\end{proof}

\section{Perspective}
We conclude by discussing some more general properties percolation of the UEG and some open problems. 

\subsection{General characterisations of percolation of the UEG}

As with any model without positive association, getting a grip of the general behaviour of the UEG is a priori a daunting task, but one might hope for some semi-robust arguments that allow one to handle large classes of graphs, as for the planar case or the case where the graph $G$ in question has a subgraph $H$ where $\Even$ separates edges and $p_c(\mathbb{P}_{p,H})<\frac{1}{2}$.

We saw in \Cref{sec:planar} that trivalence prevents the $\UEGop$ on bi-periodic planar graphs from percolating. However, as we now show, it is not true that trivalence is an obstruction to percolation in general. In light of the results in \cite{angel2021uniform} and \Cref{prop:evens and ends} , this question has a significantly different flavour depending on whether or not the graph in question is one-ended. As one anonymous referee pointed out, and as can be seen in \cite{monster_paper} in higher generality, the (wired) \UEG\; can be seen to percolate on a 3-regular tree, which is planar but has uncountably many ends. However, on one-ended graphs, the free and wired UEG coincide \cite[Lemma 3.9]{angel2021uniform}. Below, we give an example which is one-ended, but non-amenable and non-planar.

Consider the trivalent supergraph $\mathbb{J}$ of $\mathbb{N},$ adding the edges $(1,10),$ $(1,100)$ and $(n,10^{n+1})$ for all $n$ that are not powers of $10$. Denote by $\mathbb{E}(\mathbb{J})$ and $\mathbb{E}(\mathbb{N})$ the corresponding edge sets.

\begin{proposition} \label{prop:jump_graph} 
 $\UEGop_{\mathbb{J}}=\UEGop^0_{\mathbb{J}}=\UEGop^1_{\mathbb{J}}$ percolates. In particular, there exists an infinite trivalent one-ended graph for which the uniform even subgraph percolates.
\end{proposition}

\begin{proof}
We follow the construction of the uniform even subgraph from \cite{angel2021uniform} and construct a basis for $\Omega_{\emptyset}(\mathbb{J})$ from a spanning tree.
Pick the spanning tree in $\mathbb{J}$ which is simply the original graph $\mathbb{N}$ and for every edge $e\in \mathbb{E}(\mathbb{J})\setminus \mathbb{E}(\mathbb{N})$ (henceforth called 'external'), let $C_e$ be the simple loop in $\mathbb{J}$ given by joining $e$ with the unique path in $\mathbb{N}$ connecting the end-points of $e$. One may check that the $C_e$ form a finitary basis of $\Omega_{\emptyset}(\mathbb{J})$. One may note, in particular, that every external edge belongs to a unique such loop.

Thus, similar to \Cref{eq:sample-cycles}, for every external edge $e,$ we can let $\epsilon_e$ be an i.i.d. family of Bernoulli-$\frac{1}{2}$ variables and sample the UEG of $\mathbb{J}$ as 
$$
\sum_e \epsilon_e C_e.
$$

Due to the trivalence of $\mathbb{J},$ this gives the cluster of $1$ the following random walk type representation on $\mathbb{N}$: First, we set $x_0:=1$ and reveal the states of $\epsilon_{(1,10)}$ and $\epsilon_{(1,100)}.$ If both are $0$, then the cluster of $1$ is trivial and the process ends. Otherwise, we set $x_1$ equal to the largest number $j$ such that $\epsilon_{(1,j)}=1$ and set $e_1=(1,j)$.

Now, recursively, given that the process arrived at $x_j$ through the edge $e_j$, since the cluster of $1$ is a simple loop in $\mathbb{J},$ exactly one of the other edges adjacent to $x_j$ must be open. We proceed by cases:

\begin{enumerate}
    \item[$i)$] We have already revealed the state of one of the neighboring edges. In this case, we know the unique open edge among the two and set $e_{j+1}$ equal to this edge and $x_{j+1}$ equal to the other end-point of $e_{j+1}$. Either $x_{j+1}=1,$ and the process terminates, or it is not, in which case the process continues.
    \item[$ii)$] $e_j$ is an external edge and we have not revealed the state of $(x_j-1,x_j)$. In this case, we check the sum of $\epsilon_e$ for every external edge $e=(n,m)$ such that $n<x_j<m$ (these are exactly the external edges such that $(x_j-1,x_j)\in C_e$). If the sum is even, then $(x_j-1,x_j)$ is closed, so $(x_j,x_j+1)$ must be open and vice versa. Since we have not yet revealed the state of $(x_j-1,x_j),$ the conditional parity of the sum is a Bernoulli random variable. In conclusion, with probability $\frac{1}{2},$ $(x_j-1,x_j)$ is open, and we set $e_{j+1}=(x_j-1,x_j)$ and $x_{j+1}=x_{j}-1$. Otherwise, we set $e_{j+1}=(x_j,x_j+1)$ and $x_{j+1}=x_j+1$. The upshot is that the process goes left with probability $\frac{1}{2}$ and right otherwise. 
    \item[$iii)$] $e_j$ is not an external edge and the state of the unique external edge $e^{x_j}$ going through $x_j$ has not yet been revealed. In this case, if $\epsilon_{e^{x_j}}=1,$ we set $e_{j+1}=e^{x_j}$ and $x_{j+1}$ equal to the other end-point of $e^{x_j}$. Otherwise, the edge $e_{j+1}:=(x_j,x_j+(x_j-x_{j-1}))$ is open, and we set $x_{j+1}=x_j+(x_j-x_{j-1})$. The upshot is that the process takes the external edge it just arrived at with probability $\frac{1}{2}$, and otherwise it continues along $\mathbb{N}$ in the same direction as it has been traveling thus far.
\end{enumerate}
The upshot is that the walk arrives somewhere via an external edge, turns left or right with probability $\frac{1}{2}$ and keeps walking for a geometric number of steps until it encounters an open external edge. The probability that $1$ is connected to $\infty$ is then at least the probability that $x_j$ never has the form $10^n$ at times $j$ where $e_j$ is not an external edge (note that $1,$ $10$ and $100$ do have this form). This is achieved if $B_n$ happens for every $n$, where $B_n$ denotes the event that there is an open external edge in $[10^n+1,10^{n+1}-1]$. 

 Since the marginal of the open external edges is Bernoulli $\frac{1}{2}$ percolation, we get that
 $$
 \UEG_{\mathbb{J}}[B_n]\geq 1-2^{10^{n+1}-10^{n}-2}\geq 1-2^{8\cdot 10^n}
 $$
 A union bound now yields that
$$
\UEG_{\mathbb{J}}[1\cc \infty] \geq \UEG_{\mathbb{J}}[\cap_n B_n]\geq 1-\sum_{n=0}^{\infty}2^{-8\cdot 10^{n}}>0.
$$
\end{proof}

\subsection{The situation for $p=p_c$}
Our proof of infinite expectation of cluster sizes of $\ell_{p, \Z^d}$ required that $p > p_c$. For $p<p_c$ there is exponential decay by stochastic domination by $\phi_p$. Here, we briefly discuss the situation for $p=p_c.$

For $d\geq 3$ and $p=p_c,$ we do have a polynomial lower bound for connection probabilities, since by \cite[Theorem 4.8]{DC17}, there exists $c,C > 0$ such that
$$
\frac{c}{\abs{v}^{d-1}} \leq \langle \sigma_0 \sigma_v \rangle_{\beta_c, \Z^d} \leq \frac{C}{\abs{v}^{d-2}},
$$
where we note that a more general, and for $d \geq 5$ tighter, bound is proven in \cite{aizenman1988critical}. 
Therefore, $\phi_{p_c,\Z^d}[\abs{\mathcal{C}_0}] = \infty$. On the other hand, it does not percolate by continuity of the Ising phase transition \cite{DCcont}. Heuristically, a model with infinite expected cluster sizes should percolate as soon as any independent density of edges is added to it. Therefore, with the coupling from \Cref{thm:couplings} in mind, we do not expect that $\ell_{x_c, \Z^d}[\abs{\mathcal{C}_0}] = \infty$. 

However, for any $d\geq 2,$ we still believe that connection probabilities of $\ell_{x_c,\mathbb{T}_n^d}$ satisfy polynomial lower bounds in the volume of the torus. It remains, however, a difficulty to transfer the result to $\Z^d$, since we cannot use Pizstora's construction to exhibit separating surfaces in $\phi_{p_c,\Z^d}$. For $d=2$, polynomial bounds on the existence can be achieved with ordinary RSW theory, which suffices for establishing a polynomial lower bound on cluster sizes in $\ell_{x_c, \Z^2}$. In other dimensions, though, no similar tool exists in the literature to our knowledge. We summarise our expectations in the following conjecture:

\begin{conjecture}
  Let $d \geq 3$. Then we expect that $\ell_{x_c, \Z^d}[ \abs{\mathcal{C}_0}] < \infty$ and that there exists some $a,b > 0$ such that 
  $$ \frac{a}{\abs{v}^b} < \ell_{x_c, \Z^d}[0 \cc v].$$
\end{conjecture}

 \subsection{Remaining questions for the phase diagram of $\ell$ and $\Prbcur$ on $\Z^d$}
In \Cref{main theorem} and \Cref{prop:hexagonal_trivial} we managed to prove a condition akin to criticality all the way down to the random-cluster phase transition for the loop $\mathrm{O}(1)$ model on both $\Z^d$ and $\mathbb{H}$. 
However, on $\Z^2,$ the transition from exponential decay to percolation happens at one point, whereas it never happens for $\mathbb{H}$. This motivates the following question that asks whether a proper intermediate regime can exist: 

\begin{question}
Does there exist a lattice $\mathbb{L}$ such that the loop $\mathrm{O}$(1) model $\ell_{x,\mathbb{L}}$ for $x \in \lbrack 0,1 \rbrack$ has a non-trivial intermediate regime, i.e. an interval $(a,b) \subset \lbrack 0,1 \rbrack$ and points $0 < x_0 < a$ and $b < x_1 < 1$ such that $\ell_{x_0,\mathbb{L}}$ has exponential decay, $\ell_{x_1,\mathbb{L}}$ percolates, and $\ell_{x,\mathbb{L}}$ neither has exponential decay nor percolates for $x \in (a,b)$?
\end{question}

Even with our main Theorem in mind, our motivating problem of interest 
\cite[Question 1]{DC16} is still left open for $\Z^d$ with $d\geq 3$. Since the random current model stochastically dominates the loop $\mathrm{O}$(1) model, a positive answer would follow from sharpness for $\ell_x$. One reason to suspect this on $\mathbb{Z}^d$ goes as follows:

The estimate in \Cref{thm:Torus wrap} is particularly crude. We essentially only use that one may increase the size of the non-trivial clusters by taking the symmetric difference with a wrap-around which intersects only trivial clusters. However, due to the existence of vertices of degree at least 4, there are plenty of scenarios where acting by a wrap-around increases the size of the non-trivial cluster even when the two intersect. As such, we conjecture the following, which was noted for $d=2$ in \cite{GMM18} (cf. \Cref{thm:two_dimensions}):

\begin{conjecture} \label{conj:overall} 
For $\Z^d$, $d \geq 3$ it holds that $\beta_c^{\perc}(\ell) = \beta_c^{\exp}(\ell)$.  
\end{conjecture}
Of course, one may also settle for the weaker statement:
 \begin{conjecture}
For $d\geq 3$, the single random current on $\Z^d$ has a unique sharp percolative phase transition, i.e. $\beta_c^{\perc}(\Prbcur) = \beta_c^{\exp}(\Prbcur)$. 
 \end{conjecture}

 \subsection{Infinite clusters of the loop $\mathrm{O}$(1) model on the cut open lattice}\label{sec:zipper}

 The technical reason that we cannot extend our result from infinite expectation of cluster sizes to percolation is that we have difficulties in controlling the variations of where wrap-arounds occur. If $\eta$ is trivial and $\gamma$ is a wrap-around, then $\eta\triangle \gamma$ contains at least one wrap-around by \Cref{juice}, but there is no reason to suspect that this wrap-around intersects $\gamma$ at all. In an artificial setup that we will now sketch, we can overcome this barrier.

 Consider the graph $\mathcal{Z}^d$ obtained from $\Z^d$ by removing all edges from a fixed hyperplane with the exception of a single edge, $e$. We also consider a similar cut-up $\mathcal{H}$ version of the hexagonal lattice $\mathbb{H}$. Just as in the non-cut-up case, we may consider a quotient of $\Lambda_{n}\cap \mathcal{Z}^d$ as a subgraph of the torus and carry out our arguments from before (and similarly for $\Lambda_n^{\mathbb{H}}\cap \mathcal{H}$).

 \begin{theorem}
For $p > p_c(\phi_{\Z^d})$, there exists an infinite cluster of $\ell_{p,\mathcal{Z}^d}$ with positive probability. 
Even though $\mathcal{H} \subset \mathbb{H}$, there exists a $p > 0$ such that 
$$\ell_{p,\mathcal{H}} \lbrack 0 \cc \infty \rbrack  > 0 = \ell_{p,\mathbb{H}} \lbrack 0 \cc \infty\rbrack. $$
 \end{theorem}
 \begin{proof}
Since $p > p_c(\phi_{\Z^d})$ and the critical parameter for the random-cluster model on the half space is the same as for the full space (see e.g. \cite{Bod05} for a much more general result), we also have  $p_c(\phi_{\Z^d}) = p_c(\phi_{\mathcal{Z}^d})$. Now, the infinite cluster on each half plane intersects $e$ with positive probability. Thus, by the FKG inequality (\Cref{FKG}) and insertion tolerance  \Cref{thm:comparison theorem}$iii)$ the edge $e$ is part of an infinite self-avoiding path 
configuration with positive probability. Suppose that $\omega$ is such a configuration and that $\gamma\subset \omega$ is such a path. If $\eta \sim \UEG_\omega$ then also $\eta \triangle \gamma \sim \UEG_\omega$. A parity consideration shows that there is an infinite path through $e$ in one of the two configurations. Thus, $\ell_{p,\mathcal{Z}^d}$ has an infinite cluster with at positive probability. The same reasoning holds for $\mathcal{H}$.
 \end{proof}

One reason that we show the argument here is that we speculate that finer control of the variations of the wrap-arounds arising from the combinatorial argument could help shed light on the remaining questions. 
 
\section*{Acknowledgements} 
We thank Ioan Manolescu and  Peter Wildemann for helpful discussions and Aran Raoufi and Franco Severo for getting us started on the problem in the first place.
We thank the anonymous referees for helpful comments.
BK and FRK acknowledge the Villum Foundation for funding through the QMATH center of Excellence (Grant No. 10059) and the Villum Young Investigator (Grant No. 25452) programs. UTH acknowledges funding from Swiss SNF.

\bibliographystyle{abbrv}
\bibliography{bibliography}

\begin{thebibliography}{10}

\bibitem{aizenman1982geometric}
M.~Aizenman.
\newblock {Geometric analysis of $\varphi^{4}$ fields and Ising models. Parts I
  and II}.
\newblock {\em Communications in Mathematical Physics}, 86(1):1--48, 1982.

\bibitem{aizenman1987phase}
M.~Aizenman, D.~J. Barsky, and R.~Fern{\'a}ndez.
\newblock {The phase transition in a general class of Ising-type models is
  sharp}.
\newblock {\em Journal of Statistical Physics}, 47(3):343--374, 1987.

\bibitem{AllTheAuthors}
M.~Aizenman, J.~T. Chayes, L.~Chayes, J.~Fr\"{o}hlich, and L.~Russo.
\newblock On a sharp transition from area law to perimeter law in a system of
  random surfaces.
\newblock {\em Communications in Mathematical Physics}, 92(1):19--69, 1983.

\bibitem{aizenman2021marginal}
M.~Aizenman and H.~Duminil-Copin.
\newblock {Marginal triviality of the scaling limits of critical 4D Ising and
  $\varphi_4^4$ models}.
\newblock {\em Annals of Mathematics}, 194(1):163--235, 2021.

\bibitem{aizenman2015random}
M.~Aizenman, H.~Duminil-Copin, and V.~Sidoravicius.
\newblock {Random currents and continuity of Ising model’s spontaneous
  magnetization}.
\newblock {\em Communications in Mathematical Physics}, 334(2):719--742, 2015.

\bibitem{aizenman1988critical}
M.~Aizenman and R.~Fern{\'a}ndez.
\newblock Critical exponents for long-range interactions.
\newblock {\em Letters in Mathematical Physics}, 16(1):39--49, 1988.

\bibitem{AGenhancement}
M.~Aizenman and G.~Grimmett.
\newblock Strict monotonicity for critical points in percolation and
  ferromagnetic models.
\newblock {\em Journal of Statistical Physics}, 63(5-6):817--835, 1991.

\bibitem{angel2021uniform}
O.~Angel, G.~Ray, and Y.~Spinka.
\newblock Uniform even subgraphs and graphical representations of ising as
  factors of iid.
\newblock {\em Electronic Journal of Probability}, 29:1--31, 2024.

\bibitem{balister2014essential}
P.~Balister, B.~Bollob{\'a}s, and O.~Riordan.
\newblock Essential enhancements revisited.
\newblock {\em arXiv preprint arXiv:1402.0834}, 2014.

\bibitem{baxter2016exactly}
R.~J. Baxter.
\newblock {\em Exactly solved models in statistical mechanics}.
\newblock Elsevier, 2016.

\bibitem{berezinskii1971destruction}
Berezinskii.
\newblock {Destruction of long-range order in one-dimensional and
  two-dimensional systems having a continuous symmetry group I. Classical
  systems}.
\newblock {\em Sov. Phys. JETP}, 32(3):493--500, 1971.

\bibitem{Bod05}
T.~Bodineau.
\newblock Slab percolation for the {I}sing model.
\newblock {\em Probability Theory and Related fields}, 132(1):83--118, 2005.

\bibitem{broadbent1957percolation}
S.~R. Broadbent and J.~M. Hammersley.
\newblock {Percolation processes: I. Crystals and mazes}.
\newblock In {\em Mathematical proceedings of the Cambridge philosophical
  society}, volume 53,3, pages 629--641. Cambridge University Press, 1957.

\bibitem{BurtonKeane}
R.~M. Burton and M.~Keane.
\newblock Density and uniqueness in percolation.
\newblock {\em Communications in Mathematical Physics}, 121(3):501--505, 1989.

\bibitem{camia2020exponential}
F.~Camia, J.~Jiang, and C.~M. Newman.
\newblock {Exponential Decay for the Near-Critical Scaling Limit of the Planar
  Ising Model}.
\newblock {\em Communications on Pure and Applied Mathematics},
  73(7):1371--1405, 2020.

\bibitem{campanino1985upper}
M.~Campanino and L.~Russo.
\newblock An upper bound on the critical percolation probability for the
  three-dimensional cubic lattice.
\newblock {\em The Annals of Probability}, pages 478--491, 1985.

\bibitem{crawford2020macroscopic}
N.~Crawford, A.~Glazman, M.~Harel, and R.~Peled.
\newblock {Macroscopic loops in the loop $O(n)$ model via the XOR trick}.
\newblock {\em arXiv preprint arXiv:2001.11977}, 2020.

\bibitem{DC16}
H.~Duminil-Copin.
\newblock {Random current expansion of the Ising model}.
\newblock {\em Proceedings of the 7th European Congress of Mathematicians in
  Berlin}, 2016.

\bibitem{duminil2018sixty}
H.~Duminil-Copin.
\newblock Sixty years of percolation.
\newblock In {\em Proceedings of the International Congress of Mathematicians:
  Rio de Janeiro 2018}, pages 2829--2856. World Scientific, 2018.

\bibitem{DC17}
H.~Duminil-Copin.
\newblock {Lectures on the Ising and Potts models on the hypercubic lattice.}
\newblock {\em PIMS-CRM Summer School in Probability}, 2019.

\bibitem{duminil2021macroscopic}
H.~Duminil-Copin, A.~Glazman, R.~Peled, and Y.~Spinka.
\newblock {Macroscopic loops in the loop $O(n)$ model at Nienhuis' critical
  point.}
\newblock {\em Journal of the European Mathematical Society (EMS Publishing)},
  23(1), 2021.

\bibitem{duminil2020exponential}
H.~Duminil-Copin, S.~Goswami, and A.~Raoufi.
\newblock {Exponential decay of truncated correlations for the Ising model in
  any dimension for all but the critical temperature}.
\newblock {\em Communications in Mathematical Physics}, 374(2):891--921, 2020.

\bibitem{duminil2019double}
H.~Duminil-Copin and M.~Lis.
\newblock On the double random current nesting field.
\newblock {\em Probability Theory and Related Fields}, 175(3):937--955, 2019.

\bibitem{duminil2021conformal}
H.~Duminil-Copin, M.~Lis, and W.~Qian.
\newblock {Conformal invariance of double random currents II: tightness and
  properties in the discrete}.
\newblock {\em arXiv preprint arXiv: 2107.12880}, 2021.

\bibitem{duminil2021conformal2}
H.~Duminil-Copin, M.~Lis, and W.~Qian.
\newblock Conformal invariance of double random currents {I}: {I}dentification
  of the limit.
\newblock {\em Proc. Lond. Math. Soc. (3)}, 130(1):Paper No. e70022, 2025.

\bibitem{DCsharpness}
H.~Duminil-Copin, A.~Raoufi, and V.~Tassion.
\newblock {Sharp phase transition for the random-cluster and {P}otts models via
  decision trees}.
\newblock {\em Annals of Mathematics}, 189(1):75--99, 2019.

\bibitem{DCcont}
H.~Duminil-Copin, V.~Sidoravicius, and V.~Tassion.
\newblock Continuity of the phase transition for planar random-cluster and
  {P}otts models with {$1 \leq q \leq 4$}.
\newblock {\em Communications in Mathematical Physics}, 349(1):47--107, 2017.

\bibitem{duminil2016new}
H.~Duminil-Copin and V.~Tassion.
\newblock {A new proof of the sharpness of the phase transition for Bernoulli
  percolation and the Ising model}.
\newblock {\em Communications in Mathematical Physics}, 343(2):725--745, 2016.

\bibitem{edwards1988generalization}
R.~G. Edwards and A.~D. Sokal.
\newblock {Generalization of the Fortuin-Kasteleyn-Swendsen-Wang representation
  and Monte Carlo algorithm}.
\newblock {\em Physical review D}, 38(6):2009, 1988.

\bibitem{fortuin1972random}
C.~M. Fortuin and P.~W. Kasteleyn.
\newblock On the random-cluster model: I. introduction and relation to other
  models.
\newblock {\em Physica}, 57(4):536--564, 1972.

\bibitem{friedli2017statistical}
S.~Friedli and Y.~Velenik.
\newblock Statistical mechanics of lattice systems: a concrete mathematical
  introduction.
\newblock {\em Cambridge University Press}, 2017.

\bibitem{frohlich1981kosterlitz}
J.~Fr{\"o}hlich and T.~Spencer.
\newblock {The Kosterlitz-Thouless transition in two-dimensional Abelian spin
  systems and the Coulomb gas}.
\newblock {\em Communications in Mathematical Physics}, 81(4):527--602, 1981.

\bibitem{GMM18}
O.~Garet, R.~Marchand, and I.~Marcovici.
\newblock {Does Eulerian percolation on $\mathbb{Z}^2$ percolate?}
\newblock {\em ALEA, Lat. Am. J. Probab. Math. Stat.}, page 279–294, 2018.

\bibitem{griffiths1970concavity}
R.~B. Griffiths, C.~A. Hurst, and S.~Sherman.
\newblock {Concavity of magnetization of an Ising ferromagnet in a positive
  external field}.
\newblock {\em Journal of Mathematical Physics}, 11(3):790--795, 1970.

\bibitem{GJ09}
G.~Grimmet and S.~Janson.
\newblock { Random even graphs}.
\newblock {\em The Electronic Journal of Combinatorics, Volume 16, Issue 1},
  2009.

\bibitem{Gri06}
G.~Grimmett.
\newblock {The random-cluster model}.
\newblock {\em volume 333 of { Grundlehren der Mathematischen Wissenschaften
  [Fundamental Principles of Mathematical Sciences]}.}, 2006.

\bibitem{Haggstrom1996}
O.~H{\"a}ggstr{\"o}m.
\newblock The random-cluster model on a homogeneous tree.
\newblock {\em Probability Theory and Related Fields}, 104(2):231--253, 1996.

\bibitem{Unique_Forest}
N.~Halberstam and T.~Hutchcroft.
\newblock Uniqueness of the infinite tree in low-dimensional random forests.
\newblock {\em Probab. Math. Phys.}, 5(4):1185--1216, 2024.

\bibitem{hansen2022strict}
U.~T. Hansen and F.~R. Klausen.
\newblock {Strict monotonicity, continuity, and bounds on the {K}ert\'{e}sz
  line for the random-cluster model on {$\Bbb{Z}^d$}}.
\newblock {\em Journal of Mathematical Physics}, 64(1):Paper No. 013302, 22,
  2023.

\bibitem{monster_paper}
U.~T. Hansen, F.~R. Klausen, and P.~Wildemann.
\newblock {Non-uniqueness of phase transitions for graphical representations of
  the Ising model on tree-like graphs}, 2024.
\newblock arXiv, 2410.22061.

\bibitem{higuchi1993coexistence}
Y.~Higuchi.
\newblock {Coexistence of infinite (*)-clusters II. Ising percolation in two
  dimensions}.
\newblock {\em Probability theory and related fields}, 97(1):1--33, 1993.

\bibitem{hutchcroft2020continuity}
T.~Hutchcroft.
\newblock Continuity of the ising phase transition on nonamenable groups.
\newblock {\em Communications in Mathematical Physics}, pages 1--60, 2023.

\bibitem{kesten1980critical}
H.~Kesten et~al.
\newblock The critical probability of bond percolation on the square lattice
  equals 1/2.
\newblock {\em Communications in Mathematical Physics}, 74(1):41--59, 1980.

\bibitem{kitaev2003fault}
A.~Y. Kitaev.
\newblock Fault-tolerant quantum computation by anyons.
\newblock {\em Annals of Physics}, 303(1):2--30, 2003.

\bibitem{klausen2021monotonicity}
F.~R. Klausen.
\newblock On monotonicity and couplings of random currents and the
  loop-$\mathrm{O}(1)$-model.
\newblock {\em ALEA}, 19:151--161, 2022.

\bibitem{klausen2022mass}
F.~R. Klausen and A.~Raoufi.
\newblock {Mass scaling of the near-critical 2D Ising model using random
  currents}.
\newblock {\em Journal of Statistical Physics}, 188(3):1--21, 2022.

\bibitem{kosterlitz1973ordering}
J.~M. Kosterlitz and D.~J. Thouless.
\newblock Ordering, metastability and phase transitions in two-dimensional
  systems.
\newblock {\em Journal of Physics C: Solid State Physics}, 6(7):1181, 1973.

\bibitem{kozma2013lower}
G.~Kozma and V.~Sidoravicius.
\newblock Lower bound for the escape probability in the lorentz mirror model on
  $\mathbb{Z}^2$.
\newblock {\em Israel Journal of Mathematics}, 209(2):683--685, 2015.

\bibitem{kramers1941statistics}
H.~A. Kramers and G.~H. Wannier.
\newblock {Statistics of the two-dimensional ferromagnet. Part I}.
\newblock {\em Physical Review}, 60(3):252, 1941.

\bibitem{lenz1920beitrag}
W.~Lenz.
\newblock {Beitrag zum Verst{\"a}ndnis der magnetischen Erscheinungen in festen
  K{\"o}rpern}.
\newblock {\em Z. Phys.}, 21:613--615, 1920.

\bibitem{liggett1997domination}
T.~M. Liggett, R.~H. Schonmann, and A.~M. Stacey.
\newblock Domination by product measures.
\newblock {\em The Annals of Probability}, 25(1):71--95, 1997.

\bibitem{Lis}
M.~Lis.
\newblock Spins, percolation and height functions.
\newblock {\em Electronic Journal of Probability}, 27:1--21, 2022.

\bibitem{LW16}
T.~Lupu and W.~Werner.
\newblock {{A note on Ising random currents, Ising-FK, loop-soups and the
  Gaussian free field}}.
\newblock {\em Electronic Communications in Probability}, 21:7 pp., 2016.

\bibitem{nienhuis1982exact}
B.~Nienhuis.
\newblock {Exact critical point and critical exponents of $O(n)$ models in two
  dimensions}.
\newblock {\em Physical Review Letters}, 49(15):1062, 1982.

\bibitem{onsager1944crystal}
L.~Onsager.
\newblock {Crystal statistics. I. A two-dimensional model with an
  order-disorder transition}.
\newblock {\em Physical Review}, 65(3-4):117, 1944.

\bibitem{peierls1936ising}
R.~Peierls.
\newblock {On Ising's model of ferromagnetism}.
\newblock In {\em Mathematical Proceedings of the Cambridge Philosophical
  Society}, volume 32,3, pages 477--481. Cambridge University Press, 1936.

\bibitem{Pis96}
A.~Pisztora.
\newblock Surface order large deviations for {I}sing, {P}otts and percolation
  models.
\newblock {\em Probability Theory and Related Fields}, 104(4):427--466, 1996.

\bibitem{Rao_1971}
M.~Rao.
\newblock Projective limits of probability spaces.
\newblock {\em Journal of Multivariate Analysis}, 1(1):28--57, apr 1971.

\bibitem{raoufi2020translation}
A.~Raoufi.
\newblock {Translation-invariant Gibbs states of the Ising model: general
  setting}.
\newblock {\em The Annals of Probability}, 48(2):760--777, 2020.

\bibitem{She05}
S.~R. Sheffield.
\newblock {\em Random surfaces: {L}arge deviations principles and gradient
  {G}ibbs measure classifications}.
\newblock ProQuest LLC, Ann Arbor, MI, 2003.
\newblock Thesis (Ph.D.)--Stanford University.

\bibitem{tassion_notes}
V.~Tassion.
\newblock Ising model.
\newblock https://metaphor.ethz.ch/x/2021/hs/401-3822-17L/sc/lectureNotes.pdf,
  2021.
\newblock lecture notes, ETH Zürich.

\bibitem{van1941lange}
B.~L. van~der Waerden.
\newblock {Die lange Reichweite der regelm{\"a}ssigen Atomanordnung in
  Mischkristallen}.
\newblock {\em Zeitschrift f{\"u}r Physik}, 118(7-8):473--488, 1941.

\end{thebibliography}

\end{document}